\documentclass{article}

\usepackage{arxiv}

\usepackage[utf8]{inputenc} 
\usepackage[T1]{fontenc}    
\usepackage{hyperref}       
\usepackage{url}            
\usepackage{booktabs}       
\usepackage{amsfonts}       
\usepackage{nicefrac}       
\usepackage{lipsum}
\usepackage{pdflscape}
\usepackage{morefloats}
\usepackage{color}
\usepackage{multirow}
\usepackage{latexsym}
\usepackage{amsmath,amssymb,amsthm,graphicx}
\usepackage{dsfont}
\usepackage{enumitem}
\usepackage{hyperref}

\newtheoremstyle{thm}
{9pt}
{9pt}
{\itshape}
{}
{\bfseries}
{.}
{ }
{}
\theoremstyle{thm}

\newtheorem{theorem}{Theorem}[section]

\newtheorem{corollary}[theorem]{Corollary}

\newcommand{\vertk}{\stackrel{{\cal D}}{\longrightarrow}}
\newcommand{\edist}{\stackrel{{\cal D}}{=}}
\newcommand{\stk}{\stackrel{\mbox{\scriptsize $\PP$}}{\longrightarrow}}
\newcommand{\fsk}{\stackrel{{\rm a.s.}}{\longrightarrow}}

\newtheoremstyle{def}
{9pt}
{9pt}
{}
{}
{\bfseries}
{.}
{ }
{}
\theoremstyle{def}

\newenvironment{prf}{\textbf{\emph{Proof.}}}{\qed}

\newcommand{\CSP}{{{\text{CS}}^+}} 
\newcommand{\CSM}{{{\text{CS}}^-}} 
\newcommand{\ii}{{\text{i}}} 
\newcommand{\C}{\mathbb{C}} 
\newcommand{\R}{\mathbb{R}} 
\newcommand{\N}{\mathbb{N}} 
\newcommand{\ID}{{\text{I}}_d} 
\newcommand{\E}{\mathbb{E}} 
\newcommand{\PP}{\mathbb{P}} 
\newcommand{\HH}{\mathbb{H}} 

\renewcommand{\footnoterule}{%
	\kern -3.5pt
	\hrule width \textwidth height 1pt
	\kern 3.5pt
}

\makeatletter
\def\blfootnote{\xdef\@thefnmark{}\@footnotetext}
\makeatother

\title{Testing normality in any dimension by Fourier methods in a multivariate Stein equation}

\author{Bruno Ebner\\
Institute of Stochastics, \\
Karlsruhe Institute of Technology (KIT), \\
Englerstr. 2, D-76133 Karlsruhe. \\
\texttt{Bruno.Ebner@kit.edu}\\
\And
Norbert Henze\\
Institute of Stochastics, \\
Karlsruhe Institute of Technology (KIT), \\
Englerstr. 2, D-76133 Karlsruhe. \\
\texttt{Norbert.Henze@kit.edu}\\
\And David Strieder\\
Karlsruher Str. 62a,\\
D-69126 Heidelberg, \\
\texttt{david.strieder@student.kit.edu}\\
}

\begin{document}

\date{\today}
\maketitle


\begin{abstract}
We study a novel class of affine invariant and consistent tests for multivariate normality. The tests are based  on a characterization
of the standard $d$-variate normal distribution by means of the unique solution of an initial value problem connected to a partial differential equation, which is motivated by a multivariate Stein equation. The test
criterion is a suitably weighted $L^2$-statistic. We derive the limit distribution of the test statistic under the null hypothesis as well as
under contiguous and fixed alternatives to normality. A consistent estimator of the limiting variance under fixed alternatives as well as an asymptotic confidence interval of the distance of an underlying alternative with respect to the multivariate normal law is derived. In simulation studies, we show that the tests are strong in comparison with prominent competitors, and that the empirical coverage rate of the asymptotic confidence interval converges to the nominal level.
We present a real data example, and we outline topics for further research.
\end{abstract}

\section{Introduction.}
Statistical inference for a data set starts with assumptions on the underlying stochastic mechanism which determines the generation of the data. In most classical models for multidimensional data, such as multivariate linear regression models or multivariate analysis of variance, the assumption of multivariate normality of the underlying random vectors is inherent. Hence, prior to any serious statistical inference, one should check this assumption. To be specific, let $X,X_1,X_2, \ldots $ be a sequence of independent identically distributed (i.i.d.) $d$-dimensional (column) vectors that are defined on a common probability space $(\Omega,{\cal A},\PP)$. We make the basic standing assumptions that the distribution $\PP^X$ of $X$ is absolutely continuous
with respect to $d$-dimensional Lebesgue measure. In what follows, we denote by  N$_d(\mu,\Sigma)$  the $d$-variate normal distribution with expectation
vector $\mu$ and covariance matrix $\Sigma$, and  we write
\[
\mathcal{N}_d:=\{\text{N}_d(\mu,\Sigma): \mu \in \R^d,\ \Sigma \in \R^{d\times d} \text{ positive definite}\}
\]
for the class of all non-degenerate $d$-variate normal distributions. The unit matrix of order $d$ will be denoted by $\ID$. The problem of matter is testing the hypothesis
\[
H_0: \PP^X \in  \mathcal{N}_d,
\]
based on $X_1,\ldots,X_n$, against general alternatives. The purpose of this paper is to introduce and study a novel class of affine invariant and consistent tests based on a partial differential equation (PDE) that determines the characteristic function of the multivariate standard normal law. We write $\nabla$ for the gradient operator and consider for $f \in L^2(\R^d)$  the initial value problem of the PDE
\begin{align}\label{1}
\begin{cases}
(t+\nabla)f(t)&=0, \ t\in \R^d, \\
f(0)&=1.
\end{cases}
\end{align}
Note that the operator $Af(x)=(x+\nabla)f(x)$ is a multivariate Stein operator in the following sense: For a centred random vector $X$ with $\E [XX^\intercal]=\ID$, which has a differentiable density with full support $\R^d$, we have $\E [Af(X)]=\E[Xf(X)+\nabla f(X)]=0$ for each function $f$ with existing derivatives in every direction and for which all occurring expectations exist, if and only if $X$ has the normal distribution N$_d(0,\ID)$, see Theorem 3.5 in \cite{MRS:18} as well as \cite{LVY:13,L:94,S:81} for more information on the multivariate Stein lemma. Here and in the following the symbol $^\intercal$ means transposition of column vectors and matrices. In the spirit of the Stein-Tikhomirov method, see \cite{AMPS:16,FF:13}, and hence using the characteristic functions $\{\exp(\ii t^\intercal x), t\in\R^d\}$ as test functions, a simple calculation shows the equivalence of the Stein equation to the initial value problem in \eqref{1}. In the case $d=1$ the same initial value problem was motivated by a fixed point of the zero bias transform in \cite{E:20}. For more information on the zero bias transform, see \cite{GR:1997,S:2013}.
\begin{theorem}\label{thmcharac}
The characteristic function
\begin{equation}\label{defcfn01}
\psi(t)=\exp\bigg(-\frac{{\Vert t \Vert}^2}{2}\bigg), \ t\in \R^d,
\end{equation}
of the $d$-variate standard normal distribution  {\rm N}$_d(0,${\rm I}$_d)$ is the only solution of \eqref{1}.
\end{theorem}
\begin{prf}
If  $f \in L^2(\R^d)$ is an arbitrary solution of  \eqref{1}, the product rule yields
$$
\nabla\bigg(\exp\bigg(\frac{{\Vert t \Vert}^2}{2}\bigg)f(t)\bigg)=\exp\bigg(\frac{{\Vert t \Vert}^2}{2}\bigg)\bigg(t f(t)+ \nabla f(t)\bigg) =0.
$$
In view of $f(0) = 1$, we have $\exp(\Vert t \Vert^2/2) f(t)=1$, and the assertion follows.
\end{prf}

According to Theorem \ref{thmcharac}, the characteristic function (CF) of the $d$-variate standard  normal distribution  is the only CF satisfying
$\nabla \psi(t) = -t\psi(t)$.
Our test statistic will be based on this equation. To achieve affine invariance of the test statistic with respect to
full rank affine transformations of $X_1,\ldots,X_n$, let
\[
Y_{n,j} := S_n^{-1/2}(X_j - \overline{X}_n),\ j=1,...,n,
\]
denote the so-called scaled residuals, where $\overline{X}= n^{-1} \sum_{j=1}^nX_j$ and
$S_n := n^{-1} \sum_{j=1}^n(X_j-\overline{X}_n)(X_j-\overline{X}_n)^\intercal$  stand for the sample mean and the sample covariance matrix of
$X_1,\ldots,X_n$, respectively.
The matrix $S_n^{-1/2}$ is the unique symmetric positive definite square root of $S_n^{-1}$. To ensure almost sure invertibility of $S_n$, we tacitly
assume $n \ge d+1$ in what follows, see \cite{EP:73}. Writing
\begin{equation}\label{defpsint}
\psi_n(t)=\frac{1}{n}\sum_{j=1}^n \exp(\ii t^{\intercal}Y_{n,j}), \quad t \in \R^d,
\end{equation}
for the empirical CF of $Y_{n,1},\ldots,Y_{n,n}$,
our test statistic is
\begin{equation}\label{groesse}
T_{n,a} = n\int_{\R^d} {\Vert \nabla \psi_n(t) + t \psi(t) \Vert}_{\C}^2 \ w_a(t) \, \text{d}t.
\end{equation}
Here, $ w_a(t)=\exp\big(-a{\Vert t \Vert}^2\big)$, $a>0$, is a suitable weight function that depends on a positive parameter $a$, and $\Vert \cdot \Vert_{\C}$ denotes the complex Euclidean vector norm.
 Rejection of $H_0$ is for large values of $T_{n,a}$.
With this approach, we obtain a flexible class of genuine tests for multivariate normality, all of which are motivated by the result of Theorem \ref{thmcharac}.

Clearly, we propose a new approach to a well-known and widely studied problem, for a survey of affine invariant tests of multivariate normality, see \cite{H:2002}, and for recent developments with an emphasis on $L^2$ type statistics, see \cite{EH:20}. We list a short overview of different approaches: \cite{DEH:2019a,DEH:2019,HW:1997,P:2005,T:2009} consider tests connected to the empirical characteristic function, while \cite{HJG:2019,HJM:2019,HV:19} are based on the
 empirical moment generating function. The most classical approach is to consider measures of multivariate skewness and kurtosis, see, e.g., \cite{DH:2008,KTO:2007,M:70,MRS:94}, although inconsistency of those measures with regard to elliptically symmetric alternatives are known, see \cite{BH:1991,BH:1992,HEN:1994b,HEN:1994a}. Generalizations of tests for univariate normality, as in \cite{KP:2018,S:2006,AE:2009}, the examination of nonlinearity of dependence, see \cite{CS:1978,E:2012}, canonical correlations, see \cite{T:2014}, and the notion of energy, see \cite{SR:05}, are other approaches to this testing problem. Empirical competitive Monte Carlo studies can be found in \cite{EH:20,VPMV:2016}.

The rest of this paper unfolds as follows: In Section \ref{secbasprop}, we give a representation of $T_{n,a}$ that is amenable for computational purposes. Moreover,
we derive limits of $T_{n,a}$, after suitable affine transformations, as $a \to \infty$ and $a \to 0$, that hold
elementwise on the underlying probability space.
Section \ref{secnull} deals with the limit distribution of $T_{n,a}$ under the null hypothesis, and Section  \ref{secalternatives}
considers the limit behavior of $T_{n,a}$ both under contiguous and fixed alternatives to $H_0$.
Section \ref{secsimulation} presents the results of a simulation study, and  Section \ref{secrealdata} exhibits
a real data example. Section \ref{secsummary} contains a brief summary, and it indicates topics for further research.
For the sake of readability, some of the proofs have been deferred to Appendix \ref{secproofs}.

Throughout the paper, we use the following notation: The symbol $\edist$ means equality in distribution, and $\stk$ and $\fsk$ stand for
convergence in probability  and almost sure convergence, respectively. Moreover, $\vertk$ is shorthand for convergence in distribution
for random elements in whatever space (which will be clear from the context). If not stated otherwise, each limit refers to $n \to \infty$, and each unspecified integral is over $\R^d$. The stochastic Landau symbols $o_{\mathbb{P}}(1)$ and $O_{\mathbb{P}}(1)$ refer to convergence to zero in probability and stochastic boundedness, respectively.

\section{Basic properties of the test statistic}\label{secbasprop}
In this section, we provide some information on the test statistic $T_{n,a}$ defined in  \eqref{groesse}.
The first result shows that $T_{n,a}$ allows for a simple representation that is amenable to computational purposes. Moreover,
since this representation shows that $T_{n,a}$ depends on $X_1,\ldots,X_n$ only via $Y_{n,i}^\intercal Y_{n,j}$, $i,j \in \{1,\ldots,n\}$,
the statistic $T_{n,a}$ is affine invariant.

\begin{theorem}\label{themresprstat}
We have
\begin{eqnarray}\label{3}
T_{n,a} & =  & n {\bigg(\frac{\pi}{a+1}\bigg)}^\frac{d}{2} \frac{d}{2(a+1)} -2 {\bigg(\frac{2 \pi}{2a+1}\bigg)}^\frac{d}{2} \sum_{j=1}^n  \frac{{\Vert Y_{n,j} \Vert}^2}{2a+1} \exp\bigg(-\frac{{\Vert Y_{n,j} \Vert}^2}{4a+2}\bigg) \\
& & + \frac{1}{n}{\bigg(\frac{\pi}{a}\bigg)}^\frac{d}{2} \sum_{i,j=1}^n Y_{n,i}^{\intercal} Y_{n,j} \exp\bigg(-\frac{{\Vert Y_{n,i} - Y_{n,j} \Vert}^2}{4a}\bigg). \nonumber
\end{eqnarray}
\end{theorem}

Note that this representation is implemented in the \texttt{R} package \texttt{mnt}, see \cite{BE:2020}. The proof of Theorem \ref{themresprstat} is given in Appendix \ref{secproofs}.

We now consider the elementwise limits (on the underlying probability space) of $T_{n,a}$ for fixed $n$ as $a\rightarrow\infty$ and $a\rightarrow 0$.
It will bee seen that the class of tests based on $T_{n,a}$ is 'closed at the boundaries' $a \to \infty$ and $a \to 0$ in the sense that, after suitable
affine transformations, there are well-defined 'limit statistics'. Our first result refers to the limit $a \rightarrow \infty$.

\begin{theorem}\label{thmainfty}
Elementwise on the underlying probability space $(\Omega,{\cal A},\PP)$, we have
\begin{equation}\label{limainfsk}
\lim_{a\rightarrow\infty} \frac{a^{\frac{d}{2}+2}}{n\pi^{\frac{d}{2}}}16T_{n,a}=\widetilde{b}_{1,d}+2 b_{1,d}.
\end{equation}
Here, $b_{1,d}=n^{-2} \sum_{i,j=1}^n (Y_{n,i}^{\intercal} Y_{n,j})^3$ is Mardia's celebrated measure of multivariate skewness, see \cite{M:70}, and\linebreak $\widetilde{b}_{1,d}=n^{-2}\sum_{i,j=1}^n Y_{n,i}^{\intercal} Y_{n,j} \Vert Y_{n,i} \Vert^2 \Vert Y_{n,j} \Vert^2$ is a measure of multivariate skewness introduced by M\'{o}ri,  Rohatgi, and Sz\'{e}kely, see \cite{MRS:94}.
\end{theorem}

\begin{prf}
Invoking  \eqref{3}, it follows that
\begin{align*}
\frac{a^{\frac{d}{2}+2}}{n\pi^{\frac{d}{2}}}T_{n,a}&= {\bigg(\frac{a}{a+1}\bigg)}^{\frac{d}{2}+1} \frac{ad}{2} \ - \ \frac{a}{n} {\bigg(\frac{a}{a+\frac{1}{2}}\bigg)}^{\frac{d}{2}+1} \sum_{j=1}^n  {\Vert Y_{n,j} \Vert}^2 \exp\bigg(-\frac{{\Vert Y_{n,j} \Vert}^2}{4a+2}\bigg) \\
&\ \ \ \ + \frac{a^2}{n^2}\sum_{i,j=1}^n Y_{n,i}^{\intercal} Y_{n,j} \exp\bigg(-\frac{{\Vert Y_{n,i} - Y_{n,j} \Vert}^2}{4a}\bigg)\\
&=: A_n - B_n + C_n
\end{align*}
(say). We now use
\begin{align} \label{taylor1}
\bigg(\frac{a}{a+1}\bigg)^{\frac{d}{2}+1}=\bigg(1+\frac{1}{a}\bigg)^{-\frac{d}{2}-1}=1-\bigg(\frac{d}{2}+1\bigg)\frac{1}{a}+O(a^{-2})
\end{align}
as $a \to \infty$ and
\begin{align} \label{taylor2}
\exp(-x)=1-x+\frac{1}{2}x^2+O(x^3)
\end{align}
as $x  \to 0$, and we employ the identities
 $\sum_{j=1}^n Y_{n,j}=0$,  $\sum_{j=1}^n \Vert Y_{n,j} \Vert^2=nd$ as well as
\begin{align*}
\sum_{i,j=1}^n Y_{n,i}^{\intercal} Y_{n,j} \Vert Y_{n,i} -Y_{n,j} \Vert^2 &=-2\sum_{i,j=1}^n (Y_{n,i}^{\intercal} Y_{n,j})^2=-2n^2d, \\
\sum_{i,j=1}^n Y_{n,i}^{\intercal} Y_{n,j} \Vert Y_{n,i} -Y_{n,j}  \Vert^4&=2n^2\widetilde{b}_{1,d}+4 n^2 b_{1,d}- 8 \sum_{i,j=1}^n (Y_{n,i}^{\intercal} Y_{n,j})^2 \Vert Y_{n,j}\Vert^2, \\
\sum_{i,j=1}^n (Y_{n,i}^{\intercal} Y_{n,j})^2 \Vert Y_{n,j}\Vert^2&=n\sum_{j=1}^n \Vert Y_{n,j} \Vert^4.
\end{align*}
to obtain $A_n=ad/2 -d^2/4 - d/2 + o(1)$ as $a \to \infty$. Likewise,
\begin{align*}
B_n&=\frac{1}{n}\bigg(a-\bigg(\frac{d}{2}+1\bigg) \frac{1}{2}\bigg)\sum_{j=1}^n \Vert Y_{n,j} \Vert^2 \bigg(1-\frac{\Vert Y_{n,j} \Vert^2}{4a+2}\bigg)+o(1)\\
&=\bigg(da-\frac{d^2}{4}-\frac{d}{2}\bigg)-\frac{1}{4n}\sum_{j=1}^n \Vert Y_{n,j} \Vert^4+o(1), \\
C_n&=\frac{a^2}{n^2}\sum_{i,j=1}^n Y_{n,i}^{\intercal} Y_{n,j} \bigg(1-\frac{{\Vert Y_{n,i} - Y_{n,j} \Vert}^2}{4a}+\frac{{\Vert Y_{n,i} - Y_{n,j} \Vert}^4}{32a^2}\bigg)+o(1)\\
&=\frac{da}{2}+\frac{1}{16}\Big(\widetilde{b}_{1,d}+2 b_{1,d} - \frac{4}{n}\sum_{j=1}^n \Vert Y_{n,j} \Vert^4 \Big) +o(1).
\end{align*}
Upon combining, the assertion follows.
\end{prf}

Notice that the right hand side of \eqref{limainfsk} is a linear combination of two time-honored measures of multivariate skewness. Notably,
the same linear combination showed up not only for the class of BHEP tests (see Theorem 2.1 of \cite{H:97}), but also as a
limit of a related test statistic in connection with a test for multivariate normality
based on a partial differential equation for the {\em moment generating function} of the normal distribution, see \cite{HV:19}.

Regarding the limit of $T_{n,a}$ as $a \to 0$, we have the following result.

\begin{theorem}\label{thmato0}
Elementwise on the underlying probability space, we have
\begin{align*}
\lim_{a\rightarrow 0} \frac{1}{na^{\frac{d}{2}}}\left(\left(\frac{a}{\pi}\right)^{\frac{d}{2}}T_{n,a}-d\right)=\frac{d}{2}-2^{\frac{d}{2} + 1}\frac1n\sum_{j=1}^n \Vert Y_{n,j} \Vert^2 \exp\left(-\frac{\Vert Y_{n,j} \Vert^2}{2}\right).
\end{align*}
\end{theorem}

\begin{prf} From the representation \eqref{3}, it follows that
\begin{align*}
\frac{T_{n,a}}{\pi^{\frac{d}{2}}} &=   \frac{nd}{2\left(a+1\right)^{\frac{d}{2}+1}} -
{\left(\dfrac{2}{2a+1}\right)}^{\frac{d}{2}+1} \sum_{j=1}^n  \Vert Y_{n,j} \Vert^2 \exp\left(-\frac{{\Vert Y_{n,j} \Vert}^2}{4a+2}\right) \\
&\ \ \ \ + \frac{1}{na^\frac{d}{2}} \sum_{i,j=1}^n Y_{n,i}^{\intercal} Y_{n,j} \exp\left(-\frac{{\Vert Y_{n,i} - Y_{n,j} \Vert}^2}{4a}\right) \\
&= A_{n,a} - B_{n,a} + C_{n,a}
\end{align*}
(say). Now, $\lim_{a\rightarrow 0} A_{n,a} = nd/2$ and
$\lim_{a\rightarrow 0} B_{n,a} = 2^{\frac{d}{2} + 1}\sum_{j=1}^n \Vert Y_{n,j} \Vert^2 \exp\left(-\Vert Y_{n,j} \Vert^2/2\right)$, elementwise on the underlying probability space.
To tackle $C_{n,a}$, the relation $\sum_{j=1}^n \Vert Y_{n,j} \Vert^2 = nd$ yields
\begin{align*}
C_{n,a}&= \frac{1}{na^\frac{d}{2}}\sum_{j=1}^n \Vert Y_{n,j} \Vert^2 + \frac{1}{na^\frac{d}{2}} \sum_{i\neq j}^n Y_{n,i}^{\intercal} Y_{n,j} \exp\left(-\frac{{\Vert Y_{n,i} - Y_{n,j} \Vert}^2}{4a}\right) \\
&= \frac{d}{a^\frac{d}{2}} + \frac{1}{na^\frac{d}{2}} \sum_{i\neq j}^n Y_{n,i}^{\intercal} Y_{n,j} \exp\left(-\frac{{\Vert Y_{n,i} - Y_{n,j} \Vert}^2}{4a}\right),
\end{align*}
and the assertion follows.
\end{prf}

Interestingly, Theorem \ref{thmato0} means that for (very) small values of $a$, rejection of  $H_0$ for large values of $T_{n,a}$ is essentially equivalent to
the rejection of $H_0$ for {\em small} values of
\[
\frac1n\sum_{j=1}^n \Vert Y_{n,j} \Vert^2 {\text e}^{-\Vert Y_{n,j}\Vert^2/2}.
\]
This statistic, upon expanding the exponential function, comprises even powers of $\|Y_{n,j}\|$ and is thus related to Mardia's measure of multivariate kurtosis,
which is defined by $b_{2,d} = n^{-1}\sum_{j=1}^n \|Y_{n,j}\|^4$, see \cite{M:70}.

\section{The limit null distribution of $T_{n,a}$}\label{secnull}

In this section we derive the limit distribution of $T_{n,a}$ under the hypothesis $H_0$. In view of affine invariance, we assume without loss of generality
that $X$ has the standard normal distribution N$_d(0,\ID)$ in what follows. The starting point is an alternative representation of $T_{n,a}$, namely
\begin{align} \label{groesse2}
T_{n,a}=\int {\Vert Z_n(t) \Vert}^2 w_a(t) \ \text{d}t,
\end{align}
where
\begin{align}
Z_n(t) = \dfrac{1}{\sqrt{n}} \sum_{j=1}^n \bigg(Y_{n,j}\big(\cos(t^{\intercal}Y_{n,j}) + \sin(t^{\intercal}Y_{n,j})\big)-t \psi(t) \bigg). \label{Zn}
\end{align}
This assertion follows from straightforward calculations using
\begin{align} \label{symmetrie1}
\int \cos(t^{\intercal}Y_{n,j}) \sin(t^{\intercal} Y_{n,i}) w_a(t) \ \text{d}t = 0, \quad
\int \cos(t^{\intercal}Y_{n,j}) t^{\intercal}Y_{n,j} w_a(t) \ \text{d}t =0.
\end{align}

Writing $L^2:= L^2(\R^d,\mathcal{B}^d,w_a(t)\text{d}t)$ for the separable Hilbert space of (equivalence classes of) functions $f:\R^d \rightarrow \R$ that are square integrable with respect to
 $w_a(t)\text{d}t$, we regard $Z_n$ as a random element of the Hilbert space
$ \HH=L^2 \otimes \cdot \cdot \cdot\otimes L^2 $. Putting $f=(f_1,\ldots,f_d),g=(g_1,\ldots,g_d)$, the space $\HH$ is equipped with the inner product $\langle f, g \rangle_{\HH} := \langle f_1,g_1 \rangle_{L^2} + \ldots + \langle f_d, g_d \rangle_{L^2} $ and the
norm ${\Vert f \Vert}_{\HH} = \langle f, f \rangle_{\HH}^{1/2}$.  Notice that we have
 \[
T_{n,a} =  \int {\Vert Z_n(t) \Vert}^2 w_a(t) \ \text{d}t  = {\Vert Z_n \Vert}^2_{\HH}.
\]

The main theorem of this section is as follows:

\begin{theorem}\label{thmlimitnull}
Under $H_0$, there is a centred Gaussian random element $Z$ of $\HH$ having covariance matrix kernel
\begin{align} \label{kovark}
K(s,t)&=\big({\rm I}_d - (s-t)(s-t)^{\intercal}\big)\psi(s-t)\\ & \ \ \ + \Big(ss^{\intercal}  + tt^{\intercal} - t s^{\intercal} - st^{\intercal} - {\rm I}_d
+ s^{\intercal}t(ss^{\intercal} + t t^{\intercal} -st^{\intercal} - {\rm I}_d) - \dfrac{s^{\intercal}ts^{\intercal} t}{2}  s t^{\intercal} \Big)\psi(s)\psi(t) ,\ \ \ s,t \in \R^d, \nonumber
\end{align}
such that $Z_n \overset{\mathcal{D}}{\longrightarrow} Z \text{ in } \HH$,
where $Z_n$ is the random element defined in \eqref{Zn}.
\end{theorem}

Since the proof of Theorem \ref{thmlimitnull} is long and tedious, it is deferred to Appendix \ref{secproofs}.
From Theorem \ref{thmlimitnull} and the continuous mapping theorem, we obtain the following result.

\begin{corollary}\label{42}
Under $H_0$, we have
\begin{align*}
T_{n,a}\overset{\mathcal{D}}{\longrightarrow}\Vert Z \Vert_{\HH}^2 = \int \Vert Z(t) \Vert^2 w_a(t) \, {\rm d}t.
\end{align*}
\end{corollary}

It is well-known that the distribution of $T_{\infty,a}:=\Vert Z \Vert^2_{\HH}$ is that of $T_{\infty,a}\overset{\mathcal{D}}{=} \sum_{j=1}^{\infty} \lambda_j(a) N_j^2$, where $N_1,N_2, \ldots$ is a sequence of i.i.d. standard normal random variables, and $\lambda_1(a),\lambda_2(a), \ldots$ are the positive eigenvalues
associated with the integral operator
\begin{align}\label{Inteq}
\mathbb{K}f(s):=\int K(s,t)f(t)w_a(t) \ \text{d}t, \ \ s \in \R^d,
\end{align}
 $f\in \HH$. In view of the complexity of $K(s,t)$, we did not succeed in obtaining closed-form expressions for
     these eigenvalues. In our simulation study presented in Section \ref{secsimulation}, we use approximate critical values for $T_{n,a}$ that have been
      obtained by means of simulations. Some information on the limit null distribution, however, is given by the following result.

\begin{theorem}
We have
\begin{align*}
\E [T_{\infty,a}]=\Big(\dfrac{\pi}{a}\Big)^{\frac{d}{2}}d-\Big(\dfrac{\pi}{a+1}\Big)^{\frac{d}{2}}\dfrac{\big(16a^3+(8d+48)a^2+(12d+40)a+d^2+10d+16\big)d}{16(a+1)^3}.
\end{align*}
\end{theorem}
\begin{prf} From Fubini's theorem, it follows that $\E [T_{\infty,a}]=\int \E \Vert Z(t) \Vert ^2  w_a(t) \, \text{d}t$.
Moreover, writing tr for trace, we have
\begin{equation*}
\E \Vert Z(t) \Vert^2= \E[Z(t)^{\intercal}Z(t)]= \text{tr} \big(\E[Z(t)Z(t)^{\intercal}]\big) = \text{tr}\big(K(t,t)\big)=d -\Big(d  +d\Vert t \Vert^2 - \Vert t \Vert^4 + \frac{\Vert t \Vert^6}{2}\Big)\exp\big(-\Vert t\Vert^2\big).
\end{equation*}
Since
\[
\int \Vert t \Vert^4 {\text e}^{-a\Vert t \Vert^2}\, \text{d}t=\Big(\dfrac{\pi}{a}\Big)^\frac{d}{2}\dfrac{d}{4a^2}(d+2) \quad\mbox{and}\quad
\int \Vert t \Vert^6 {\text e}^{-a\Vert t \Vert^2}\, \text{d}t=\Big(\dfrac{\pi}{a}\Big)^\frac{d}{2} \dfrac{d}{8a^3}(d^2+6d+8),
\]
  the assertion follows by straightforward computations.
\end{prf}

In the univariate case, which is deliberately not excluded from our study,
 we have been able to calculate the first four cumulants of $T_{\infty,a}$. By the methods presented in Chapter 5 of \cite{SW:86} the $m$th cumulant of $T_{\infty,a}$ is derived by
\begin{align*}
\kappa_m(a)=2^{m-1} (m-1)!\int_{\R} h_m(t,t)w_a(t)\ \text{d}t.
\end{align*}
Here, $h_1(s,t) =K(s,t)$, and $h_m(s,t):=\int_{\R} h_{m-1}(s,u)K(u,t)w_a(u) \, {\text  d}u$ if $m \ge 2$. In order to calculate $\kappa_m(a)$, $m \in \{1,2,3,4\}$, we used the computer algebra system Maple, see \cite{Maple2019}.
The formulae for $\kappa_3(a)$ and $\kappa_4(a)$ are given in the appendix.

For $\kappa_1(a)$ and $\kappa_2(a)$, we obtain
\begin{align*}
\kappa_1(a)&=\int_{\R} \big(1-(1+t^2-t^4+\dfrac{t^6}{2})\exp(-t^2)\big)\exp(-at^2) \ \text{d}t \\
&=\sqrt{\dfrac{\pi}{a}}-\sqrt{\dfrac{\pi}{a+1}}-\sqrt{\dfrac{\pi}{a+1}}\dfrac{1}{2(a+1)}+\sqrt{\dfrac{\pi}{a+1}}\dfrac{3}{4(a+1)^2}-\sqrt{\dfrac{\pi}{a+1}}\dfrac{15}{16(a+1)^3} \\ &=
\dfrac{(-16a^3 - 56a^2 - 52a - 27)\sqrt{\dfrac{\pi}{a + 1}} + 16\sqrt{\dfrac{\pi}{a}}(a + 1)^3}{16(a + 1)^3}
\end{align*}
and
\begin{align*}
\kappa_2(a)=&\frac {7260811 \pi}{8 \left( a+2 \right) ^{5/2} \left( 4 {a}^{2}+8a+3 \right) ^{5/2}\sqrt {a} \left( 2a+3 \right) ^{2} \left( a+1 \right) ^{7}}  \bigg(  \bigg(  \bigg( \frac {1024}{7260811}{a}^{\frac{29}{2}}\\
&+\frac {15360}{7260811}{a}^{\frac{27}{2}}+\frac {108032}{7260811}{a}^{\frac{25}{2}}+\frac {473856}{7260811}{a}^{\frac{23}{2}}+\frac {1449216}{7260811}{a}^{\frac{21}{2}}+{\frac {3263232}{7260811}}{a}^{\frac{19}{2}}\\
&+\frac {5559908}{7260811}{a}^{\frac{17}{2}}+\frac {7254348}{7260811}{a}^{\frac{15}{2}}+{a}^{\frac{13}{2}}+\frac {5535906}{7260811}{a}^{\frac{11}{2}}+\frac {160113}{367636}{a}^{\frac{9}{2}}\\
&+\frac {5253759}{29043244}{a}^{
\frac{7}{2}}+\frac {6017409}{116172976}{a}^{\frac{5}{2}}+\frac {266733}{29043244}{a}^{\frac{3}{2}}+\frac {22113}{29043244}\sqrt {a} \bigg) \sqrt {a+2}\\
&+\frac {1024 \left( a+3/2 \right) ^{4} \left( a
+1 \right) ^{7} \left( a+1/2 \right) ^{2} \left( {a}^{2}+2a+3 \right) }{7260811} \bigg) \sqrt {4{a}^{2}+8a+3}\\
&-\frac {51420992\sqrt {a+2}}{7260811} \bigg( \frac {64}{803453}{a}^{\frac{31}{2}}+\frac {1024}{803453}{a}^{\frac{29}{2}}+\frac {1104}{
114779}{a}^{\frac{27}{2}}+\frac {36544}{
803453}{a}^{\frac{25}{2}}\\
&+\frac {121054}{803453}{a}^{\frac{23}{2}}+\frac {297018}{803453}{a}^{\frac{21}{2}}+\frac {556163}{803453}{a}^{\frac{19}{2}}+\frac {807017}{803453}{a}^{\frac{17}{2}}+\frac {912747}{803453}{a}^{\frac{15}{2}}+{a}^{\frac{13}{2}}\\
&+\frac {545801}{803453}{a}^{\frac{11}{2}}+\frac {281319}{803453}{a}^{\frac{9}{2}}+\frac {106779}{803453}{a}^{\frac{7}{2}}+\frac {28293}{803453}{a}^{\frac{5}{2}}+\frac {96}{16397}{a}^{\frac{3}{2}}+\frac {372}{803453}\sqrt {a}
 \bigg)  \bigg).
\end{align*}
From these cumulants, we obtain the expectation, the variance as well as the skewness $\beta_1$ and the kurtosis $\beta_2$
of $T_{\infty,a}$ for the case $d=1$ (see Table \ref{schief}), since
\[
\E [T_{\infty,a}] =\kappa_1(a),  \ \
\text{Var}[T_{\infty,a}] =\kappa_2(a), \ \
\beta_1(a) = \frac{\kappa_3(a)}{\kappa_2(a)^{3/2}}, \ \
\beta_2(a) =3+\dfrac{\kappa_4(a)}{\kappa_2(a)^2}.
\]
By complete analogy with  \cite{E:20,H:1990}, we can now approximate the distribution of
 $T_{\infty,a}$ by that of a member of the system of Pearson distributions which has the same first four moments as $T_{\infty,a}$.
To this end, we used the statistic software \texttt{R}, see \cite{R:2019}, and the package \texttt{PearsonDS}, see \cite{BK:2017}.
Table \ref{qu} shows the quantiles of the fitted Pearson distribution, which serve as approximations of the corresponding
quantiles of the distribution of $T_{\infty,a}$.

\begin{table}[t]
\centering
\renewcommand{\arraystretch}{1.3}
\begin{tabular}{ r|rrrrrr }
\hline\hline
 $a$ & 0.1 & 0.5 & 1 & 2 & 5 & 10  \\ \hline
$\E[T_{\infty,a}]$ & 3.0040 & 0.6574 & 0.2939 & 0.1092 & 0.0207 & 0.0047 \\
$\text{Var}[T_{\infty,a}]$ & 2.8028 & 0.2686& 0.0742 & 0.0133 & 0.0006 & 0.0000\\
$\beta_1(a)$ & 1.3737 & 1.9098 & 2.1996 & 2.4619 & 2.7090 & 2.7938\\
$\beta_2(a)$ & 6.0366 & 8.8662 & 10.7047 & 12.5510 & 14.3071 & 19.4464 \\
\hline\hline
\end{tabular}
\smallskip
\caption{Expectation, variance, skewness and kurtosis  of $T_{\infty,a}$,  $d=1$}
\label{schief}
\end{table}

\renewcommand{\arraystretch}{1.3}
\begin{table}[t]
\centering
\begin{tabular}{ r|rrrrrr }
\hline\hline
 $q\setminus a$ & 0.1 & 0.5 & 1 & 2 & 5 & 10 \\ \hline
0.01 & 0.6857 & 0.0903 & 0.0331 & 0.0110 & 0.0018 & -0.0013\\
0.05 & 0.9970 & 0.1299 & 0.0435 & 0.0130 & 0.0020 & -0.0009\\
0.1 & 1.2382 & 0.1712 & 0.0573 & 0.0165 & 0.0023 &  -0.0005\\
0.5 & 2.6510 & 0.5137 & 0.2091 & 0.0700 & 0.0115 & 0.0030\\
0.9 & 5.2211 & 1.3283 & 0.6405 & 0.2529 & 0.0511 & 0.0119\\
0.95 & 6.2138 & 1.6743 & 0.8329 & 0.3384 & 0.0705 & 0.0162\\
0.99 & 8.4485 & 2.4904 & 1.2956 & 0.5470 & 0.1182 & 0.0275\\
\hline\hline
\end{tabular}
\smallskip
\caption{Approximate quantiles of $T_{\infty,a}$  in the case $d=1$}
\label{qu}
\end{table}

\section{Limit behavior of $T_{n,a}$ under alternatives}\label{secalternatives}
In this section, we assume that $H_0$ does not hold, and we will derive limit distributions for $T_{n,a}$ both under contiguous and fixed alternatives
to $H_0$. To define the setting for a triangular array of contiguous alternatives, we assume that, for each $n \ge d+1$,  $X_{n,1}, ... , X_{n,n}$ are i.i.d.
$d$-variate random vectors having Lebesgue density
\[
f_n(x)=\varphi(x)\Big(1+ \dfrac{g(x)}{\sqrt{n}}\Big), \ x \in \R^d.
\]
Here, $\varphi(x)=(2\pi)^{-d/2}\exp(-\Vert x \Vert^2/2)$, $x\in \R^d$, is the density of the distribution N$_d(0,\ID)$,
and $g$ is a bounded measurable function satisfying $\int g(x)\varphi(x) \, {\rm d}x = 0$. Notice that $f_n$ is nonnegative for sufficiently large
$n$ due to the boundedness of $g$.  To derive the limit distribution of $T_{n,a}$ under
this sequence of alternatives, we employ the representation \eqref{groesse2}, which comprises the random element $Z_n$
as defined in  \eqref{Zn}. For repeated later use, we put
\begin{equation}\label{CSp}
\CSP(s,t) = \cos(s^{\intercal}t) + \sin(s^{\intercal}t), \quad  \CSM(s,t) = \cos(s^{\intercal}t) - \sin(s^{\intercal}t), \quad  s,t \in \R^d.
\end{equation}

\begin{theorem}\label{thmcontig} Under the sequence of alternatives $(X_{n,1}, ... , X_{n,n})_{n \ge d+1}$, we have
\[
Z_n \vertk  Z+c \text{ in } \ \HH.
\]
Here, $Z_n$ is defined in  \eqref{Zn}, $Z$ is the centred Gaussian random element of $\HH$ figuring in Theorem \ref{thmlimitnull}, and the shift function $c(\cdot)$ is given by
\begin{equation}\label{defshiftfkt}
c(t)=\int Z^{**}(x,t)g(x)\varphi(x) \, {\rm d} x,\quad  t \in \R^d,
\end{equation}
where
\[
Z^{**}(x,t)=x {\rm CS}^+(t,x)-\big(t+x+ (2{\rm I}_d -tt^{\intercal}) \dfrac{1}{2}(xx^{\intercal}-{\rm I}_d )t - t^{\intercal}xt\big)\psi(t), \ \ x,t \in \R^d.
\]
\end{theorem}

\begin{prf} We write $\lambda^d$ for $d$-dimensional Lebesgue measure, and we put
 $\mathbb{P}^{(n)}:= \otimes (\varphi \lambda^d)$, $Q^{(n)}:= \otimes (f_n \lambda^d)$.
Furthermore, let $L_n := {\text d}Q^{(n)}/{\text d}\mathbb{P}^{(n)}$. The boundedness of $g$ and a Taylor expansion then give
\begin{equation}\label{contig1}
\log(L_n(X_{n,1},...,X_{n,n}))=\sum_{j=1}^n \log\Big(1+\dfrac{g(X_{n,j})}{\sqrt{n}}\Big)=\sum_{j=1}^n \Big(\dfrac{g(X_{n,j})}{\sqrt{n}} - \dfrac{g(X_{n,j})^2}{2n} \Big) + o_{\mathbb{P}^{(n)}}(1).
\end{equation}
In the following we write $\sigma^2 =\int g(x)^2 \varphi(x) \, {\text d}x <\infty$. Since, under $\mathbb{P}^{(n)}$, expectation and
variance of the sum figuring in \eqref{contig1} converge to $-\sigma^2/2$ and $\sigma^2$, respectively,
the Lindeberg--Feller central limit theorem and Slutsky's lemma yield
\begin{align} \label{logn}
\log(L_n) \vertk {\rm N}\Big{(}-\dfrac{\sigma^2}{2},\sigma^2\Big{)}  \text{ under } \mathbb{P}^{(n)}.
\end{align}
Notice that  the boundedness of $g$ ensures the validity of the Lindeberg condition.
In view of Le Cam's first lemma (see, e.g., \cite{LB:19}, p. 297), the probability measures
      $Q^{(n)}$ and $\mathbb{P}^{(n)}$ are mutually contiguous.
According to Theorem \ref{thmlimitnull}, the auxiliary process $Z_n^*$ introduced in \eqref{Zn*} is tight
under $\mathbb{P}^{(n)}$ and thus, in view of contiguity, also under $Q^{(n)}$. Let $\{e_k, \ k\geq 1 \}$
be an arbitrary complete orthonormal system of $\HH$. It remains to
show that, for each $\ell \ge 1$, we have $\Pi_{\ell}(Z_n) \vertk \Pi_{\ell}(Z+c)$ under $Q^{(n)}$, where
$\Pi_\ell$ denotes the orthogonal projection onto the linear subspace of $\HH$ spanned by $e_1,\ldots,e_\ell$.
We first consider
\[
\Pi_{\ell}(Z_n^*)=\sum_{j=1}^{\ell} \langle Z_n^*,e_j \rangle_{\HH} e_j,
\]
where $Z_n^*$ is given in \eqref{Zn*}, with the only difference that $X_j$ is throughout replaced with $X_{n,j}$.
In view of Theorem \ref{thmlimitnull}, the asymptotic distribution of $Z_n^*$ under  $\mathbb{P}^{(n)}$ is centred Gaussian
 with a covariance operator $\mathbb{K}$ given by the covariance matrix kernel $K(s,t)$, whence
$\langle Z_n^*,e_j \rangle_{\HH} \vertk  {\text N}(0,\langle \mathbb{K}e_j,e_j\rangle_{\HH})$ under $\mathbb{P}^{(n)}$.
In view of \eqref{logn} we have
\begin{align*}
 \big(\langle Z_n^*,e_1 \rangle_{\HH}, ... , \langle Z_n^*,e_{\ell} \rangle_{\HH}, \log (L_n)\big)^{\intercal} \vertk  {\text N}_{\ell+1} \Bigg(
 (0, \ldots,0,-\sigma^2/2)^{\intercal},
 \begin{bmatrix}
 \Sigma && \widetilde{c} \\ \widetilde{c}^{\intercal} && \sigma^2 \end{bmatrix}\Bigg)
\end{align*}
under $\mathbb{P}^{(n)}$ for each $\ell \ge 1$. Here,
 $\Sigma:=\big(\langle \mathbb{K}e_i, e_j \rangle_{\HH}\big)_{1\leq i,j\leq \ell} \in \R^{\ell\times \ell} $ and $\widetilde{c}=\big(\widetilde{c}_1, ..., \widetilde{c}_\ell\big)^{\intercal} \in \R^{\ell}$,
 where, by Fubini's theorem, $\widetilde{c}_j:=\lim_{n\rightarrow \infty} \E \big[\langle Z_n^*,e_j \rangle_{\HH},\log(L_n)\big]=\langle c, e_j \rangle_{\HH}$,
and $c$ is given in \eqref{defshiftfkt}. According to
Le Cam's third Lemma (see, e.g., \cite{LB:19}, p. 300), it follows that
$
 \big(\langle Z_n^*,e_1 \rangle_{\HH}, ... , \langle Z_n^*,e_{\ell} \rangle_{\HH}\big)^{\intercal} \vertk  {\text N}_{\ell} (\widetilde{c},\Sigma)
$
under $Q^{(n)}$. Since, for the centred Gaussian random element figuring in Theorem \ref{thmlimitnull}, we have
\begin{align*}
\big(\langle Z + c,e_1 \rangle_{\HH}, ... , \langle Z + c,e_{\ell} \rangle_{\HH}\big)^{\intercal} \overset{\mathcal{D}}{=} N_{\ell} (\widetilde{c},\Sigma),
\end{align*}
it follows that
\begin{align}\label{verschoben}
  \big(\langle Z_n^*,e_1 \rangle_{\HH}, ... , \langle Z_n^*,e_{\ell} \rangle_{\HH}\big)^{\intercal} \vertk  \big(\langle Z + c,e_1 \rangle_{\HH}, ... , \langle Z + c,e_{\ell} \rangle_{\HH}\big)^{\intercal}
\end{align}
under $\mathbb{Q}^{(n)}$. Now, let $\Psi: \R^{\ell} \rightarrow \HH$ be defined by  $\Psi(x):= \sum_{j=1}^\ell x_j e_j$,
  $x=(x_1, ... , x_{\ell})^{\intercal}$. The continuous mapping theorem and \eqref{verschoben} then yield
\begin{align*}
\Pi_{\ell}(Z_n^*)=\Psi\Big(\big(\langle Z_n^*,e_1 \rangle_{\HH}, ... , \langle Z_n^*,e_{\ell} \rangle_{\HH}\big)^{\intercal}\Big) \overset{\mathcal{D}}{\longrightarrow} &\Psi \Big(\big(\langle Z + c,e_1 \rangle_{\HH}, ... , \langle Z + c,e_{\ell} \rangle_{\HH}\big)^{\intercal} \Big) \\&= \Pi_{\ell}(Z+c)
\end{align*}
under $\mathbb{Q}^{(n)}$. In view of the tightness of $Z_n^*$ unter $Q^{(n)}$ we conclude
$
Z_n^* \overset{\mathcal{D}}{\longrightarrow} Z + c \text{ under }Q^{(n)}$.
The assertion now follows from Slutsky's lemma since, in view of \ref{Schritt2} and \ref{Schritt3}, ${\Vert Z_n - Z_n^* \Vert}_\HH$ is asymptotically negligible under $\mathbb{P}^{(n)}$ and thus,
because of contiguity, also under $Q^{(n)}$.
\end{prf}

As a corollary, we have the following result.

\begin{corollary} Under the conditions of Theorem  \ref{thmcontig}, we have
\begin{align*}
T_{n,a}\overset{\mathcal{D}}{\longrightarrow}\Vert Z + c \Vert_{\HH}^2= \int \big \Vert Z(t) + c(t) \big \Vert^2 w_a(t) \, {\rm d}t.
\end{align*}
\end{corollary}

We now consider fixed alternatives to $H_0$, and we suppose that the underlying distribution, in addition to being absolutely continuous,
satisfies $\E\Vert X \Vert^4<\infty$. In view of affine invariance, we assume
$\E [X] = 0$ and $\E [XX^{\intercal}]=\ID$. Our first result is a strong limit of $T_{n,a}/n$ as $n \to \infty$.

\begin{theorem}\label{themaslimit}
If $\E \Vert X \Vert^4 < \infty$, we have
\begin{align*}
\dfrac{T_{n,a}}{n} \fsk \Delta_a,
\end{align*}
where
\begin{equation}\label{defdeltaa}
\Delta_a := \int \Vert \mu(t) - t\psi(t) \Vert^2 w_a(t) \, {\rm  d}t
\end{equation}
and $\mu(t)=\E [X {\rm CS}^+(t,X)]$.
\end{theorem}

\begin{prf}
Invoking \eqref{groesse2}, we have $n^{-1}T_{n,a} = \Vert n^{-1/2} Z_n \Vert_{\HH}^2$, where $Z_n$ is given in \eqref{Zn}. Putting
$
Z^0_n(t)= n^{-1/2}\sum_{j=1}^n \big(X_j\CSP(t,X_j)-t \psi(t)\big)$,
the strong law of large numbers in Hilbert spaces yields $\Vert n^{-1/2} Z_n^0 \Vert_{\HH}^2 \fsk 0$, and thus it remains
to prove $\Vert n^{-1/2} (Z_n - Z_n^0) \Vert_{\HH} \fsk 0$. To this end, notice that
\begin{align*}
\frac{1}{\sqrt{n}} \big(Z_n(t) - Z_n^0(t)\big)&=\dfrac{1}{n} \sum_{j=1}^n \Big(X_j\big(\CSP(t,Y_{n,j})-\CSP(t,X_j)\big) + \Delta_{n,j} \CSP(t,Y_{n,j})\Big).
\end{align*}
Since $\CSP(t,Y_{n,j})=\CSP(t,X_j)+\varepsilon_{n,j}(t)+\eta_{n,j}(t)$,
where $\max(|\varepsilon_{n,j}(t)|,|\eta_{n,j}(t)|) \leq \Vert t \Vert \Vert \Delta_{n,j} \Vert$, it follows that
\begin{align*}
\Big\Vert \dfrac{1}{n}\sum_{j=1}^n X_j\big(\CSP(t,Y_{n,j})-\CSP(t,X_j)\big) \Big\Vert &\leq \dfrac{2}{n} \sum_{j=1}^n \Vert X_j \Vert \Vert t \Vert \Vert \Delta_{n,j} \Vert \\
&\leq 2 \Vert t \Vert n^{-1/4} \max_{j=1,...,n} \Vert X_j \Vert n^{1/4} \max_{j=1,...,n} \Vert \Delta_{n,j} \Vert.
\end{align*}
Since $\E \Vert X \Vert^4 < \infty$, Theorem 5.2 of \cite{B:63} yields $n^{-1/4} \max_{j=1,...,n} \Vert X_j \Vert \fsk 0$, and from
Proposition A.1 of \cite{DEH:2019}, we have  $n^{1/4} \max_{j=1,...,n} \Vert \Delta_{n,j} \Vert \fsk 0$.
Consequently,
$
\Vert n^{-1}\sum_{j=1}^n X_j\big(\CSP(t,Y_{n,j})-\CSP(t,X_j)\big) \Vert_{\HH} \fsk 0$.
Furthermore, $\Vert n^{-1} \sum_{j=1}^n \Delta_{n,j} \CSP(t,Y_{n,j}) \Vert \leq 2n^{-1} \sum_{j=1}^n \Vert \Delta_{n,j} \Vert$.
Since the right hand side converges to 0 $\PP$-almost surely according to  Proposition A.1 of \cite{DEH:2019},
it follows that $\Vert n^{-1} \sum_{j=1}^n \Delta_{n,j} \CSP(t,Y_{n,j}) \Vert_{\HH} \fsk 0$.
The remaining assertion $\Vert n^{-1/2} (Z_n - Z_n^0) \Vert_{\HH} \fsk 0$ now follows from the triangle inequality.
\end{prf}

As a corollary, we obtain the following result.

\begin{corollary} The test for multivariate normality based on $T_{n,a}$ is consistent against each alternative distribution
satisfying $\E \Vert X \Vert^4 < \infty$.
\end{corollary}

\begin{prf}
Let $\psi_X(t) = \E[\exp(\ii t^{\intercal}X)]$ be the CF of $X$. By straightforward calculations, we have
\begin{align*}
\Delta_ a = \int {\Vert \nabla \psi_X (t) - \nabla \psi(t) \Vert}_{\C}^2 w_a(t) \, {\text d}t,
\end{align*}
where $\Delta_a$ is given in \eqref{defdeltaa}. Since $\Delta_a =0$ if and only if $X \edist {\text N}_d(0,\ID)$ (recall the standing assumption
$\E[X] =0$ and $\E[XX^{\intercal}] = \ID$), the assertion follows.
\end{prf}

Notice that, for each $a >0$, $\Delta_a$ may be regarded as a measure of deviation from normality. The following result
sheds some more light on $\Delta_a$.

\begin{theorem}\label{thmdeltatoinfty} If $E\Vert X \Vert^6 < \infty$ then, under the standing assumptions $\E [X] = 0$ and $\E [XX^{\intercal}]={\rm I}_d$, we have
\begin{equation}\label{limitdeltaainfty}
\lim_{a\rightarrow\infty} 16a^2\Big(\dfrac{a}{\pi} \Big)^{\frac{d}{2}} \Delta_{a} = \E [ X_1^{\intercal} X_2 \Vert X_1 \Vert^2 \Vert X_2 \Vert^2] + 2 \E [(X_1^{\intercal}X_2)^3],
\end{equation}
as well as
\begin{align*}
\lim_{a\rightarrow 0} \pi^{-\frac{d}{2}}\Delta_a=\frac{d}{2}-2^{\frac{d}{2}+1}\E\left[\Vert X_1\Vert^2\exp\left( -\frac{\Vert X_1 \Vert^2}{2}\right)\right].
\end{align*}
\end{theorem}

\begin{prf}
Straightforward calculations give $\Delta_a = I_{a,1}- I_{a,2} + I_{a,3}$, where
\begin{align*}
I_{a,1} &= \int \E [X_1 \CSP(t,X_1)]^{\intercal} \E [X_2 \CSP(t,X_2)] w_a(t) \, {\text d}t \\
I_{a,2} & = 2 \int \E [X_1 \CSP(t,X_1)]^{\intercal} t \psi(t) w_a(t) \,  {\text d}t, \qquad
I_{a,3}  = \int t^{\intercal} t \psi(t)^2 w_a(t) \, {\text d}t.
\end{align*}
Using addition theorems for the sine and the cosine function as well as
 \eqref{symmetrie1} and \eqref{integrale1},\eqref{integrale2} und \eqref{integrale3}, it follows that
\begin{align*}
I_{a,1}& = \Big(\dfrac{\pi}{a} \Big)^{\frac{d}{2}} \E\Big[X_1^{\intercal} X_2 \exp\Big(-\dfrac{\Vert X_1-X_2 \Vert^2 }{4a}\Big)\Big], \\
I_{a,2} &= 2\Big(\dfrac{2\pi}{2a+1} \Big)^{\frac{d}{2}} \E\Big[ \dfrac{\Vert X_1 \Vert^2}{2a+1} \exp\Big(-\dfrac{\Vert X_1 \Vert^2}{4a+2}\Big)\Big], \qquad
I_{a,3} \Big(\dfrac{\pi}{a+1} \Big)^{\frac{d}{2}} \dfrac{d}{2a+2} .
\end{align*}
\noindent
The Taylor expansions \eqref{taylor1} und \eqref{taylor2}, together with $\E[X] = 0$, $\E[XX^{\intercal}]=\ID$ and  $\E\Vert X \Vert^6<\infty$ then yield
\begin{align*}
a^2\Big(\dfrac{a}{\pi} \Big)^{\frac{d}{2}}I_{a,1}&= a^2\E[X_1^{\intercal} X_2 ] - a\E\Big[X_1^{\intercal} X_2 \dfrac{\Vert X_1-X_2 \Vert^2 }{4} \Big] + \E \Big[X_1^{\intercal} X_2 \dfrac{\Vert X_1-X_2 \Vert^4 }{32} \Big ] + O(a^{-1}) \\
&= \dfrac{ad}{2}+ \dfrac{1}{16} \E [ X_1^{\intercal} X_2\Vert X_1 \Vert^2 \Vert X_2 \Vert^2 ] + \dfrac{2}{16} \E [(X_1^{\intercal} X_2)^3 ] - \dfrac{1}{4} \E \Vert X_1 \Vert^4  + O(a^{-1}),\\
a^2\Big(\dfrac{a}{\pi} \Big)^{\frac{d}{2}}I_{a,2}&= a\Big(\dfrac{a}{a+\frac{1}{2}} \Big)^{\frac{d}{2}+1} \E\Big[ \Vert X_1 \Vert^2 \exp\Big(-\dfrac{\Vert X_1 \Vert^2}{4a+2}\Big)\Big]= \Big(ad-\dfrac{d^2}{4}-\dfrac{d}{2}\Big)- \dfrac{1}{4}\E\Vert X_1 \Vert^4+O(a^{-1}), \\
a^2\Big(\dfrac{a}{\pi} \Big)^{\frac{d}{2}}I_{a,3}&=a\Big(\dfrac{a}{a+1} \Big)^{\frac{d}{2}+1} \dfrac{d}{2} =
\dfrac{ad}{2}-\dfrac{d^2}{4}-\dfrac{d}{2}+O(a^{-1}).
\end{align*}
Upon summarizing, the assertion follows. The second statement is proved following similar arguments.\end{prf}

We remark in passing that the first term on the right hand side of \eqref{limitdeltaainfty} is the population measure of
multivariate skewness in the sense of M\'{o}ri, Rohatgi, and Sz\'{e}kely \cite{MRS:94}, and $\E [(X_1^{\intercal}X_2)^3]$ is population
skewness in the sense of Mardia \cite{M:70}. Thus, Theorem \ref{thmdeltatoinfty} can be regarded as the 'population counterpart' of
Theorems  \ref{thmainfty} and \ref{thmato0}.

\cite{BEH:2017} observed that, in the context of goodness-of-fit testing of a general parametric hypothesis $\widetilde{H}_0$ (say),
weighted $L^2$-statistics have a normal limit under fixed alternatives to $\widetilde{H}_0$. To state such a theorem in our case, we first
introduce some notation. Again, we write  $\psi_X(t) = \E[\exp(\ii t^{\intercal}X)]$ for the CF of $X$ and put
$\psi_X^\pm(t) :=\ \text{Re }\psi_X(t) \pm \text{Im }\psi_X(t)$,
\begin{equation}\label{defwsx}
w(t,X) = X\CSP(t,X)-X\psi_X^+(t)-t^{\intercal}X \nabla \psi_X^+(t)+ \frac{1}{2}\big((XX^{\intercal}+\ID)\nabla \psi_X^-(t) - \E[XX^{\intercal} \CSM(t,X)] (XX^{\intercal}-\ID)t\big).
\end{equation}
Moreover, let
\begin{equation}\label{LST}
L(s,t) := \E\big{[}w(s,X)w(t,X)^{\intercal}\big{]}, \quad s,t \in \R^d.
\end{equation}

We then have the following result.

\begin{theorem}\label{thmnormallimit} If $\E \Vert X \Vert^4<\infty$, we have
\begin{align*}
\sqrt{n}\Big(\dfrac{T_{n,a}}{n}-\Delta_a\Big)\overset{\mathcal{D}}{\longrightarrow} {{\rm N}}(0,\sigma^2_a),
\end{align*}
where
\begin{align}\label{sigmaa}
\sigma_a^2:=4\iint z(s)^{\intercal}L(s,t)z(t)w_a(s)w_a(t) \, {\rm d}s\, {\rm d}t.
\end{align}
Here,
\begin{equation}\label{defzvont}
z(t):= \mu(t)-t\psi(t),
\end{equation}
and $L(s,t)$ is defined in \eqref{LST}.
\end{theorem}

\begin{prf} The basic observation is that, with $Z_n$ defined in \eqref{Zn} and $z(t):=\mu(t)-t\psi(t)$, we have
\begin{align}\nonumber
\sqrt{n}\Big(\dfrac{T_{n,a}}{n}-\Delta_a\Big)&=\sqrt{n}\big(\Vert n^{-1/2} Z_n \Vert_{\HH}^2 - \Vert z \Vert_{\HH}^2\big)
 = \sqrt{n} \langle n^{-1/2} Z_n -z, 2z+ n^{-1/2} Z_n -z \rangle_{\HH} \\ \label{secrow}
&=2\langle Z_n- \sqrt{n}z,z \rangle_{\HH} + n^{-1/2} \Vert Z_n - \sqrt{n}z \Vert_{\HH}^2.
\end{align}
Letting $V_n(t) := Z_n(t) - \sqrt{n}z(t) = n^{-1/2} \sum_{j=1}^n \big(Y_{n,j}\CSP(t,Y_{n,j})-\mu(t)\big)$,
the next step is to show that
\begin{align} \label{VnAsym}
V_n  \overset{\mathcal{D}}{\longrightarrow} V \text{ in } \HH
\end{align}
for some centred Gaussian random element $V$ of $\HH$ having covariance matrix kernel $L(s,t)$ given in \eqref{LST}.
The proof of \eqref{VnAsym} is completely analogous to that of Theorem \ref{thmlimitnull} and is therefore omitted.
In view of \eqref{VnAsym}, the second summand in \eqref{secrow} is $o_\PP(1)$, and the  first converges in
distribution to $2\langle V,z\rangle_\HH$ by the continuous mapping theorem. The distribution of
$2\langle V,z\rangle_\HH$ is the normal distribution ${\text N}(0,\sigma^2_a)$. \end{prf}

Using Slutsky's lemma, Theorem \ref{thmnormallimit} yields the following asymptotic confidence interval for $\Delta_a$.

\begin{corollary}\label{corconfint} For $\alpha \in (0,1)$, let $z_{1-\alpha/2}$ denote  the $(1-\alpha/2)$-quantile of the standard normal distribution.
If $\widehat{\sigma}^2_{n,a}$ is a consistent sequence of estimators for $\sigma^2_a$, and if $\sigma_a^2 >0$, then
\begin{align*}
I_{n,a,\alpha}:=\Big[\dfrac{T_{n,a}}{n}-\dfrac{\widehat{\sigma}_{n,a}}{\sqrt{n}}z_{1-\alpha/2},\dfrac{T_{n,a}}{n}+\dfrac{\widehat{\sigma}_{n,a}}{\sqrt{n}}z_{1-\alpha/2}\Big]
\end{align*}
is an asymptotic confidence interval with level $1-\alpha$ for $\Delta_a$.
\end{corollary}

A necessary and sufficient condition for $\sigma_a^2 >0$ is that the function $\R^d \ni s \mapsto \int L(s,t)z(t)w_a(t)\, {\text d}t$
does not vanish $\lambda^d$-almost everywhere, see Remark 1 of \cite{BEH:2017}.

To construct a consistent sequence of estimators for $\sigma^2_a$, we replace $z(s)$, $z(t)$ and $L(s,t)$ figuring in
\eqref{sigmaa} with suitable empirical counterparts. In view of \eqref{LST} and \eqref{defwsx} and the fact that
$\nabla \psi_X^+(t)=  \E[X{\rm CS}^-(t,X)]$, $\nabla \psi_X^-(t)= -\E[X{\rm CS}^+(t,X)]$, let
\begin{align}\label{Ln}
L_n(s,t):=\dfrac{1}{n} \sum_{j=1}^n W_{n,j}(s) W_{n,j}(t)^{\intercal},
\end{align}
where
\begin{equation}\label{defwnjt}
W_{n,j}(t) := Y_{n,j}\CSP(t,Y_{n,j})-Y_{n,j}\Psi_{1,n}(t)-t^{\intercal}Y_{n,j}\Psi_{2,n}(t)- \textstyle{\frac{1}{2}} (Y_{n,j}Y_{n,j}^{\intercal} + \ID)\Psi_{3,n}(t)- \textstyle{\frac{1}{2}}\Psi_{4,n}(t)(Y_{n,j}Y_{n,j}^{\intercal}-\ID)t,
\end{equation}
and
\begin{align}\label{defpsi12}
\Psi_{1,n}(t) &:=\dfrac{1}{n}\sum_{j=1}^n \CSP(t,Y_{n,j}), \quad \Psi_{2,n}(t):=\dfrac{1}{n}\sum_{j=1}^n Y_{n,j}\CSM(t,Y_{n,j}),\\ \label{defpsi34}
\Psi_{3,n}(t) & :=\dfrac{1}{n}\sum_{j=1}^n Y_{n,j}\CSP(t,Y_{n,j}), \quad
\Psi_{4,n}(t):=\dfrac{1}{n}\sum_{j=1}^n Y_{n,j}Y_{n,j}^{\intercal}\CSM(t,Y_{n,j}).
\end{align}
Furthermore, let
\begin{align}\label{zn}
z_n(t):=\dfrac{1}{n} \sum_{j=1}^n Y_{n,j} \CSP(t,Y_{n,j})-t\psi(t).
\end{align}

We then have the following result.

\begin{theorem}\label{thmconstest}  Let
\begin{equation*}
\widehat{\sigma}^2_{n,a}:=4\iint z_n(s)^{\intercal}L_n(s,t)z_n(t)w_a(s)w_a(t) \ {\rm d}s \ {\rm d}t,
\end{equation*}
where $L_n(s,t)$ and $z_n(t)$ are defined in \eqref{Ln} and \eqref{zn}, respectively. If $\E \|X\|^4 < \infty$, then
$(\widehat{\sigma}^2_{n,a})$  is a consistent sequence of estimators for $\sigma^2_a$, i.e., we have
$\widehat{\sigma}^2_{n,a} \overset{\mathbb{P}}{\longrightarrow} \sigma^2_a$. Moreover,
\begin{align} \label{intfrei}
\widehat{\sigma}^2_{n,a}=\sum_{i,j=1}^5 \widehat{\sigma}^{i,j}_{n,a},
\end{align}
where $\widehat{\sigma}^{i,j}_{n,a}$ is given in \eqref{defsigij}.
\end{theorem}

Since the proof of Theorem \ref{thmconstest} is long and tedious, it is deferred to Appendix \ref{secproofs}.
We stress that the representation \eqref{intfrei} does not comprise any integral, which means that
$\widehat{\sigma}^2_{n,a}$ is a feasible estimator.

We close this section with an example that illustrates the feasibility of the asymptotic confidence interval.
To this end, we consider the following  standardized symmetric alternatives to normality.
Firstly, let $X \edist {\text U}(-\sqrt{3},\sqrt{3})^d$ have the uniform distribution on the cube $(-\sqrt{3},\sqrt{3})^d$. In this case, we have
\[
\varphi_X(t)=\prod_{i=1}^d \dfrac{\sin(\sqrt{3}t_i)}{\sqrt{3}t_i}, \qquad  \nabla \varphi_X(t)^{(j)}=\dfrac{3\cos(\sqrt{3}t_j)t_j-\sqrt{3}\sin(\sqrt{3}t_j)}{3t_j^2}\prod_{i\neq j}^d \dfrac{\sin(\sqrt{3}t_i)}{\sqrt{3}t_i},
\]
where $\nabla \varphi(t)^{(j)}$ is the $j$th component of $\nabla \varphi(t)$. Secondly, we consider a Laplace distribution with i.i.d.
marginals, denoted by Laplace$(0,1/\sqrt{2})^d$, for which
\[
\varphi_X(t)=\prod_{i=1}^d \dfrac{2}{2+t_i^2}, \qquad  \nabla \varphi_X(t)^{(j)}=-\dfrac{4t_j}{(2+t_j^2)^2}\prod_{i\neq j}^d  \dfrac{2}{2+t_i^2}.
 \]
Finally, let $X$ have a logistic distribution with i.i.d. marginals, denoted by Logistic$(0,3/\pi)^d$. In this case, we obtain
\[
\varphi_X(t)=\prod_{i=1}^d \dfrac{\sqrt{3}t_i}{\sinh(\sqrt{3}t_i)}, \qquad  \nabla \varphi_X(t)^{(j)}=\dfrac{\sqrt{3}\sinh(\sqrt{3}t_j)-3t_j\cosh(\sqrt{3}t_j)}{\sinh(\sqrt{3}t_j)^2}\prod_{i\neq j}^d \dfrac{\sqrt{3}t_i}{\sinh(\sqrt{3}t_i)}.
\]
In each case, $\Delta_a$ has been computed by numerical integration. The resulting values are displayed in Table \ref{deltatab}.
\begin{table}[t]
\centering
\begin{tabular}{ rr|rrrr }
\hline\hline
& $d \setminus a$ & 0.5 & 1 & 2 & 5 \\ \hline
\multirow{2}{*}{U$(-\sqrt{3},\sqrt{3})^d$} & 1 & 0.029273 & 0.011432 &  0.002911 & 0.000259  \\
& 2 & 0.090821 & 0.027841 & 0.005709 & 0.000365  \\
\hline
\multirow{2}{*}{Laplace$(0,1/\sqrt{2})^d$} & 1 & 0.026076 & 0.013968 & 0.005230 & 0.000778  \\
& 2 & 0.071014 & 0.032525 & 0.010141 & 0.001097 \\
\hline
\multirow{2}{*}{Logistic$(0,\sqrt{3}/\pi)^d$} & 1 & 0.005014 & 0.002688 & 0.001005 & 0.000144  \\
& 2 & 0.013664 & 0.006226 & 0.001942 & 0.000202 \\
\hline\hline
\end{tabular}
\smallskip
\caption{Values of $\Delta_a$}
\label{deltatab}
\end{table}

By means of a Monte Carlo study, we estimated the probability of coverage of the confidence interval
$I_{n,a,\alpha}$ figuring in corollary  \ref{corconfint} for $a \in \{0.5,1,2,5\}$, $d \in \{1,2\}$,  and the  sample sizes
$n \in \{10,20,30,50,100,200\}$. The nominal level is $0.95$, and the number of replications is $10000$. Simulations have been carried out with the statistic software  R, see \cite{R:2019}. In particular, we used
  the package \texttt{extraDistr}, see \cite{W:19}, to generate variates from the Laplace distribution. The results are displayed in Table \ref{konfi}. As one can see the empirical coverage is converging to the nominal level, while it is obviously slower in higher dimensions. For larger values of the tuning parameter $a$ the confidence interval tends to be too wide, so we conjecture that an improvement of the asymptotic interval can be found.
\begin{table}[t]
\centering
\begin{tabular}{ r|rr|rrrr }
\hline\hline
& $d$ & $n \setminus a$ & 0.5 & 1 & 2 & 5  \\ \hline
\multirow{12}{*}{U$(-\sqrt{3},\sqrt{3})^d$} & \multirow{6}{*}{1} & 10 & 93.63 & 96.85 & 97.95 & 99.19 \\
& & 20 & 94.63 & 96.45 & 97.76 & 99.21 \\
& & 30 & 94.75 & 95.99 & 98.30 & 99.20 \\
& & 50 & 94.56 & 95.99 & 97.67 & 99.31 \\
& & 100 & 94.53 & 95.35 & 97.38 & 98.87 \\
& & 200 & 94.52 & 95.00 & 96.67 & 98.52 \\ \cline{2-7}
& \multirow{6}{*}{2} & 10 & 34.60 & 61.41 & 76.38 & 89.37 \\
& & 20 & 65.42 & 79.95 & 90.07 & 95.36 \\
& & 30 & 75.86 & 84.20 & 92.62 & 96.86 \\
& & 50 & 82.40 & 87.88 & 93.78 & 98.13 \\
& & 100 & 87.23 & 89.54 & 93.42 & 98.47 \\
& & 200 & 90.82 & 90.91 & 92.30 & 97.74 \\
\hline
\multirow{12}{*}{Laplace$(0,1/\sqrt{2})^d$} & \multirow{6}{*}{1} & 10 & 90.19 & 90.22 & 92.60 & 92.19 \\
& & 20 & 94.57 & 90.10 & 86.57 & 85.00 \\
& & 30 & 93.46 & 89.83 & 85.96 & 85.94 \\
& & 50 & 93.42 & 89.81 & 87.65 & 87.05 \\
& & 100 & 94.38 & 91.21 & 88.89 & 88.61 \\
& & 200 & 94.54 & 92.62 & 90.75 & 89.29 \\  \cline{2-7}
& \multirow{6}{*}{2} & 10 & 18.61 & 46.33 & 64.25 & 86.19 \\
& & 20 & 60.43 & 78.48 & 81.34 & 80.55 \\
& & 30 & 75.27 & 89.04 & 91.34 & 88.04 \\
& & 50 & 85.17 & 94.23 & 96.03 & 95.15 \\
& & 100 & 90.56 & 96.70 & 97.16 & 97.54 \\
& & 200 & 93.20 & 97.25 & 97.58 & 97.63 \\
\hline
\multirow{12}{*}{Logistic$(0,\sqrt{3}/\pi)^d$} & \multirow{6}{*}{1} & 10 & 91.92 & 94.40 & 96.85 & 99.29 \\
& & 20 & 97.82 & 98.24 & 98.01 & 97.27 \\
& & 30 & 98.91 & 98.66 & 97.84 & 97.20 \\
& & 50 & 99.13 & 98.22 & 96.78 & 96.80 \\
& & 100 & 98.26 & 95.90 & 94.82 & 95.04 \\
& & 200 & 96.49 & 94.53 & 93.20 & 94.42 \\ \cline{2-7}
& \multirow{6}{*}{2} & 10 & 12.13 & 37.85 & 56.30 & 77.48 \\
& & 20 & 43.14 & 69.10 & 79.29 & 81.10 \\
& & 30 & 56.73 & 81.06 & 89.30 & 90.50 \\
& & 50 & 69.40 & 89.60 & 94.57 & 95.62 \\
& & 100 & 80.17 & 94.59 & 98.13 & 98.63 \\
& & 200 & 85.47 & 96.55 & 98.85 & 99.39 \\
\hline\hline
\end{tabular}
\smallskip
\caption{Empirical coverage probability of  $I_{n,a,0.95}$ for $\Delta_a$ ($10000$ replications, nominal level 0.95)}
\label{konfi}
\end{table}

\section{Simulations}\label{secsimulation}
This section presents the results of a Monte Carlo study, with the aim to compare the power of the proposed test with respect to that of prominent
competitors against selected alternatives. We used the statistic software \texttt{R}, see \cite{R:2019}, and we employed the package \texttt{MonteCarlo}, see \cite{L:19}, which allows
 for parallel computing. In addition, we used the package \texttt{expm}, see \cite{GDMFSS:19}, for the standardization of the data.
Critical values for the test statistic have been  estimated by means of extensive simulations (100000 replications), and they are displayed in
 Table \ref{crit} for the weight parameters  $a \in \{0.5,1,2,5,10,\infty\}$ and the sample sizes $n \in \{20,50,100\}$. Throughout, the level of
 significance is $\alpha=0.05$. For the sake of comparison, Table \ref{crit} displays the approximate critical values of $T_{\infty,a}$ in the special case
 $d=1$, which have been obtained in Section \ref{secnull} by choosing a distribution of the Pearson family by equating the first four moments.
 As already mentioned in Section \ref{secbasprop}, the test statistic $T_{n,\infty}$ is a linear combination of skewness in the sense of Mardia
  \cite{M:70} and skewness in the sense of M\'{o}ri, Rohatgi und Sz\'{e}kely \cite{MRS:94}, and it equals the
  statistic $HV_{\infty}$ of Henze--Visagie, see \cite{HV:19}.

\begin{table}[t]
\centering
\begin{tabular}{ rr|rrrrrr }
\hline\hline
$d$ & $n/a$ & 0.5 & 1 & 2 & 5 & 10 & $\infty$  \\
\hline
\multirow{ 4 }{*}{ 1 } & 20 & 2.57 & 7.12 & 15.90 & 30.72 & 39.98 & 53.38 \\
 & 50 & 2.64 & 7.42 & 16.82 & 34.00 & 45.48 & 62.93 \\
 & 100 & 2.65 & 7.46 & 17.08 & 34.88 & 47.28 & 65.19 \\
 & $\infty$ & 2.67 & 7.52 & 17.28 & 35.56 & 46.23 & - \\
\hline
\multirow{ 3 }{*}{ 2 } & 20 & 5.77 & 15.94 & 35.47 & 70.27 & 93.10 & 125.90 \\
 & 50 & 5.83 & 16.27 & 37.16 & 76.41 & 102.65 & 145.38 \\
 & 100 & 5.87 & 16.19 & 37.35 & 77.40 & 106.51 & 151.15 \\
\hline
\multirow{ 3 }{*}{ 3 } & 20 & 9.43 & 27.03 & 61.74 & 125.52 & 167.47 & 230.75 \\
 & 50 & 9.57 & 27.37 & 64.02 & 135.16 & 186.80 & 267.89 \\
 & 100 & 9.58 & 27.47 & 64.38 & 137.79 & 190.30 & 276.76 \\
\hline
\multirow{ 3 }{*}{ 5 } & 20 & 17.89 & 55.55 & 137.20 & 296.36 & 407.65 & 581.08 \\
 & 50 & 18.03 & 56.21 & 141.10 & 319.59 & 452.61 & 681.00 \\
 & 100 & 18.05 & 56.32 & 141.21 & 323.19 & 462.59 & 704.12 \\
\hline\hline
\end{tabular}
\smallskip
\caption{Empirical $0.95$-quantiles for $a^{d/2+2}\pi^{-d/2}16T_{n,a}$ under $H_0$ (100000 replications) }
\label{crit}
\end{table}

\subsection{Univariate normal distribution}
In the univariate case $d=1$, we compared the power of our novel test statistics
with several competitors, which are
\begin{itemize}
\item the Cram\'{e}r--von Mises test (CvM),
\item the Anderson--Darling test (AD),
\item the  Shapiro--Wilk test (SW),
\item the Baringhaus--Henze--Epps--Pulley test (BHEP),
\item the Henze--Visagie test (HV).
\end{itemize}

The first three of these tests are well-known. The CvM-test and the AD-test have been implemented with the
        \texttt{R}-package \texttt{nortest}, see \cite{GL:15}, which contains the functions \texttt{cvm.test} and \texttt{ad.test}, and for the SW-test we used the
         function \texttt{shapiro.test} of the \texttt{stats}-package. The test statistics BHEP and HV will be explained in
        \eqref{BHEPT} and \eqref{HVT}, respectively.

For the BHEP-test and the HV-test, critical values have been simulated with $100 000$ replications.
These values and those of Table \ref{crit} for the novel test statistics
have been employed to assess the power of the various tests against several alternatives.
Table \ref{d1} exhibits percentages of rejection based on 100000 replications. An asterisk denotes power of 100\% and the best performing test for each alternative is marked in boldface.
The choice of alternatives orients itself towards those used in \cite{HV:19}.
The acronym NMix1 denotes a mixture of the normal distributions N$(0,1)$ and N$(3,1)$ with weights $0.9$ and $0.1$, respectively.

The novel tests outperform the selected competitors for the t$_3$-distribution, the  $\chi^2(15)$-distribution and the distribution NMix1,
and they keep up with the other procedures against the remaining alternatives. For most of the alternatives, power does not change much with
varying the weight parameter $a$. A notable exception is the uniform distribution U$(-\sqrt{3},\sqrt{3})$, against which power breaks down for larger tuning parameters, a feature
shared by the HV-test.

\renewcommand{\arraystretch}{1.0}
\begin{table}[p!]
\centering
\begin{tabular}{ rr|rrrrr|rrrrrr }
\hline\hline
 & $n$ & CvM & AD & SW & BHEP$_1$ & HV$_5$ & $T_{0.5}$ & $T_1$ & $T_2$ & $T_5$ & $T_{10} $ & $T_{\infty}$\\
\hline
\multirow{3}{*}{N$(0,1)$} & 20 & 5 & 5 & 5 & 5 & 5 &5 & 5 & 5 & 5 & 5& 5\\& 50 & 5 & 5 & 5 & 5 & 5 & 5 & 5 & 5 & 5 & 5& 5\\& 100 & 5 & 5 & 5 & 5 & 5 &5 & 5 & 5 & 5 & 5& 5\\ \hline
\multirow{3}{*}{NMix1} & 20 & 20 & 23 & 25 & 26 & 25 & 27 & \textbf{28} & \textbf{28} & \textbf{28} &27& 27\\& 50 & 45 & 50 & 56 & 55 & 52 & 58 & 60 & \textbf{61} & \textbf{61} & 60 & 59\\ & 100 & 75 & 81 & 85 & 84 & 82 &87 & 88 & \textbf{89}  & \textbf{89} & 88 & 88\\ \hline
\multirow{3}{*}{t$_3(0,1)$} & 20 & 30 & 33 & 34 & 33 & \textbf{36}&\textbf{36} & \textbf{36} & 35 & 35 &34& 35\\ & 50 & 57 & 61 & 64 & 61 & 63 &\textbf{66} & 65 & 63 & 59 &56& 52\\ & 100 & 83 & 85  &\textbf{88} & 86 & 84 &\textbf{88} & \textbf{88} & 86 & 80 &76& 64\\ \hline
\multirow{3}{*}{t$_5(0,1)$} & 20 & 15 & 17 & 19 & 18& \textbf{22} &20 & 20 & 20 & 20 &20 & 20\\ & 50 & 27 & 30 & 35 & 31 & \textbf{39} & 36 & 36 & 35 & 34 & 33& 32\\ & 100 & 43 & 48 & \textbf{57} & 50& 56 &55 & 56 & 53 & 49 & 45&40\\ \hline
\multirow{3}{*}{t$_{10}(0,1)$} & 20 & 8 & 9 & 10 & 9& \textbf{12} &11 & 11 & 11 & 11 & 11 &11 \\ & 50 & 11 & 12 & 15 & 13& \textbf{19} &15 & 16 & 16 & 16 & 16 & 16\\ & 100 & 14 & 16 & 23 & 17& \textbf{27}& 21 & 22 & 22 & 21 &20& 20\\ \hline
\multirow{3}{*}{$\chi^2(5)$} & 20 & 34 & 38 & \textbf{44} & 42 & 35 &42 & 43 & 43 & 42& 41 & 40\\ & 50 & 73 & 80 & \textbf{89} & 83& 74 &86 & 86 & 87 & 86 & 85 & 83\\ & 100 & 97 & 99 & \textbf{*} & 99& 97 &99 & 99 & \textbf{*} & 99 & 99 &99\\ \hline
\multirow{3}{*}{$\chi^2(15)$} & 20 & 14 & 15 & 17 & 17& 16 & 18 & \textbf{19} & \textbf{19} & \textbf{19} & \textbf{19}& 18\\
& 50 & 30 & 33 & 42 & 39& 37 &40 & 43 & \textbf{45} & \textbf{45} & \textbf{45}& 44\\ & 100 & 54 & 61 & 75 & 68 & 65 &71 & 74 & 76 & \textbf{77} & \textbf{77} &76\\ \hline
\multirow{3}{*}{Logistic$(0,1)$} & 20 & 10 & 11 & 11 & 11 & \textbf{14} &13 & 13 & 13 & 13 & 13 & 13\\
& 50 & 14 & 16 & 20 & 17& \textbf{23} &20 & 20 &20 & 19 &19& 19\\
& 100 & 21 & 24 & 31 & 25 & \textbf{32} &30 & 30 & 28 & 26& 24 & 23\\ \hline
\multirow{3}{*}{U$(-\sqrt{3},\sqrt{3})$} & 20 & 14 & 17  & \textbf{20} & 12& 0& 10& 4&2 &1 &1 &1\\ & 50 & 44 & 58 & \textbf{75} & 55 & 0& 55& 33& 5& 1& 0&0\\ & 100 & 84 & 95 & \textbf{*} &94 & 0&96 &90 & 48& 2& 1&0\\ \hline
\multirow{3}{*}{P$_{VII}(5)$} & 20 & 15 & 17 & 19 & 18& \textbf{22} & 20 & 20 & 20 & 20 & 20 & 21 \\
& 50 & 27 & 30 & 35 & 31& \textbf{39}& 36 & 36 & 35 & 34 & 33& 32\\
& 100 & 43 & 48 & \textbf{57}& 50 & 56 &55 & 56 & 53 & 49 & 45& 41\\ \hline
\multirow{3}{*}{P$_{VII}(10)$}& 20 & 8 & 9 & 10 & 9 & \textbf{12} &11 & 11 & 11 & 11 & 11 & 11\\
& 50 & 11 & 12 & 16 & 12 & \textbf{19} &15 & 16 & 16 & 16 & 16 & 16\\
& 100 & 14 & 16 & 23 & 17 & \textbf{27} &21 & 22 & 22 & 20 &20& 20\\
\hline\hline
\end{tabular}
\smallskip
\caption{Empirical power ($d=1$, $\alpha=0.05$, $100000$ replications)}
\label{d1}
\end{table}

\subsection{Multivariate normal distribution}
For the dimensions $d=2$, $d=3$ and $d=5$, we compared the novel test statistic with the following procedures:
\begin{itemize}
\item the test of Baringhaus--Henze--Epps--Pulley (BHEP),
\item the test of Henze--Zirkler (HZ),
\item the test of Henze--Visagie (HV),
\item the energy test (EN).
\end{itemize}
A recent synopsis of tests for multivariate normality is given in \cite{EH:20}.
Just as the novel procedure, the  BHEP-test (see \cite{HW:1997}) is based on the empirical characteristic function (ECF).
More precisely, it employs the test statistic
\begin{equation} \label{BHEPT}
{\text{BHEP}}_a = \int |\psi_n(t) - \psi(t)|^2 \varphi_a(t)\, {\text d}t,
\end{equation}
where $\varphi_a(t) = (2\pi a^2)^{-d/2} \exp(-\|t\|^2/(2a^2))$, and $\psi_n(t)$ and $\psi(t)$ are given in \eqref{defpsint} and \eqref{defcfn01}, respectively.
An alternative representation for BHEP$_a$ is
\begin{equation*}
{\text{BHEP}}_a= \dfrac{1}{n^2} \sum_{i,j=1}\exp\Big(-\dfrac{a^2}{2} \big\Vert Y_{n,i} - Y_{n,j}\big\Vert^2\Big)-2(1+a^2)^{-\frac{d}{2}} \dfrac{1}{n}\sum_{j=1}^n \exp\bigg(-\dfrac{a^2\big\Vert Y_{n,j}\big\Vert^2}{2(1+a^2)}\bigg)+(1+2a^2)^{-\frac{d}{2}}. \nonumber
\end{equation*}
In our study, we used the special value $a=1$.

The test HZ of Henze--Zirkler (cf. \cite{HZ:90}) originates if we choose $a=1/\sqrt{2}\left((2d+1)n/4\right)^\frac{1}{d+4}$ in the BHEP test. The \texttt{R}-package \texttt{HZ}, see \cite{KGZ:14}, contains the function \texttt{mvn}, which calculates the statistic of the HZ-test.

The recent test of Henze--Visagie, see \cite{HV:19}, is the 'moment generating function analog' of our novel test statistic.
It employs the test statistic
\[
{\text{HV}}_a = n \int \| \nabla M_n(t) - t M_n(t)\|^2 w_a(t)\, {\text d}t,
\]
where $M_n(t) = n^{-1}\sum_{j=1}^n \exp(t^{\intercal}Y_{n,j})$ is the empirical moment generating function of the scaled residuals.
An alternative representation of HV$_a$ is
\begin{align} \label{HVT}
{\text{HV}}_a=\dfrac{1}{n}  \Big(\dfrac{\pi}{a}\Big)^{\frac{d}{2}} \sum_{i,j=1}^n \exp\Big(\dfrac{\Vert Y_{n,i}+ Y_{n,j}\Vert^2}{4a}\Big) \Big(Y_{n,i}^{\intercal}Y_{n,j} + \Vert Y_{n,i} + Y_{n,j} \Vert^2 \Big(\dfrac{1}{4a^2}-\dfrac{1}{2a}\Big)+ \dfrac{d}{2a}\Big).
\end{align}
In our comparative study, we put $a=5$, as recommended in \cite{HV:19}.

The rationale of the energy test of Sz\'{e}kely and Rizzo, see \cite{SR:05}, is based on the fact that, if  $X$ and $Y$ are independent integrable $d$-dimensional random vectors and
 $X',Y'$ denote independent copies of $X$ and $Y$, respectively, then
\begin{align*}
2\E\Vert X-Y \Vert - \E \Vert X - X' \Vert - \E \Vert Y - Y' \Vert \geq 0.
\end{align*}
Here, equality holds if and only if $X \edist Y$. The statistic of the energy test for multivariate normality is
\begin{align*}
{\text{EN}} = n \Big(\dfrac{2}{n}\sum_{j=1}^n \E \Vert \widetilde{Y}_{n,j} - Z_1 \Vert - \E \Vert Z_1 - Z_2 \Vert -\dfrac{1}{n^2}\sum_{i,j=1}^n \E\Vert \widetilde{Y}_{n,i}-\widetilde{Y}_{n,j}\Vert \Big).
\end{align*}
Here, $\widetilde{Y}_{n,j}=\sqrt{n/(n-1)}Y_{n,j}$, and $Z_1, Z_2$ are i.i.d. with the normal distribution N$_d(0,\ID)$, which are also independent
of $Y_{n,1}, \ldots, Y_{n,n}$. To calculate EN, notice that
$\E \Vert Z_1 - Z_2 \Vert = 2\Gamma(\frac{d+1}{2})/\Gamma(\frac{d}{2})$ and
\[
\E\Vert a- Z \Vert = \sqrt{2} \dfrac{\Gamma(\frac{d+1}{2})}{\Gamma(\frac{d}{2})} +\sqrt{\dfrac{2}{\pi}} \sum_{k=0}^{\infty}\dfrac{(-1)^k}{k!2^k}\dfrac{\Vert a \Vert^{2k+2}}{(2k+1)(2k+2)}\dfrac{2\Gamma(\frac{d+1}{2})\Gamma(k+\frac{3}{2})}{\Gamma(k+\frac{d}{2}+1)}.
\]
The \texttt{R}-package \texttt{energy} \cite{RS:19} contains the function \texttt{mvnorm.etest} to calculate EN. Note that all of the mentioned procedures are also implemented in the \texttt{R}-package \texttt{mnt}, see \cite{BE:2020}.

Just as done in the case $d=1$, we first simulated critical values with 100000 replications.
With the same number of replications, we then simulated the power of the tests under discussion against selected alternatives.
Again, the choice of alternatives orients itself towards those used in \cite{HV:19}.
Tables \ref{d2}, \ref{d3} and \ref{d5} display percentages of rejection of $H_0$ for dimensions $d=2$, $d=3$ and $d=5$, respectively,
and an asterisk again denotes power 100\%. To generate pseudo random numbers, we used the  \texttt{R}-packages \texttt{mvtnorm}, see \cite{GBMMLST:19}, and \texttt{PearsonDS}, see \cite{BK:2017}.
Suppressing the dimension $d$, the distribution NMix1 is a mixture of the normal distributions N$_d(0,\ID)$ and N$_d(3,\ID)$ with mixing proportions $0.9$ and $0.1$, respectively.
Here, $3$ stands for the $d$-dimensional vector that contains $3$ in each component.
Likewise,  NMix2 denotes a mixture of the normal distributions N$_d(0,\ID)$ and N$_d(0,B_d)$ with mixing proportions 0.1 and 0.9, respectively. Here,
 $B_d$ is a $d\times d$-matrix with $1$ for each diagonal entry and  $0.9$ for each off-diagonal entry.

The novel tests outperform their competitors for some alternatives, notably for the $\chi^2$-, the $\Gamma$-, and the NMix-distribution, but they can also
keep up for the other alternatives. However, just as in the univariate case,
          power is extremely low against the uniform distribution U$(-\sqrt{3},\sqrt{3})$, a feature shared by the HV-test.
Based on the results of this simulation study, we recommend the choice $a=5$ for the tuning parameter, since it leads to
competitive power against nearly each of the alternatives considered.

\renewcommand{\arraystretch}{1.0}
\begin{table}[p!]
\centering
\begin{tabular}{ rr|rrrr|rrrrrr }
\hline\hline
 & $n$ & BHEP$_1$ & HZ & HV$_5$ & EN & $T_{0.5}$ & $T_1$ & $T_2$ & $T_5$ & $T_{10}$ & $T_{\infty}$\\
\hline
\multirow{3}{*}{N$_2(0,{\text I}_2)$} & 20 & 5 & 5 & 5 & 5 & 5 & 5 & 5 & 5 & 5 & 5\\& 50 & 5 & 5 & 5 & 5 & 5 & 5 & 5 & 5 & 5 &5\\& 100 & 5 & 5 & 5 & 5 & 5 & 5 & 5 & 5 & 5 &5\\ \hline
\multirow{3}{*}{NMix1} & 20 & 39&34&32&37&38&\textbf{41}&\textbf{41}&40&39&38\\& 50 & 83&74&68&82&85&88&\textbf{89}&88&88&86\\ & 100 & 99&96&97&99&99&99&\textbf{*}&\textbf{*}&\textbf{*}&\textbf{*}\\ \hline
\multirow{3}{*}{NMix2} & 20 & 20&17&\textbf{27}&20&23&24&25&25&25&25\\& 50 & 38&30&\textbf{53}&39&45&48&49&48&47&44\\ & 100 & 60&47&\textbf{77}&61&68&72&72&70&66&55\\ \hline
\multirow{3}{*}{t$_3(0,I_2)$} & 20 & 47&45&54&49&49&51&\textbf{53}&\textbf{53}&\textbf{53}&52\\& 50 & 83&80&\textbf{85}&84&82&84&83&83&81&78\\ & 100 & \textbf{98}&97&97&97&97&\textbf{98}&\textbf{98}&97&95&90\\ \hline
\multirow{3}{*}{t$_5(0,I_2)$} & 20 & 25&22&\textbf{32}&26&27&29&30&31&31&31\\& 50 & 49&42&\textbf{59}&50&49&53&55&54&54&52\\ & 100 & 75&67&\textbf{81}&76&71&76&77&75&72&66\\ \hline
\multirow{3}{*}{t$_{10}(0,I_2)$} & 20 & 11&10&\textbf{16}&12&12&14&14&15&15&\textbf{16}\\& 50 & 17&14&\textbf{29}&18&19&22&24&25&25&25\\ & 100 & 27&20&\textbf{43}&28&26&31&33&34&33&33\\ \hline
\multirow{3}{*}{$(\chi^2(5))^2$} & 20 & 48&44&38&46&46&48&\textbf{50}&48&47&46\\& 50 & 93&87&80&92&93&94&\textbf{95}&\textbf{95}&94&93\\ & 100 & \textbf{*}&\textbf{*}&99&\textbf{*}&\textbf{*}&\textbf{*}&\textbf{*}&\textbf{*}&\textbf{*}&\textbf{*}\\ \hline
\multirow{3}{*}{$(\chi^2(15))^2$} & 20 & 18&16&17&17&17&19&\textbf{20}&\textbf{20}&\textbf{20}&19\\& 50 & 45&35&39&42&43&49&53&\textbf{55}&54&52\\ & 100 & 78&62&71&77&78&84&88&\textbf{89}&88&88\\ \hline
\multirow{3}{*}{$(\chi^2(20))^2$} & 20 & 15&13&14&14&14&15&\textbf{16}&\textbf{16}&\textbf{16}&\textbf{16}\\& 50 & 34&27&31&33&33&38&41&\textbf{43}&\textbf{43}&42\\ & 100 & 64&47&58&63&64&71&76&\textbf{78}&77&77\\ \hline
\multirow{3}{*}{$\Gamma(5,1)^2$} & 20 & 26&23&23&24&24&27&\textbf{28}&27&27&26\\& 50 & 64&53&53&61&62&68&71&\textbf{72}&71&69\\ & 100 & 93&84&87&93&94&96&97&\textbf{98}&97&97\\ \hline
\multirow{3}{*}{$\Gamma(4,2)^2$} & 20 & 32&28&27&30&30&33&\textbf{34}&33&33&32\\& 50 & 75&64&61&73&73&79&\textbf{81}&\textbf{81}&80&79\\ & 100 & 98&92&93&97&98&\textbf{99}&\textbf{99}&\textbf{99}&\textbf{99}&\textbf{99}\\ \hline
\multirow{3}{*}{Logistic$(0,1)^2$} & 20 & 11&10&\textbf{16}&12&13&14&15&\textbf{16}&15&\textbf{16}\\& 50 & 18&15&\textbf{29}&20&20&23&24&25&25&25\\ & 100 & 29&23&\textbf{42}&31&29&34&35&34&33&31\\ \hline
\multirow{3}{*}{U$(-\sqrt{3},\sqrt{3})^2$} & 20 & 12&\textbf{18}&0&11&6&3&1&1&0&0\\& 50 & 60&\textbf{67}&0&52&32&13&3&0&0&0\\ & 100 & \textbf{98}&\textbf{98}&0&96&92&80&24&1&0&0\\ \hline
\multirow{3}{*}{P$_{VII}(5)^2$} & 20 & 20&18&\textbf{28}&21&22&24&26&26&26&27\\& 50 & 39&32&\textbf{51}&40&41&45&47&46&46&45\\ & 100 & 63&53&\textbf{73}&64&62&67&68&66&62&58\\ \hline
\multirow{3}{*}{P$_{VII}(10)^2$} & 20 & 10&8&\textbf{13}&10&11&11&12&\textbf{13}&\textbf{13}&\textbf{13}\\& 50 & 13&11&\textbf{23}&14&15&18&19&20&20&20\\ & 100 & 19&14&\textbf{35}&21&20&24&26&27&27&26\\ \hline
\multirow{3}{*}{P$_{VII}(20)^2$} & 20 & 7&6&\textbf{8}&7&7&7&7&\textbf{8}&\textbf{8}&\textbf{8}\\& 50 & 7&7&\textbf{12}&8&8&9&10&11&11&11\\ & 100 & 8&7&\textbf{17}&9&9&11&11&13&12&13\\
\hline\hline
\end{tabular}
\smallskip
\caption{Empirical power ($d=2$, $\alpha=0.05$, $100000$ replications)}
\label{d2}
\end{table}

\renewcommand{\arraystretch}{1.0}
\begin{table}[p!]
\centering
\begin{tabular}{ rr|rrrr|rrrrrr }
\hline\hline
 & $n$ & BHEP$_1$ & HZ & HV$_5$ & EN & $T_{0.5}$ & $T_1$ & $T_2$ & $T_5$ & $T_{10}$ & $T_{\infty}$\\
\hline
\multirow{3}{*}{N$_3(0,{\text I}_3)$} & 20 & 5 & 5 & 5 & 5 & 5 & 5 & 5 & 5 &5&5\\& 50 & 5 & 5 & 5 & 5 & 5 & 5 & 5 & 5&5 &5\\& 100 & 5 & 5 & 5 & 5 & 5 & 5 & 5 & 5&5 &5\\ \hline
\multirow{3}{*}{NMix1} & 20 & 39&35&33&41&40&43&\textbf{44}&43&41&40\\& 50 & 89&81&66&91&91&94&\textbf{95}&\textbf{95}&93&92\\ & 100 & \textbf{*}&98&95&\textbf{*}&\textbf{*}&\textbf{*}&\textbf{*}&\textbf{*}&\textbf{*}&\textbf{*}\\ \hline
\multirow{3}{*}{NMix2} & 20 & 28&24&\textbf{43}&33&34&38&40&41&41&41\\& 50 & 59&49&\textbf{80}&66&65&72&75&75&75&73\\ & 100 & 85&74&\textbf{96}&88&87&92&93&94&92&87\\ \hline
\multirow{3}{*}{t$_3(0,I_3)$} & 20 & 56&53&65&62&58&63&65&\textbf{66}&65&65\\& 50 & 93&90&\textbf{94}&\textbf{94}&89&93&93&93&92&91\\ & 100 & \textbf{*}&\textbf{*}&\textbf{*}&98&99&\textbf{*}&\textbf{*}&99&99&98\\ \hline
\multirow{3}{*}{t$_5(0,I_3)$} & 20 & 29&26&\textbf{41}&35&32&37&39&\textbf{41}&40&\textbf{41}\\& 50 & 62&54&\textbf{73}&67&57&67&70&70&70&69\\ & 100 & 90&83&\textbf{92}&91&80&88&90&89&88&84\\ \hline
\multirow{3}{*}{t$_{10}(0,I_3)$} & 20 & 12&11&\textbf{20}&15&14&17&18&19&19&\textbf{20}\\& 50 & 22&17&\textbf{38}&26&22&28&32&34&35&35\\ & 100 & 37&28&\textbf{57}&42&30&40&46&48&48&47\\ \hline
\multirow{3}{*}{$(\chi^2(5))^3$} & 20 & 48&43&38&49&46&50&\textbf{51}&50&49&48\\& 50 & 95&89&82&96&94&\textbf{97}&\textbf{97}&\textbf{97}&\textbf{97}&96\\ & 100 & \textbf{*}&\textbf{*}&99&\textbf{*}&\textbf{*}&\textbf{*}&\textbf{*}&\textbf{*}&\textbf{*}&\textbf{*}\\ \hline
\multirow{3}{*}{$(\chi^2(15))^3$} & 20 & 17&15&17&18&16&18&\textbf{19}&\textbf{19}&\textbf{19}&\textbf{19}\\& 50 & 45&34&38&48&44&51&56&\textbf{58}&57&56\\ & 100 & 82&64&69&84&81&88&92&\textbf{93}&\textbf{93}&92\\ \hline
\multirow{3}{*}{$(\chi^2(20))^3$} & 20 & 13&12&14&14&13&14&\textbf{16}&15&15&15\\& 50 & 34&25&30&36&31&39&43&\textbf{45}&44&44\\ & 100 & 67&48&56&70&65&75&81&\textbf{83}&\textbf{83}&82\\ \hline
\multirow{3}{*}{$\Gamma(5,1)^3$} & 20 & 25&22&23&25&23&26&\textbf{28}&27&27&26\\& 50 & 65&53&53&68&64&71&\textbf{76}&\textbf{76}&75&74\\ & 100 & 96&86&87&97&96&98&\textbf{99}&\textbf{99}&\textbf{99}&\textbf{99}\\ \hline
\multirow{3}{*}{$\Gamma(4,2)^3$} & 20 & 30&27&27&32&29&32&\textbf{34}&33&33&32\\& 50 & 77&65&62&79&76&82&85&\textbf{86}&85&83\\ & 100 & 99&94&93&99&99&\textbf{*}&\textbf{*}&\textbf{*}&\textbf{*}&\textbf{*}\\ \hline
\multirow{3}{*}{Logistic$(0,1)^3$} & 20 & 11&10&\textbf{17}&13&13&15&16&\textbf{17}&\textbf{17}&\textbf{17}\\& 50 & 18&14&\textbf{32}&22&19&24&27&29&29&29\\ & 100 & 31&23&\textbf{48}&36&27&35&39&39&39&38\\ \hline
\multirow{3}{*}{U$(-\sqrt{3},\sqrt{3})^3$} & 20 & 11&\textbf{15}&0&6&5&2&1&0&0&0\\& 50 & 58&\textbf{65}&0&39&20&8&2&0&0&0\\ & 100 & \textbf{98}&\textbf{98}&0&94&79&51&12&1&0&0\\ \hline
\multirow{3}{*}{P$_{VII}(5)^3$} & 20 & 20&17&\textbf{30}&24&23&27&29&\textbf{30}&\textbf{30}&\textbf{30}\\& 50 & 41&34&\textbf{58}&47&42&50&54&55&54&53\\ & 100 & 69&57&\textbf{81}&73&63&72&76&75&73&69\\ \hline
\multirow{3}{*}{P$_{VII}(10)^3$} & 20 & 9&8&\textbf{14}&11&11&12&13&\textbf{14}&\textbf{14}&\textbf{14}\\& 50 & 13&10&\textbf{26}&16&14&18&21&23&23&23\\ & 100 & 20&14&\textbf{39}&24&18&24&29&31&31&31\\ \hline
\multirow{3}{*}{P$_{VII}(20)^3$} & 20 & 6&6&\textbf{9}&7&7&7&7&8&8&8\\& 50 & 7&6&\textbf{13}&8&8&9&10&11&12&12\\ & 100 & 8&7&\textbf{17}&10&8&10&12&13&14&14\\
\hline\hline
\end{tabular}
\smallskip
\caption{Empirical power ($d=3$, $\alpha=0.05$, $100000$ replications)}
\label{d3}
\end{table}

\renewcommand{\arraystretch}{1.0}
\begin{table}[p!]
\centering
\begin{tabular}{ rr|rrrr|rrrrrr }
\hline\hline
 & $n$ & BHEP$_1$ & HZ & HV$_5$ & EN  & $T_{0.5}$ & $T_1$ & $T_2$ & $T_5$ & $T_{10}$ & $T_{\infty}$\\
\hline
\multirow{3}{*}{N$_5(0,{\text I}_5)$} & 20 & 5 & 5 & 5 & 5 & 5 & 5 & 5 & 5 &5&5\\& 50 & 5 & 5 & 5 & 5 & 5 & 5 & 5 & 5&5 &5\\& 100 & 5 & 5 & 5 & 5 & 5 & 5 & 5 & 5&5 &5\\ \hline
\multirow{3}{*}{NMix1} & 20 & 25&22&31&32&27&33&\textbf{36}&34&34&33\\& 50 & 85&74&50&94&87&94&\textbf{95}&92&90&86\\ & 100 & \textbf{*}&98&77&\textbf{*}&\textbf{*}&\textbf{*}&\textbf{*}&\textbf{*}&\textbf{*}&\textbf{*}\\ \hline
\multirow{3}{*}{NMix2} & 20 & 32&27&\textbf{62}&48&40&51&56&58&59&59\\& 50 & 76&67&\textbf{96}&89&79&89&93&94&94&94\\ & 100 & 96&92&\textbf{*}&99&96&99&99&\textbf{*}&\textbf{*}&\textbf{*}\\ \hline
\multirow{3}{*}{t$_3(0,I_5)$} & 20 & 62&59&79&76&67&76&79&\textbf{81}&\textbf{81}&80\\& 50 & 98&97&99&99&99&\textbf{*}&\textbf{*}&99&99&99\\ & 100 & \textbf{*}&\textbf{*}&\textbf{*}&\textbf{*}&\textbf{*}&\textbf{*}&\textbf{*}&\textbf{*}&\textbf{*}&\textbf{*}\\ \hline
\multirow{3}{*}{t$_5(0,I_5)$} & 20 & 31&28&54&47&37&48&52&54&54&\textbf{55}\\& 50 & 77&71&\textbf{89}&88&68&82&88&\textbf{89}&\textbf{89}&\textbf{89}\\ & 100 & 98&96&\textbf{99}&\textbf{99}&88&96&98&\textbf{99}&98&98\\ \hline
\multirow{3}{*}{t$_{10}(0,I_5)$} & 20 & 12&11&\textbf{26}&20&15&21&24&25&\textbf{26}&\textbf{26}\\& 50 & 28&23&\textbf{55}&44&26&39&48&52&54&53\\ & 100 & 54&44&\textbf{78}&69&36&54&67&72&73&72\\ \hline
\multirow{3}{*}{$(\chi^2(5))^5$}& 20 & 39&35&36&\textbf{48}&39&46&\textbf{48}&\textbf{48}&47&45\\& 50 & 94&87&80&\textbf{98}&94&97&\textbf{98}&\textbf{98}&\textbf{98}&97\\ & 100 & \textbf{*}&\textbf{*}&99&\textbf{*}&\textbf{*}&\textbf{*}&\textbf{*}&\textbf{*}&\textbf{*}&\textbf{*}\\ \hline
\multirow{3}{*}{$(\chi^2(15))^5$} & 20 & 13&12&15&16&13&15&\textbf{17}&\textbf{17}&\textbf{17}&\textbf{17}\\& 50 & 38&29&35&52&37&49&56&\textbf{58}&\textbf{58}&56\\ & 100 & 78&60&64&90&77&89&94&\textbf{95}&\textbf{95}&94\\ \hline
\multirow{3}{*}{$(\chi^2(20))^5$} & 20 & 11&9&12&13&11&12&13&\textbf{14}&13&13\\& 50 & 28&22&28&39&27&36&42&\textbf{45}&44&43\\ & 100 & 61&43&51&77&60&74&83&\textbf{86}&\textbf{86}&85\\ \hline
\multirow{3}{*}{$\Gamma(5,1)^5$} & 20 & 18&16&21&24&18&22&24&\textbf{25}&24&24\\& 50 & 59&47&20&74&58&71&78&\textbf{79}&78&76\\ & 100 & 95&85&83&99&95&99&99&\textbf{*}&\textbf{*}&99\\ \hline
\multirow{3}{*}{$\Gamma(4,2)^5$} & 20 & 23&20&25&29&23&28&\textbf{30}&\textbf{30}&\textbf{30}&29\\& 50 & 72&60&59&84&71&83&87&\textbf{88}&87&85\\ & 100 & 99&94&91&\textbf{*}&99&\textbf{*}&\textbf{*}&\textbf{*}&\textbf{*}&\textbf{*}\\ \hline
\multirow{3}{*}{Logistic$(0,1)^5$} & 20 & 9&8&\textbf{17}&13&11&14&16&\textbf{17}&\textbf{17}&\textbf{17}\\& 50 & 15&13&\textbf{34}&26&17&24&30&33&\textbf{34}&\textbf{34}\\ & 100 & 29&22&\textbf{53}&42&23&34&43&47&47&47\\ \hline
\multirow{3}{*}{U$(-\sqrt{3},\sqrt{3})^5$} & 20 & 9&\textbf{11}&0&2&4&2&1&0&0&0\\& 50 & 50&\textbf{51}&0&12&12&4&1&0&0&0\\ & 100 & \textbf{96}&95&0&75&49&20&5&0&0&0\\ \hline
\multirow{3}{*}{P$_{VII}(5)^5$} & 20 & 16&14&\textbf{33}&25&20&27&30&\textbf{33}&32&32\\& 50 & 39&32&\textbf{67}&56&39&54&62&65&65&65\\ & 100 & 71&60&\textbf{89}&83&59&77&84&86&85&83\\ \hline
\multirow{3}{*}{P$_{VII}(10)^5$} & 20 & 8&7&\textbf{14}&11&9&11&13&\textbf{14}&\textbf{14}&\textbf{14}\\& 50 & 11&9&\textbf{28}&19&12&18&23&26&27&27\\ & 100 & 18&13&\textbf{44}&28&16&24&32&37&38&38\\ \hline
\multirow{3}{*}{P$_{VII}(20)^5$} & 20 & 6&5&\textbf{8}&7&7&7&\textbf{8}&\textbf{8}&\textbf{8}&\textbf{8}\\& 50 & 7&6&\textbf{13}&9&8&9&11&12&12&12\\ & 100 & 7&7&\textbf{19}&11&8&10&13&15&16&16\\
\hline\hline
\end{tabular}
\smallskip
\caption{Empirical power ($d=5$,  $\alpha=0.05$, $100000$ replications)}
\label{d5}
\end{table}

\section{A real data example}\label{secrealdata}
\begin{figure}[t]
\centering
\includegraphics[scale=0.5]{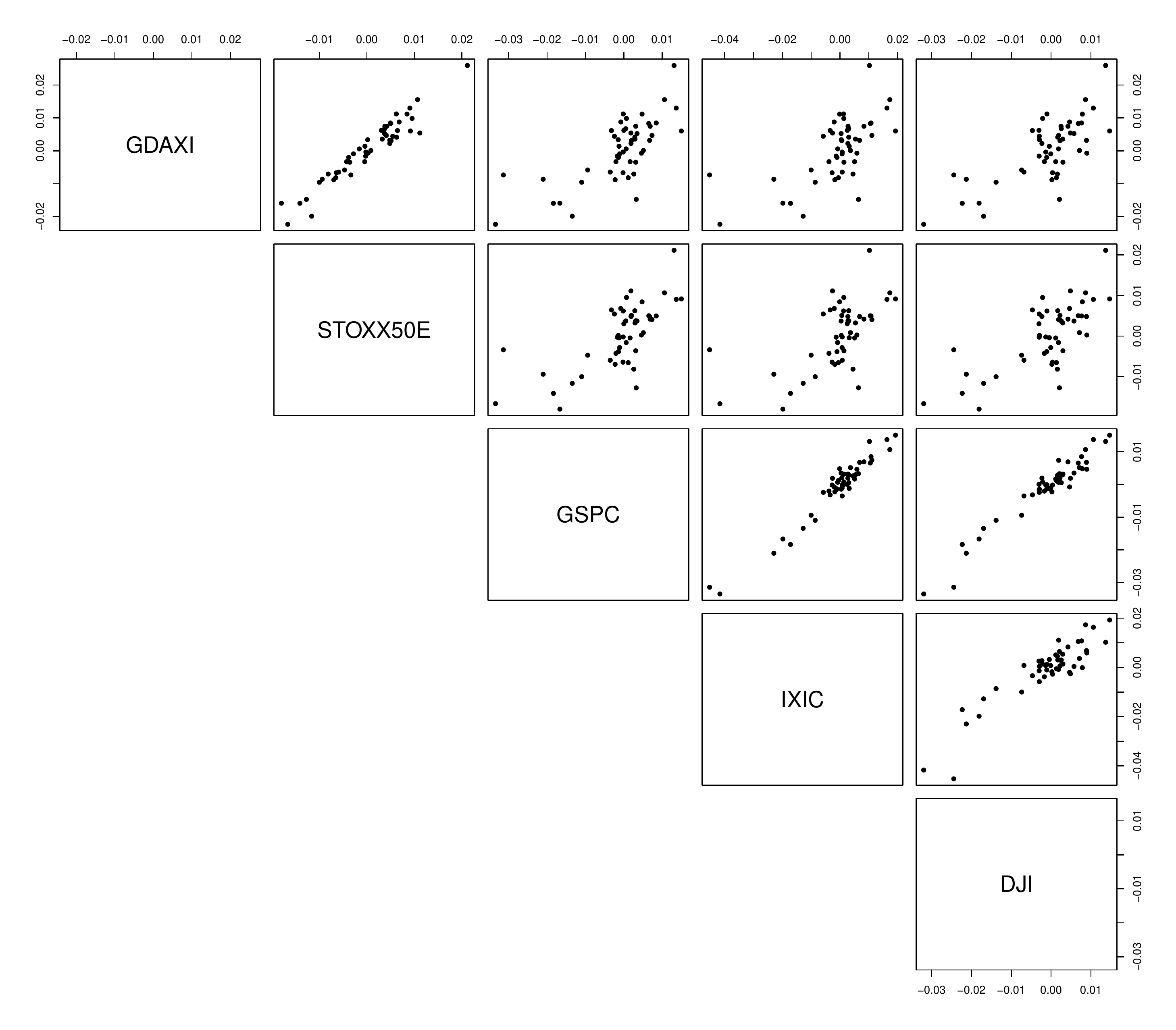}
\caption{2D projections of the log returns of the indexes .}\label{fig:data}
\end{figure}

The Black-Scholes-Merton model is a stochastic model for the dynamics of a financial market that contains derivative investment instruments.
One of the basic assumptions of this model is the normality  of the log returns of stocks
and indexes. To test the hypothesis joint normality of log returns of several indexes, we consider the five stock  indexes Standard \& Poor 500 (\textasciicircum GSPC), Dow Jones Industrial Average
(\textasciicircum DJI), NASDAQ Composite (\textasciicircum IXIC), DAX Perfomance Index (\textasciicircum GDAXI), and EURO STOXX 50 (\textasciicircum STOXX50E), over a period
 of  50 trading days, starting July 1st, 2017. The data (daily closing prices of the stocks) were obtained by means of the
 \texttt{R}-package \texttt{quantmod}, see \cite{RU:19}. To model the independence assumption between the realisations, we ignored a time span of 10 trading days between each of the five dimensional observations. Figure \ref{fig:data} shows a plot of the two-dimensional projections of the log returns.  For each value $a \in \{0.5, 1, 2, 5, 10\}$ of the weight parameter $a$, we performed a Monte Carlo simulation
 based on 100000  replications, in order to estimate the p-value of the observations.
The empirical p-values are displayed in Table \ref{pval}. As can be seen, the hypothesis of a multivariate
normality of the log returns of the selected stock prices is rejected at the 1\%-level, for each of the choices of the weight parameter $a$.

\renewcommand{\arraystretch}{1.2}
\begin{table}[t]
\centering
\begin{tabular}{ r|rrrrr }
\hline\hline
$a$ & 0.5 & 1 & 2 & 5 & 10 \\ \hline
p-value & 0.0002 & 0.0001 & 0.0002 & 0.0003 & 0.0003 \\
\hline\hline
\end{tabular}
\smallskip
\caption{Empirical p-value ($100000$ replications)}
\label{pval}
\end{table}

\section{Summary and outlook}\label{secsummary}
We propose a novel class of tests of normality based on an initial value problem connected to a multivariate Stein equation, which characterises the multivariate standard normal law. We derived asymptotic theory under the null hypothesis as well as under contiguous and fixed alternatives. Moreover, we proved consistency against each alternative distribution that satisfies a weak moment condition, and we provided insights into the structure of the behaviour of the test statistic under fixed alternatives by calculating asymptotic confidence intervals for $\Delta_a$, and by providing a consistent estimator for the limiting variance $\sigma_a^2$. Monte Carlo simulations show that the methods operate as expected, and that the new family of tests is a strong class of competitors to established procedures.

A first open question for further research is to find explicit formulae or numerical stable approximations for the eigenvalues $\lambda_j(a)$, $j=1,2,\ldots$ connected to the integral operator $\mathbb{K}$ in \eqref{Inteq}. We also leave as an open problem the calculation of higher cumulants of $T_{\infty,a}$ for dimensions $d>1$. Results of this kind would open ground to efficient approximation methods for the computation of critical values that avoid Monte Carlo simulations and efficiency statements, since the largest eigenvalue has a crucial influence on the approximate Bahadur efficiency, see \cite{B:1960,N:1995}. An promising new field of interest in connection with tests of multivariate normality is to consider their behaviour in high-dimensional settings, i.e., to answer the question whether one can find a suitable rescaling and shifting of the test statistic to obtain a non trivial limit distribution under a suitable limiting regime, under which, e.g., $n,d\rightarrow\infty$ such that $d/n\rightarrow\tau\in[0,\infty]$. For first results, see \cite{chen}. As a starting point, we conjecture that for a sequence $(n_d)_{d \in \N}$, where $n_d \geq d+1$ and $n_d=o\Big(\big(\frac{2a}{2a+1}\big)^{-\frac{d}{2}}\Big)$, we have under $H_0$ as $d \rightarrow \infty$
\begin{align*}
\Big(\frac{a}{\pi}\Big)^{\frac{d}{2}}\frac{T_{n_d,a}}{d}\fsk1.
\end{align*}
Finally, it would be of interest to consider a related family of test statistics, which is given by
\begin{equation*}
S_{n,a}= n\int_{\R^d} {\Vert \nabla \psi_n(t) + t \psi_n(t) \Vert}_{\C}^2 \ w_a(t) \, \text{d}t.
\end{equation*}
Thus, the theoretical CF in $T_{n,a}$ has been replaced by the empirical counterpart. Note that in the univariate case, this family is extensively studied in \cite{E:20}, but the generalisation to higher dimensions is still open. We conjecture that similar results as derived in Sections 2 to 4 hold for $S_{n,a}$.

\section*{Acknowledgment}
The authors thank Yvik Swan for sharing his knowledge of multivariate Stein operators and Stein characterisations.

\begin{appendix}
\section{Proofs}\label{secproofs}

\subsection{Proof of Theorem \ref{themresprstat}}

\begin{prf}
Putting $t= (t_1,\ldots,t_d)^\intercal \in \R^d$ and $Y_{n,j} = (Y_{n,j}^{(1)}, \ldots, Y_{n,j}^{(d)})^\intercal$, some algebra (using symmetry and the addition theorem for the cosine function)
yields
\begin{align*}
T_{n,a}&=n\int {\big\Vert \nabla \psi_n(t) + t \psi(t) \big\Vert}_{\C}^2 \ w_a(t) \ \text{d}t\\
&=n \int \Big\Vert \frac{1}{n} \sum_{j=1}^n \ii Y_{n,j} \exp\big(\ii t^{\intercal} Y_{n,j}\big) + t \psi(t) \Big\Vert^2_{\C} w_a(t) \ \text{d}t \\
&=n \int \Big\Vert \frac{1}{n} \sum_{j=1}^n \Big{\{} t \psi(t) - Y_{n,j} \sin(t^{\intercal}Y_{n,j}) + \ii Y_{n,j} \cos(t^{\intercal}Y_{n,j}) \Big{\}}\Big\Vert^2_{\C} w_a(t) \ \text{d}t\\
&=n \int \sum_{k=1}^d \Bigg{\{} \bigg( \frac{1}{n} \sum_{j=1}^n t^{(k)} \psi(t)  - Y_{n,j}^{(k)} \sin(t^{\intercal}Y_{n,j})\bigg)^2
 +  \bigg(\frac{1}{n} \sum_{j=1}^n Y_{n,j}^{(k)} \cos(t^{\intercal}Y_{n,j})\bigg)^2 \Bigg{\}} w_a(t) \ \text{d}t\\
&=n \int \sum_{k=1}^d \Big{\{} t^{(k)}t^{(k)} \exp\big(-{\Vert t \Vert}^2) - \frac{2}{n} \sum_{j=1}^n t^{(k)} \psi(t) Y_{n,j}^{(k)} \sin(t^{\intercal} Y_{n,j})\\
& \ \ \ \ \ \ \ \ \ \ \ \ \ \ \ + \frac{1}{n^2} \sum_{i,j=1}^n Y_{n,j}^{(k)} Y_{n,i}^{(k)} \cos(t^\intercal(Y_{n,i}-Y_{n,j})) \Big{\}} w_a(t) \ \text{d}t.
\end{align*}
We thus have
\begin{align*}
T_{n,a}&=n \int \bigg{\{} {\Vert t \Vert}^2 \exp\big(-(a+1) {\Vert t \Vert}^2\big)  - \frac{2}{n} \sum_{j=1}^n t^{\intercal} Y_{n,j} \sin(t^{\intercal} Y_{n,j}) \exp\Big(-\Big(a+\frac{1}{2}\Big){\Vert t \Vert}^2\Big) \\
& \ \ \ \ \ \ \ \ \ \ \ \ \ \ \ + \frac{1}{n^2} \sum_{i,j=1}^n Y_{n,i}^{\intercal} Y_{n,j} \cos\big(t^{\intercal}(Y_{n,i}-Y_{n,j})\big)\exp\big(-a {\Vert t \Vert}^2\big) \bigg{\}} \ \text{d}t.
\end{align*}
Using
\begin{align} \label{integrale1}
\int  {\Vert t \Vert}^2 \exp\big(-a {\Vert t \Vert}^2\big) \ \text{d}t &= {\bigg(\frac{\pi}{a}\bigg)}^\frac{d}{2} \frac{d}{2a}, \\
\int \cos(t^{\intercal} c) \exp\big(-a {\Vert t \Vert}^2\big) \ \text{d}t &= {\bigg(\frac{\pi}{a}\bigg)}^\frac{d}{2} \exp\bigg(-\frac{{\Vert c \Vert}^2}{4a}\bigg)  \label{integrale2}, \\
\int t^{\intercal}c \sin(t^{\intercal}c) \exp\big(-a{\Vert t \Vert}^2\big) \ \text{d}t &= {\bigg(\frac{\pi}{a}\bigg)}^\frac{d}{2} \frac{{\Vert c \Vert}^2}{2a} \exp\bigg(-\frac{{\Vert c \Vert}^2}{4a}\bigg), \label{integrale3}
\end{align}
the assertion follows readily.
\end{prf}

\subsection{Proof of Theorem \ref{thmlimitnull}}

\begin{prf}
Recall that, in view of invariance, there is no loss of generality if we assume $X \edist {\text N}_d(0,\ID)$.
With the notation in \eqref{CSp}, $Z_n$ defined in \eqref{Zn} takes the form
\[
Z_n(t) = \dfrac{1}{\sqrt{n}} \sum_{j=1}^n \big(Y_{n,j}\CSP(t,Y_{n,j})-t \psi(t)\big).
\]
To prove Theorem \ref{thmlimitnull}, we use  a central limit theorem for Hilbert space valued random elements, see, e.g., Theorem 2.7 of \cite{B:00}.
Since $Z_n$ does not comprise independent summands, we approximate $Z_n$ by a sum of i.i.d. random elements
of $\HH$. To this end, we introduce the auxiliary random elements
\begin{align}
\widetilde{Z}_n(t) &:= \dfrac{1}{\sqrt{n}} \sum_{j=1}^n \big( (X_j+ \Delta_{n,j})\CSP(t,X_j)-t\psi(t) + X_j  \CSM(t,X_j) t^{\intercal} \Delta_{n,j}\big), \nonumber \\
Z_n^*(t)&:= \dfrac{1}{\sqrt{n}} \sum_{j=1}^n \Big(X_j \CSP(t,X_j)-\big(t+X_j + (2\ID -tt^{\intercal}) \dfrac{1}{2}(X_jX_j^{\intercal}-\ID )t - t^{\intercal}X_jt\big)\psi(t)\Big)  \label{Zn*}\\
&=: \dfrac{1}{\sqrt{n}} \sum_{j=1}^n Z_j^{**}(t) \nonumber
\end{align}
(say), where
\begin{align} \label{delta}
 \Delta_{n,j}= Y_{n,j} - X_j = (S_n^{-\frac{1}{2}}-\ID)X_j - S_n^{-\frac{1}{2}} \overline{X}_n.
\end{align}
The proof of Theorem \ref{thmlimitnull} comprises 3 steps. We show
\begin{align}& \ Z_n^* \overset{\mathcal{D}}{\longrightarrow} Z \text{ in } \HH, \label{Schritt1}\\
& \ {\Vert Z_n - \widetilde{Z}_n \Vert}_{\HH} \overset{\mathbb{P}}{\longrightarrow} 0, \label{Schritt2}\\
& \ {\Vert \widetilde{Z}_n - Z_n^* \Vert}_{\HH} \overset{\mathbb{P}}{\longrightarrow} 0. \label{Schritt3}\end{align}
The assertion then follows from Slutsky's lemma. To prove \eqref{Schritt1}, notice that $Z_1^{**}, Z_2^{**}, \ldots $ is a
sequence of i.i.d. random elements of $\HH$. These elements are centred, since
\begin{align*}
\E [Z_1^{**} (t)] &= \E \Big[X \CSP(t,X)-\big(t+X + (2\ID -tt^{\intercal}) \dfrac{1}{2}(XX^{\intercal}-\ID )t - t^{\intercal}Xt\big)\psi(t)\Big]\\
&= \E \big[X \CSP(t,X)-t\psi(t)\big] = 0, \quad t\in \R^d.
\end{align*}
The covariance matrix kernel $\E\big{[}Z_n^*(s)Z_n^*(t)^{\intercal}\big{]} = \E\big{[}Z_1^{**}(s)Z_1^{**}(t)^{\intercal}\big{]} = K(s,t)$ (say), where $s,t \in \R^d$,
is given by
\begin{align*}
K(s,t)
&=\E \Big[ \Big(X \CSP(s,X)- \big(s + X + (2\ID -ss^{\intercal}) \dfrac{1}{2}(XX^{\intercal}-\ID )s - s^{\intercal}Xs\big) \psi(s)\Big) \\ & \ \ \ \ \ \ \ \ \Big(X \CSP(t,X)- \big(t + X + (2\ID -tt^{\intercal}) \dfrac{1}{2}(XX^{\intercal}-\ID )t - t^{\intercal}Xt\big) \psi(t)\Big)^{\intercal}\Big].
\end{align*}
In view of $\E[X]=0$ and $\E[XX^{\intercal}]=\ID$, tedious but straightforward calculations yield
\begin{align*}
K(s,t) &=\E\big[XX^{\intercal} \CSP(s,X)\CSP(t,X)\big] - s\psi(s)\E\big[X^{\intercal} \CSP(t,X)\big] - \psi(s) \E\big[XX^{\intercal} \CSP(t,X)\big] \\
& \ \ \ \ -\psi(s)\E \big[\big((2\ID -ss^{\intercal}) \dfrac{1}{2}(XX^{\intercal}-\ID ) - s^{\intercal}X\big)s X^{\intercal} \CSP(t,X) \big] \\
& \ \ \ \ - \E\big[ X \CSP(s,X)\big]t^{\intercal}\psi(t) +st^{\intercal} \psi(s) \psi(t) - \E\big[ X X^{\intercal} \CSP(s,X) \big] \psi(t) + \ID \psi(s) \psi(t) \\
& \ \ \ \ + \E \big[\big((2\ID -ss^{\intercal}) \dfrac{1}{2}(XX^{\intercal}-\ID ) - s^{\intercal}X\big) s X^{\intercal}\big] \psi(s) \psi(t) \\
& \ \ \ \ - \E \big[X \CSP(s,X)t^{\intercal}\big((2\ID -tt^{\intercal}) \dfrac{1}{2}(XX^{\intercal}-\ID ) - t^{\intercal}X\big)^{\intercal}\big] \psi(t)\\
& \ \ \ \ + \E \big[ X t^{\intercal}\big((2\ID -tt^{\intercal}) \dfrac{1}{2} (XX^{\intercal}-\ID ) - t^{\intercal}X\big)^{\intercal}\big] \psi(s) \psi(t) \\
& \ \ \ \ + \E \big[ \big((2\ID -ss^{\intercal}) \dfrac{1}{2}(XX^{\intercal}-\ID ) - s^{\intercal}X\big)s t^{\intercal}\big((2\ID -tt^{\intercal}) \dfrac{1}{2}(XX^{\intercal}-\ID ) - t^{\intercal}X\big)^{\intercal}\big] \psi(s) \psi(t).
\end{align*}
Since the occurring expectations are given by
\begin{align*}
\E \big[\CSP(t,X)\big]&= \psi(t), \\
\E \big[X \CSP(t,X)\big]&=t \psi(t), \\
\E \big[X \CSM(t,X)\big]&=-t\psi(t),  \\
\E \big[X X^{\intercal} \CSP(t,X)\big]&=(\ID - tt^{\intercal})\psi(t), \\
\E \big[X X^{\intercal} \CSM(t,X)\big]&=(\ID - tt^{\intercal})\psi(t), \\
\E \big[s^{\intercal}XXX^{\intercal}\CSP(t,X)\big]&= \big(s^{\intercal}t(\ID - tt^{\intercal})+st^{\intercal}+ts^{\intercal}\big)\psi(t), \\
\E \Big[XX^{\intercal} \CSP(s,X)\CSP(t,X)\Big]&=\E \Big[X X^{\intercal} \big(\sin(s+t) + \cos(s-t)\big) \Big]\\
&= \big(\ID - (s-t)(s-t)^{\intercal}\big) \psi(s-t),\\
\E \big[s^{\intercal}X X X^{\intercal}\big]&= 0\in\R^{d\times d},\\
\E \big[(X X^{\intercal}-\ID)st^{\intercal} (X X^{\intercal}-\ID)\big]&=ts^{\intercal} + s^{\intercal}t \ID,
\end{align*}
some algebra shows that $K(s,t)$  takes the form given in \eqref{kovark}. Thus, by the central limit theorem in Hilbert spaces, \eqref{Schritt1} follows. To prove \eqref{Schritt2}, notice that
\begin{align*}
\cos(t^{\intercal}Y_{n,j}) &= \cos(t^{\intercal}X_j) - \sin(t^{\intercal}X_j)t^{\intercal}\Delta_{n,j} + \varepsilon_{n,j}(t),\\
\sin(t^{\intercal}Y_{n,j}) &= \sin(t^{\intercal}X_j) + \cos(t^{\intercal}X_j)t^{\intercal}\Delta_{n,j} + \eta_{n,j}(t),
\end{align*}
where
\begin{align}
\max(|\varepsilon_{n,j}(t)|,|\eta_{n,j}(t)|) \le  \Vert t \Vert^2 \Vert \Delta_{n,j}\Vert^2. \label{norme}
\end{align}
Hence
\begin{align*}
\CSP(t,Y_{n,j})=\CSP(t,X_j)+\CSM(t,X_j)t^{\intercal}\Delta_{n,j} + \varepsilon_{n,j}(t) + \eta_{n,j}(t),
\end{align*}
and some algebra gives
\[
Z_n(t)-\widetilde{Z}_n(t)
 =\dfrac{1}{\sqrt{n}} \sum_{j=1}^n \big((X_j+\Delta_{n,j}) (\varepsilon_{n,j}(t) + \eta_{n,j}(t)) + \Delta_{n,j}\CSM(t,X_j)t^{\intercal} \Delta_{n,j}\big).
\]
Putting
\[
A_n = \dfrac{1}{\sqrt{n}}\sum_{j=1}^n 2 \Vert X_j \Vert {\Vert \Delta_{n,j}\Vert}^2,\quad
B_n = \dfrac{1}{\sqrt{n}}\sum_{j=1}^n 2 {\Vert \Delta_{n,j} \Vert}^2, \quad
C_n = \dfrac{1}{\sqrt{n}}\sum_{j=1}^n 2 {\Vert \Delta_{n,j} \Vert}^3,
\]
\eqref{norme} and the  Cauchy--Schwarz inequality yield
\[
\Vert Z_n(t) - \widetilde{Z}_n(t)\Vert \le  A_n {\Vert t \Vert}^2 + B_n \Vert t \Vert + C_n {\Vert t \Vert}^2.
\]
By Theorem 5.2 of \cite{B:63}, we have  $n^{-1/4} \max_{j=1,...,n} \Vert X_j \Vert \fsk 0$. Invoking  Proposition A.1 of
\cite{DEH:2019}, according to which  $n^{1/4} \max_{j=1,...,n} \Vert \Delta_{n,j} \Vert \fsk 0$ and $\sum_{j=1}^n \Vert \Delta_{n,j} \Vert^2 = O_{\mathbb{P}}(1)$,
it is readily seen that each of the expressions $A_n$, $B_n$ and $C_n$ converges to zero in probability as $n \to \infty$.
In view of
\[
{\Vert Z_n - \widetilde{Z}_n \Vert}_{\HH}^2 \leq \int \big(A_n{\Vert t \Vert}^2+ B_n\Vert t \Vert + C_n {\Vert t\Vert}^2\big)^2 w_a(t)\, {\text d}t
\]
the proof of \eqref{Schritt2} is finished. To prove \eqref{Schritt3}, we put
\begin{align*}
A_n(t)& = \dfrac{1}{\sqrt{n}} \sum_{j=1}^n \Big(\Delta_{n,j}\CSP(t,X_j) + \big(X_j + \dfrac{1}{2}(X_j X_j^{\intercal} - \ID ) t \big)\psi(t)\Big), \\
B_n(t)& = \dfrac{1}{\sqrt{n}} \sum_{j=1}^n \Big( X_j \CSM(t,X_j) t^{\intercal} \Delta_{n,j} + \big((\ID - tt^{\intercal}) \dfrac{1}{2}(X_jX_j^{\intercal}-\ID )t - t^{\intercal}X_jt\big)\psi(t)\Big).
\end{align*}
Using the triangle inequality, some calculations give ${\Vert \widetilde{Z}_n - Z_n^* \Vert}_{\HH}  \leq {\Vert A_n \Vert}_{\HH}\ + {\Vert B_n \Vert}_{\HH}$, and thus \eqref{Schritt3} follows if we can show that ${\Vert A_n \Vert}_{\HH} = o_\PP(1)$ and ${\Vert B_n \Vert}_{\HH} =o_\PP(1)$.
We only prove ${\Vert A_n \Vert}_{\HH} = o_\PP(1)$, since the reasoning for ${\Vert B_n \Vert}_{\HH} =o_\PP(1)$ is completely similar.
From the definition of $\Delta_{n,j}$ in \eqref{delta}, we have
\begin{align*}
A_n(t)&= (S_n^{-\frac{1}{2}}-\ID) \dfrac{1}{\sqrt{n}} \sum_{j=1}^n \big(X_j \CSP(t,X_j) - t \psi(t)\big) - S_n^{-\frac{1}{2}} \overline{X}_n \dfrac{1}{\sqrt{n}} \sum_{j=1}^n \big( \CSP(t,X_j) - \psi(t) \big)\\
& \ \ \ \ - \psi(t)\big(S_n^{-\frac{1}{2}} - \ID\big) \sqrt{n} \overline{X}_n  + \Big(\sqrt{n}(S_n^{-\frac{1}{2}}-\ID) + \dfrac{1}{2 \sqrt{n}} \sum_{j=1}^n \big(X_j X_j^{\intercal} - \ID\big)\Big) t \psi(t) \\
&= A_{n,1}(t) - A_{n,2}(t) - A_{n,3}(t) + A_{n,4}(t),
\end{align*}
say, and thus it remains to prove that each of ${\Vert A_{n,k} \Vert}_{\HH}$, $k \in \{1,2,3,4\}$, is $o_\PP(1)$.
Letting  ${\Vert \cdot \Vert}_2$  denote the spectral norm, it follows that
\[
{\Vert A_{n,1} \Vert}_{\HH}^2 \leq \big\Vert \sqrt{n} (S_n^{-\frac{1}{2}}-\ID) \big\Vert^2_2 \ \Big\Vert \dfrac{1}{n} \sum_{j=1}^n \big(X_j \CSP(t,X_j) - t \psi(t)\big) \Big\Vert^2_{\HH}.
\]
Here, the first factor on the right  hand side is $O_\PP(1)$, and the second converges to zero almost surely because of the strong law of
large numbers in $\HH$. As for ${\Vert A_{n,2} \Vert}_{\HH}^2$, it holds that
\[
{\Vert A_{n,2} \Vert}_{\HH}^2 \leq \big\Vert S_n^{-\frac{1}{2}} \big\Vert^2_2 {\big\Vert \sqrt{n}\overline{X}_n \big\Vert}^2 \Big\Vert \dfrac{1}{n} \sum_{j=1}^n \big( \CSP(t,X_j) - \psi(t) \big)  \Big\Vert_{L^2}^2.
\]
Here, each of the first two factors on the right hand side are $O_\PP(1)$, and the last one converges to zero almost surely because of the strong law of
large numbers in $L^2$. The term ${\Vert A_{n,3} \Vert}_{\HH}^2$ is bounded from above by
\[
{\Vert A_{n,3} \Vert}_{\HH}^2 \leq \Vert \sqrt{n}(S_n^{-\frac{1}{2}}-\ID) \Vert^2_2 \Vert \overline{X}_n \Vert^2  \int \exp(-{\Vert t \Vert}^2) w_a(t) \, {\text d}t.
\]
Hence ${\Vert A_{n,3} \Vert}_{\HH}^2 = o_\PP(1)$ since  $\Vert \overline{X}_n \Vert^2 = o_\PP(1)$. Finally, we have
\[
{\Vert A_{n,4} \Vert}_{\HH}^2 \leq \Big{\|} \sqrt{n}(S_n^{-\frac{1}{2}}-\ID)+\dfrac{1}{2\sqrt{n}}\sum_{j=1}^n(X_j X_j^{\intercal} - \ID) \Big{\|}^2_2 \int {\Vert t \Vert}^2 \exp(-{\Vert t \Vert}) w_a(t) \, {\text d}t.
\]
From display (2.13) of \cite{HW:1997}, the factor preceding the integral is $o_\PP(1)$, and thus  ${\Vert A_{n,4} \Vert}_{\HH}^2 = o_\PP(1)$.
The proof of Theorem \ref{thmlimitnull} is completed.
\end{prf}

\subsection{Proof of Theorem \ref{thmconstest}}

\begin{prf}
Since the proof is analogous to that given in \cite{DEH:2019}, it will only be sketched. The first observation is that
the quantities $\Psi_{\ell,n}(t)$, $\ell \in \{1,2,3,4\}$, defined in \eqref{defpsi12}, \eqref{defpsi34} have the following almost sure limits:
\[
\Psi_{1,n}(t) \fsk \psi^+_X(t), \ \Psi_{2,n}(t) \fsk  \nabla\psi^+_X(t), \ \Psi_{3,n}(t) \fsk  -\nabla\psi^-_X(t),\  \Psi_{4,n}(t) \fsk \E[XX^{\intercal}\CSM(t,X)].
\]
Here, the convergence of $\Psi_{3,n}(t)$ is assertion a) of Lemma 6.6 of \cite{DEH:2019}, and the remaining claims
follow mutatis mutandis the reasoning given in the proof of Lemma 6.6. of \cite{DEH:2019}.
From \eqref{Ln} and \eqref{defwnjt}, we have
\begin{align}\label{zerlln}
L_n(s,t)=\sum_{i,j=1}^5 L_n^{i,j}(s,t),
\end{align}
where $L_n^{i,j}(s,t)=L_n^{j,i}(t,s)^{\intercal}$ and -- putting ${\rm I}^\pm_{n,j} := Y_{n,j}Y_{n,j}^{\intercal} \pm \ID$ --
\begin{align*}
L_n^{1,1}(s,t)&= \dfrac{1}{n}\sum_{j=1}^n Y_{n,j}\CSP(s,Y_{n,j})Y_{n,j}^{\intercal}\CSP(t,Y_{n,j}), \quad
L_n^{1,2}(s,t)= -\dfrac{1}{n}\sum_{j=1}^n Y_{n,j}\CSP(s,Y_{n,j})Y_{n,j}^{\intercal}\Psi_{1,n}(t), \\
L_n^{1,3}(s,t) &= -\dfrac{1}{n}\sum_{j=1}^n  Y_{n,j}\CSP(s, Y_{n,j})t^{\intercal} Y_{n,j}\Psi_{2,n}(t)^{\intercal}, \quad
L_n^{1,4}(s,t) = -\dfrac{1}{2n}\sum_{j=1}^n  Y_{n,j} \CSP(s, Y_{n,j}) \Psi_{3,n}(t)^{\intercal}{\rm I}_{n,j}^+ , \\
L_n^{1,5}(s,t)& = -\dfrac{1}{2n}\sum_{j=1}^n Y_{n,j} \CSP(s,Y_{n,j})t^{\intercal}{\rm I}_{n,j}^- \Psi_{4,n}(t), \quad
L_n^{2,2} (s,t) = \dfrac{1}{n} \sum_{j=1}^n Y_{n,j} \Psi_{1,n}(s)  Y_{n,j}^{\intercal} \Psi_{1,n}(t), \\
L_n^{2,3} (s,t)& = \dfrac{1}{n} \sum_{j=1}^n Y_{n,j} \Psi_{1,n}(s) t^{\intercal}Y_{n,j}\Psi_{2,n}(t)^{\intercal}, \quad
L_n^{2,4} (s,t) = \dfrac{1}{2n} \sum_{j=1}^n Y_{n,j} \Psi_{1,n}(s) \Psi_{3,n}(t)^{\intercal}{\rm I}_{n,j}^+, \\
L_n^{2,5} (s,t)&= \dfrac{1}{2n} \sum_{j=1}^n Y_{n,j} \Psi_{1,n}(s) t^{\intercal}{\rm I}_{n,j}^- \Psi_{4,n}(t), \quad
L_n^{3,3} (s,t) = \dfrac{1}{n} \sum_{j=1}^n s^{\intercal}Y_{n,j} \Psi_{2,n}(s)  t^{\intercal} Y_{n,j} \Psi_{2,n}(t)^{\intercal}, \\
L_n^{3,4} (s,t) &= \dfrac{1}{2n} \sum_{j=1}^n s^{\intercal}Y_{n,j} \Psi_{2,n}(s) \Psi_{3,n}(t)^{\intercal}{\rm I}_{n,j}^+, \quad
L_n^{3,5} (s,t) = \dfrac{1}{2n} \sum_{j=1}^n s^{\intercal}Y_{n,j} \Psi_{2,n}(s) t^{\intercal}{\rm I}_{n,j}^- \Psi_{4,n}(t), \\
L_n^{4,4} (s,t) &= \dfrac{1}{4n} \sum_{j=1}^n {\rm I}_{n,j}^+ \Psi_{3,n}(s) \Psi_{3,n}(t)^{\intercal} {\rm I}_{n,j}^+, \quad
L_n^{4,5} (s,t) = \dfrac{1}{4n} \sum_{j=1}^n {\rm I}_{n,j}^+ \Psi_{3,n}(s) t^{\intercal}{\rm I}_{n,j}^- \Psi_{4,n}(t), \\
L_n^{5,5} (s,t) &= \dfrac{1}{4n} \sum_{j=1}^n \Psi_{4,n}(s){\rm I}_{n,j}^-s t^{\intercal}{\rm I}_{n,j}^- \Psi_{4,n}(t).
\end{align*}
From \eqref{zerlln}, it follows that $\widehat{\sigma}^2_{n,a}=\sum_{i,j=1}^5 \widehat{\sigma}^{i,j}_{n,a}$, where
\begin{align}\label{hatsigij}
\widehat{\sigma}^{i,j}_{n,a} =4\iint z_n(s)^{\intercal}L_n^{i,j}(s,t)z_n(t)w_a(s)w_a(t) \, {\text d}s \, {\text d}t.
\end{align}
Notice that $\widehat{\sigma}^{i,j}_{n,a}=\widehat{\sigma}^{j,i}_{n,a}$.
In view of \eqref{LST} and \eqref{defwsx}, we have $L(s,t) = \sum_{i,j=1}^5 L^{i,j}(s,t)$,
where $L^{i,j}(s,t) = \E[w_i(s,X)w_j(t,X)^{\intercal}]$, and
\begin{align*}
w_1(t,X) &= X\CSP(t,X), \quad w_2(t,X) =  -X\psi_X^+(t), \quad w_3(t,X) = -t^{\intercal}X \nabla \psi_X^+(t),\\
w_4(t,X) &= \dfrac{1}{2}(XX^{\intercal}+\ID)\nabla \psi_X^-(t), \quad w_5(t,X) = - \dfrac{1}{2} \E[XX^{\intercal}\CSM(t,X)](XX^{\intercal}-\ID)t.
\end{align*}
Therefore, $\sigma^2_a=\sum_{i,j=1}^5 \sigma^{i,j}_a$, where
\begin{align*}
\sigma^{i,j}_a =4\iint z(s)^{\intercal}L^{i,j}(s,t)z(t)w_a(s)w_a(t) \, {\text d}s \, {\text d}t
\end{align*}
and, by symmetry, $L^{i,j}(s,t)=L^{j,i}(t,s)^{\intercal}$ and hence $\sigma_a^{i,j}=\sigma_a^{j,i}$.
 We thus have to prove $ \widehat{\sigma}^{i,j}_{n,a} \stk \sigma^{i,j}_a$ for each choice of $i,j \in \{1, \ldots, 5\}$.
To this end, we proceed in two steps. The first one is to replace $L_n^{i,j}(s,t)$ in \eqref{hatsigij} with $L_{n,0}^{i,j}(s,t)$.
Here, $L_{n,0}^{i,j}(s,t)$ originates from $L_n^{i,j}(s,t)$ by throughout replacing $Y_{n,j}$ with $X_j$, and this replacement also affects
the quantities $\Psi_{\ell,n}(t)$, $\ell \in \{1,\ldots,4\}$. Moreover, we replace $z_n(t)$ with $z_{n,0}(t) = n^{-1}\sum_{j=1}^n X_j\CSP(t,X_j)-t \psi(t)$.
Putting
\[
\widehat{\sigma}_{n,0,a}^{i,j} = 4\iint z_{n,0}(s)^{\intercal}L_{n,0}^{i,j}(s,t)z_{n,0}(t)w_a(s)w_a(t) \, {\text d}s \, {\text d}t,
\]
it follows from Fubini's theorem that $\widehat{\sigma}_{n,0,a}^{i,j} \stk \sigma^{i,j}_a$. The second, much more technical step is
to prove $\widehat{\sigma}_{n,a}^{i,j} -\widehat{\sigma}_{n,0,a}^{i,j} = o_{\mathbb{P}}(1)$. To this end, notice that
\begin{align}
z_n(s)^{\intercal}L_n^{i,j}(s,t)z_n(t)-z_{n,0}(s)^{\intercal}L_{n,0}^{i,j}(s,t)z_{n,0}(t)=&z_n(s)^{\intercal}\big(L_n^{i,j}(s,t)-L_{n,0}^{i,j}(s,t)\big)z_n(t) \nonumber \\
& + \big(z_n(s)-z_{n,0}(s)\big)^{\intercal}L_{n,0}^{i,j}(s,t)z_n(t) \nonumber \\ &+ z_{n,0}(s)^{\intercal}L_{n,0}^{i,j}(s,t)\big(z_n(t)-z_{n,0}(t)\big), \label{differenz}
\end{align}
where
\begin{align*}
\big| \big(z_n(s)-z_{n,0}(s)\big)^{\intercal}L_{n,0}^{i,j}(s,t)z_n(t) \big| & \leq \big\Vert z_n(s) - z_{n,0}(s)\big\Vert \big\Vert L_{n,0}^{i,j}(s,t) \big\Vert_2 \big\Vert z_n(t) \big\Vert,\\
\big| z_{n,0}(s)^{\intercal}L_{n,0}^{i,j}(s,t)\big(z_n(t)-z_{n,0}(t)\big) \big| & \leq \big\Vert z_{n,0}(s)\big\Vert \big\Vert L_{n,0}^{i,j}(s,t) \big\Vert_2 \big\Vert z_n(t)-z_{n,0}(t) \big\Vert.
\end{align*}
We have $\Vert z_{n,0}(t) \Vert \leq 2n^{-1}\sum_{j=1}^n \Vert X_j \Vert + \Vert t \Vert \psi (t)$, and a Taylor expansion yields
\begin{eqnarray*}
\Vert z_{n}(t) \Vert \! & \! \leq \! & \! \dfrac{2}{n}\sum_{j=1}^n \big(\Vert X_j \Vert + \Vert X_j \Vert \Vert t \Vert \Vert \Delta_{n,j} \Vert +\Vert \Delta_{n,j} \Vert + \Vert t \Vert \Vert \Delta_{n,j} \Vert^2\big) + \Vert t \Vert \psi (t), \\
\Vert z_n(t) - z_{n,0}(t) \Vert \! & \! \leq \! & \!   \dfrac{2}{n} \sum_{j=1}^n \Vert \Delta_{n,j} \Vert + \dfrac{2\Vert t \Vert}{n} \sum_{j=1}^n \Vert \Delta_{n,j} \Vert \Vert X_j \Vert.
\end{eqnarray*}
Notice that each of the terms $\Vert L_{n,0}^{i,j}(s,t)\Vert_2 $ are bounded from above by terms of the type
     $2^k \Vert s \Vert^\ell \Vert t \Vert^m$, multiplied with finitely many products of the type  $n^{-1} \sum_{j=1}^n \Vert X_j \Vert^{\beta}$,
with  $k\leq 2$, $\ell,m \in \{0,1\}$,  and $\beta \in \{1,2,3,4\}$.
In view of the condition $\E \Vert X \Vert^4<\infty$ and the fact that $n^{-1} \sum_{j=1}^n \Vert \Delta_{n,j}\Vert^k \Vert X_k\Vert^\ell \fsk 0$ (see Proposition A.2 of \cite{DEH:2019}),
it follows that
\begin{align*}
 \iint \big| \big(z_n(s)-z_{n,0}(s)\big)^{\intercal}L_{n,0}^{i,j}(s,t)z_n(t) \big| w_a(s)w_a(t)\, {\text d}s \, {\text d}t& \stk 0, \\
 \iint \big| z_{n,0}(s)^{\intercal}L_{n,0}^{i,j}(s,t)\big(z_n(t)-z_{n,0}(t)\big) \big| w_a(s) w_a(t) \, {\text d}s \, {\text d}t \stk 0.
\end{align*}
As a consequence, we only have  to consider the first term on the right hand side  of \eqref{differenz}.
To this end, notice that
\begin{align*}
\big|z_n(s)^{\intercal}\big(L_n^{i,j}(s,t)-L_{n,0}^{i,j}(s,t)\big)z_n(t) \big|\leq \big\Vert z_n(s) \big\Vert \big\Vert L_n^{i,j}(s,t)-L_{n,0}^{i,j}(s,t)\big \Vert_2 \big \Vert z_n(t) \big \Vert .
\end{align*}
To find an upper bound for $\Vert L^{i,j}_n(s,t) - L^{i,j}_{n,0}(s,t) \Vert_2$, we have to consider each case $i,j \in \{1,\ldots,5\}$ such that $i \le j$ separately.
We will elaborate on the case $i=j=1$; the other cases are treated similarly. We have
\begin{equation*}
\big\Vert L^{1,1}_n(s,t) - L^{1,1}_{n,0}(s,t) \big\Vert_2 = \Big\Vert  \dfrac{1}{n}\sum_{j=1}^n \big(Y_{n,j}Y_{n,j}^{\intercal}\CSP(s,Y_{n,j})\CSP(t,Y_{n,j})-  X_j X_j^{\intercal}\CSP(s,X_j)\CSP(t,X_j)\big) \Big\Vert_2,
\end{equation*}
and a Taylor expansion yields
\begin{align*}
\Vert L^{1,1}_n(s,t) - L^{1,1}_{n,0}(s,t) \Vert_2 \leq & \dfrac{4}{n} \sum_{j=1}^n \Vert X_j \Vert^2 \big(\Vert t  \Vert \Vert \Delta_{n,j} \Vert + \Vert s  \Vert \Vert \Delta_{n,j} \Vert\big)+\dfrac{4}{n} \sum_{j=1}^n \Vert X_j \Vert^2 \big(\Vert s \Vert \Vert t \Vert \Vert \Delta_{n,j} \Vert^2\big)  \\
&+\dfrac{8}{n} \sum_{j=1}^n \Vert X_j \Vert \Vert \Delta_{n,j} \Vert \big(1 + \Vert s \Vert \Vert \Delta_{n,j} \Vert\big)\big(1 + \Vert t \Vert \Vert \Delta_{n,j} \Vert\big) \\
&+\dfrac{4}{n} \sum_{j=1}^n \Vert \Delta_{n,j} \Vert^2 \big(1 + \Vert s \Vert \Vert \Delta_{n,j} \Vert\big)\big(1 + \Vert t \Vert \Vert \Delta_{n,j} \Vert\big).
\end{align*}
From Proposition A.2 of \cite{DEH:2019}, it follows that $\Vert L^{1,1}_n(s,t) - L^{1,1}_{n,0}(s,t) \Vert_2 \fsk 0$.

To prove \eqref{intfrei},
we need the integrals
\begin{align*}
L_{1,a} (x)&:= \int t\psi(t) \CSP(t,x) w_a(t) \, {\text d}t = \dfrac{(2\pi)^{\frac{d}{2}}}{(2a+1)^{\frac{d}{2}+1}} x \exp\Big(-\dfrac{\Vert x \Vert^2}{4a+2}\Big),\\
L_{2,a} (x)&:= \int t t^{\intercal}x \psi(t) \CSM(t,x) w_a(t) \, {\text d}t \\ &\ = \dfrac{(2\pi)^{\frac{d}{2}}}{(2a+1)^{\frac{d}{2}+2}} \big((2a+1)x-\Vert x \Vert^2 x\big) \exp\Big(-\dfrac{\Vert x \Vert^2}{4a+2}\Big),\\
I_{1,a} (x,y)&:= \int \CSP(t,x)\CSP(t,y) w_a(t) \, {\text d}t = \Big(\dfrac{\pi}{a}\Big)^{\frac{d}{2}}\exp\Big(-\dfrac{\Vert x-y \Vert^2}{4a}\Big),  \\
I_{2,a} (x,y)&:= \int t\CSP(t,x)\CSM(t,y) w_a(t) \, {\text d}t = \Big(\dfrac{\pi}{a}\Big)^{\frac{d}{2}}\dfrac{(x-y)}{2a}\exp\Big(-\dfrac{\Vert x-y \Vert^2}{4a}\Big). \\
\end{align*}
Putting
\begin{align*}
P_{1,a}^{i,j}:=&Y_{n,i}^{\intercal}Y_{n,j}I_{1,a}(Y_{n,i},Y_{n,j})-L_{1,a}(Y_{n,j})^{\intercal}Y_{n,j}, \\
P_{2,a}^{i,j,k}:=&Y_{n,i}^{\intercal}Y_{n,j}I_{1,a}(Y_{n,i},Y_{n,k})-L_{1,a}(Y_{n,k})^{\intercal}Y_{n,j}, \\
P_{3,a}^{i,j,k}:=&Y_{n,i}^{\intercal}Y_{n,k} Y_{n,j}^{\intercal}I_{2,a}(Y_{n,i},Y_{n,k})-Y_{n,j}^{\intercal}L_{2,a}(Y_{n,k}), \\
P_{4,a}^{i,j,k} :=&Y_{n,i}^{\intercal}(Y_{n,j}Y_{n,j}^{\intercal}+\ID)Y_{n,k}I_{1,a}(Y_{n,i},Y_{n,k}) -Y_{n,k}^{\intercal}(Y_{n,j}Y_{n,j}^{\intercal}+\ID)L_{1,a}(Y_{n,k}), \\
P_{5,a}^{i,j,k} :=&Y_{n,i}^{\intercal}Y_{n,k}Y_{n,k}^{\intercal}(Y_{n,j}Y_{n,j}^{\intercal}-\ID)I_{2,a}(Y_{n,i},Y_{n,k}) -Y_{n,k}^{\intercal}(Y_{n,j}Y_{n,j}^{\intercal}-\ID)L_{2,a}(Y_{n,k}),
\end{align*}
straightforward calculations give
\begin{align}\label{defsigij}
\widehat{\sigma}_{n,a}^{1,1}=&\dfrac{4}{n^3}\sum_{i,j,k=1}^n P_{1,a}^{i,j} P_{1,a}^{k,j}, \quad
\widehat{\sigma}_{n,a}^{1,2}= -\dfrac{4}{n^4} \sum_{i,j,k,\ell=1}^n P_{1,a}^{i,j}P_{2,a}^{\ell,j,k}, \\ \nonumber
\widehat{\sigma}_{n,a}^{1,3}=&-\dfrac{4}{n^4}\sum_{i,j,k,\ell=1}^n P_{1,a}^{i,j}P_{3,a}^{\ell,j,k}, \quad
\widehat{\sigma}_{n,a}^{1,4}=-\dfrac{2}{n^4} \sum_{i,j,k,\ell=1}^n P_{1,a}^{i,j} P_{4,a}^{\ell,j,k}, \\ \nonumber
\widehat{\sigma}_{n,a}^{1,5}=&-\dfrac{2}{n^4} \sum_{i,j,k,\ell=1}^n P_{1,a}^{i,j} P_{5,a}^{\ell,j,k}, \quad
\widehat{\sigma}_{n,a}^{2,2}= \dfrac{4}{n^5}\sum_{i,j,k,\ell,m=1}^n P_{2,a}^{i,j,k} P_{2,a}^{m,j,\ell}, \\ \nonumber
\widehat{\sigma}_{n,a}^{2,3}=&\dfrac{4}{n^5}\sum_{i,j,k,\ell,m=1}^n P_{2,a}^{i,j,k} P_{3,a}^{m,j,\ell},\quad
\widehat{\sigma}_{n,a}^{2,4}=\dfrac{2}{n^5}\sum_{i,j,k,l,m=1}^n P_{2,a}^{i,j,k} P_{4,a}^{m,j,\ell},\\ \nonumber
\widehat{\sigma}_{n,a}^{2,5}=&\dfrac{2}{n^5}\sum_{i,j,k,l,m=1}^n P_{2,a}^{i,j,k} P_{5,a}^{m,j,\ell},\quad
\widehat{\sigma}_{n,a}^{3,3}=\dfrac{4}{n^5}\sum_{i,j,k,\ell,m=1}^n P_{3,a}^{i,j,k} P_{3,a}^{m,j,\ell},\\ \nonumber
\widehat{\sigma}_{n,a}^{3,4}=&\dfrac{2}{n^5}\sum_{i,j,k,\ell,m=1}^n P_{3,a}^{i,j,k} P_{4,a}^{m,j,\ell},\quad
\widehat{\sigma}_{n,a}^{3,5}=\dfrac{2}{n^5}\sum_{i,j,k,\ell,m=1}^n P_{3,a}^{i,j,k} P_{5,a}^{m,j,\ell},\\ \nonumber
\widehat{\sigma}_{n,a}^{4,4}=&\dfrac{1}{n^5}\sum_{i,j,k,\ell,m=1}^n P_{4,a}^{i,j,k} P_{4,a}^{m,j,\ell},\quad
\widehat{\sigma}_{n,a}^{4,5}=\dfrac{1}{n^5}\sum_{i,j,k,\ell,m=1}^n P_{4,a}^{i,j,k} P_{5,a}^{m,j,\ell},\\ \nonumber
\widehat{\sigma}_{n,a}^{5,5}=&\dfrac{1}{n^5}\sum_{i,j,k,\ell,m=1}^n P_{5,a}^{i,j,k} P_{5,a}^{m,j,\ell}.
\end{align}
\end{prf}

\newpage

\begin{landscape}
\begin{align*}
\kappa_3(a)=&\frac {65536{\pi}^{3/2}}{ \left( 4a^{2}+8a+3 \right)^{2} \left( 2{a}^{2}+4a+1 \right)^{11/2}\sqrt {a} \left( 2a+3
 \right)^{7/2} \left( 2a+1 \right)^{3/2} \left( a+1 \right)^{21/2}} \Big(  \Big(  \Big( \frac {117321111 \sqrt {2\,a+3}\sqrt {2\,a
+1}}{65536} \Big( \frac {1362117632}{39107037}{a}^{\frac{43}{2}}+
\frac {32768}{39107037}{a}^{\frac{53}{2}}\\
&+\frac {9003008}{39107037}{a}^{\frac{49}{2}}+\frac {262144}{13035679}{a}^{\frac{51}{2}}+\frac{43947161792}{39107037}{a}^{\frac{35}{2}}+\frac {
8042308544}{13035679}{a}^{\frac{37}{2}}+\frac {68156023424}{
39107037}{a}^{\frac{33}{2}}+\frac {103075274992}{39107037}{a}^{
\frac{29}{2}}+\frac {100834016080}{39107037}{a}^{\frac{27}{2}}+
\frac {84594302576\,}{39107037}{a}^{\frac{25}{2}}\\
&+\frac {
60691508044\,}{39107037}{a}^{\frac{23}{2}}+\frac {37074189596\,}{
39107037}{a}^{\frac{21}{2}}+\frac {6389004348\,}{13035679}{a}^{\frac{19}{2}}+\frac {
2773185824\,}{13035679}{a}^{\frac{17}{2}}+\frac {3000075335\,}{
39107037}{a}^{\frac{15}{2}}+\frac {886275569\,}{39107037}{a}^{\frac{13}{2}}+\frac {70150049
\,}{13035679}{a}^{\frac{11}{2}}+{a}^{\frac{9}{2}}\\
&+\frac {1822035\,}{
13035679}{a}^{\frac{7}{2}}+\frac {179757\,}{13035679}{a}^{\frac{5}{2}}+\frac {85\,
}{99509}{a}^{\frac{3}{2}}+\frac {325\,}{13035679}\sqrt {a}+\frac {90474975952}{
39107037}{a}^{\frac{31}{2}}+\frac {341954560}{39107037}{a}^{\frac
{45}{2}}+\frac {65572864}{39107037}{a}^{\frac{47}{2}}+\frac {
4316055040}{39107037}{a}^{\frac{41}{2}}+\frac {11180065088}{
39107037}{a}^{\frac{39}{2}} \Big) \\
&+\frac { \left( 16\,{a}^{4}+
48\,{a}^{3}+72\,{a}^{2}+84\,a+45 \right)  \left( 2\,a+1 \right) ^{3}
 \left( 2\,{a}^{2}+4\,a+1 \right) ^{5} \left( a+1 \right) ^{21/2}}{
4096} \Big) \sqrt {4\,{a}^{2}+8\,a+3}-\frac {4557280077\,\sqrt {2
\,a+3}\sqrt {2\,a+1}}{4194304} \Big( \frac {237774045184}{1519093359
}{a}^{\frac{43}{2}}\\
&+\frac {104857600}{4557280077}{a}^{\frac{53}{2
}}+\frac {4194304}{4557280077}{a}^{\frac{55}{2}}+\frac {
9500622848}{4557280077}{a}^{\frac{49}{2}}+\frac {1250951168}{
4557280077}{a}^{\frac{51}{2}}+\frac {13772677078016}{4557280077}{a
}^{\frac{35}{2}}+\frac {8413599202304}{4557280077}{a}^{\frac{37}{
2}}+\frac {19360446462976}{4557280077}{a}^{\frac{33}{2}}\\
&+\frac
{24491194230784}{4557280077}{a}^{\frac{29}{2}}+\frac {
3150868184512}{651040011}{a}^{\frac{27}{2}}+\frac {17092449262400
\,}{4557280077}{a}^{\frac{25}{2}}+\frac {3787061754352\,
}{1519093359}{a}^{\frac{23}{2}}+\frac {6446977759604\,}{4557280077}{a}^{\frac{21}{2}}+\frac
{3103648635392\,}{4557280077}{a}^{\frac{19}{2}}\\
&+\frac {1257161559490\,}{4557280077}{a}^{
\frac{17}{2}}+\frac {141288311512\,}{1519093359}{a}^{\frac{15}{2}}+
\frac {13030129762\,}{506364453}{a}^{\frac{13}{2}}+\frac {322428332\,}{56262717}{a}^{\frac{11}{2}}+{a}^{\frac{9}{2}}+\frac {117347\,}{893059}{a}^{\frac{7}{2}}+\frac {
43583\,}{3572236}{a}^{\frac{5}{2}}+\frac {636\,}{893059}{a}^{\frac{3}{2}}\\
&+\frac {
981\,}{50011304}\sqrt {a}+\frac {2605134412544}{506364453}{a}^{
\frac{31}{2}}+\frac {71725285376}{1519093359}{a}^{\frac{45}{2}}+
\frac {7386431488}{651040011}{a}^{\frac{47}{2}}+\frac {
645474689024}{1519093359}{a}^{\frac{41}{2}}+\frac {1462751346688}{
1519093359}{a}^{\frac{39}{2}} \Big)  \Big) \sqrt {2\,{a}^{2}+4
\,a+1}\\
&-\frac {63481963221\,\sqrt {2\,a+3}\sqrt {2\,a+1}\sqrt {2}}{
16777216} \Big( \frac {397764591616}{783727941}{a}^{\frac{43}{2}}
+\frac {5431623680}{21160654407}{a}^{\frac{53}{2}}+\frac {
436207616}{21160654407}{a}^{\frac{55}{2}}+\frac {246551674880}{
21160654407}{a}^{\frac{49}{2}}+\frac {16777216}{21160654407}{a}^{
\frac{57}{2}}\\
&+\frac {43182456832}{21160654407}{a}^{\frac{51}{2}}
+\frac {133654614470656}{21160654407}{a}^{\frac{35}{2}}+\frac {
89697147781120}{21160654407}{a}^{\frac{37}{2}}+\frac {
172252205867008}{21160654407}{a}^{\frac{33}{2}}+\frac {
186440146500608}{21160654407}{a}^{\frac{29}{2}}+\frac {
17390571188992}{2351183823}{a}^{\frac{27}{2}}\\
&+\frac {
12621837363200\,}{2351183823}{a}^{\frac{25}{2}}+\frac {
876956342720\,}{261242647}{a}^{\frac{23}{2}}+\frac {38087012979520\,}{21160654407}{a}^{\frac{21}{
2}}+\frac {17401201904624\,}{21160654407}{a}^{\frac{19}{2}}+
\frac {6719460948704\,}{21160654407}{a}^{\frac{17}{2}}\\
&+\frac {723312731296\,
}{7053551469}{a}^{\frac{15}{2}}+\frac {192635967404\,}{7053551469}{a}^{\frac{13}{2}}
+\frac {1537776695\,}{261242647}{a}^{\frac{11}{2}}+{a}^{\frac{9}{2}}+\frac {
33641318\,}{261242647}{a}^{\frac{7}{2}}+\frac {3084786\,}{261242647
}{a}^{\frac{5}{2}}+\frac {179523\,}{261242647}{a}^{\frac{3}{2}}\\
&+\frac {4995\,}{
261242647}\sqrt {a}+\frac {10127956063232}{1113718653}{a}^{\frac{31}{2}}+
\frac {139213668352}{783727941}{a}^{\frac{45}{2}}+\frac {
1078067462144}{21160654407}{a}^{\frac{47}{2}}+\frac {2848598917120
}{2351183823}{a}^{\frac{41}{2}}+\frac {51845834817536}{21160654407
}{a}^{\frac{39}{2}} \Big)  \Big)
\end{align*}
\begin{align*}
\kappa_4(a)=&\frac {4238729565{\pi}^{2}}{\sqrt {a}
 ( 4\,{a}^{2}+8\,a+3 ) ^{11/2}\sqrt {{a}^{2}+3\,a+2}
 ( a+3/2 ) ^{4} ( {a}^{2}+2\,a+1/2 ) ^{6}
 ( {a}^{2}+5/2\,a+5/4 ) ^{7} ( a+2 ) ^{4}
 ( a+1 ) ^{14} ( {a}^{2}+3/2\,a+1/4 ) ^{5}}
 \\
&\bigg(-\frac {185906111765413888\,\sqrt {4\,{a}^{2}+8\,a+3}\sqrt {
16\,{a}^{4}+64\,{a}^{3}+84\,{a}^{2}+40\,a+5}\sqrt {{a}^{2}+3\,a+2}}{
189423662804146875} \bigg( \frac {19802228366444765118464}{
46529006931033345}{a}^{\frac{111}{2}} \\
&+\frac {
142124794512189343251826622464}{46529006931033345}{a}^{\frac{81}{2}}
+\frac {1258096727350689610082\,}{3101933795402223}{a}^{\frac{23}{2}}
+
\frac {5826979313322568108292964352}{46529006931033345}{a}^{\frac{91}
{2}}\\
&+\frac {15750308417458890175577472592}{15509668977011115}{a}^{
\frac{49}{2}}+\frac {138611297001425344382631876536}{
46529006931033345}{a}^{\frac{53}{2}}+{a}^{\frac{11}{2}}+
\frac {
27942998903860339337720467216}{15509668977011115}{a}^{\frac{51}{2}}\\
&+\frac {704490872277734804863483698176}{46529006931033345}{a}^{\frac
{67}{2}}+\frac {4463906604\,}{1418351072429}{a}^{\frac{7}{2}}+\frac {
5679222185915318272}{5169889659003705}{a}^{\frac{117}{2}}+\frac {
161688996\,}{1418351072429}{a}^{\frac{5}{2}}+\frac {
14538236644323150732396593152}{3101933795402223}{a}^{\frac{79}{2}}\\
&+\frac {6847822034336568046765835552}{1033977931800741}{a}^{\frac{57}
{2}}+\frac {542999500682333909942272}{46529006931033345}{a}^{\frac
{107}{2}}+\frac {266964236560762308395008}{5169889659003705}{a}^{
\frac{105}{2}}+\frac {36658643371530037231616}{15509668977011115}{a
}^{\frac{109}{2}}\\
&+\frac {86496238983532270324100694016}{46529006931033345}{a}^{\frac{83}{2}}+\frac {
68164570063096914777497\,}{31019337954022230}{a}^{\frac{25}{2}}+
\frac {3139144012127128256512}{46529006931033345}{a}^{\frac{113}{2}}
+\frac {43848\,}{1418351072429}\sqrt {a}\\
&+\frac {206569108354665664014319566848}{15509668977011115}{a}^{\frac{71}{2}}+\frac {12765734064419898295898341376}{46529006931033345}{a}^{\frac
{89}{2}}+\frac {2069793471986052845657\,}{
31019337954022230}{a}^{\frac{21}{2}}\\
&+\frac {75850264590928714946676895744}{
5169889659003705}{a}^{\frac{69}{2}}+\frac {
415069977156299500403300811584}{46529006931033345}{a}^{\frac{59}{2}}+\frac {25961321024761185580881805312}{46529006931033345}{a}^{\frac
{87}{2}}+\frac {326117721828947\,}{25530319303722}{a}^{\frac{13}{2}}\\
&+\frac {549755813888}{46529006931033345}{a}^{\frac{129}{2}}+\frac {
432833812682392993792}{46529006931033345}{a}^{\frac{115}{2}}+\frac {616767662358185891720952447616}{46529006931033345}{a}^{\frac{
63}{2}}+\frac {136218522380459322233433523840}{9305801386206669}{a}
^{\frac{65}{2}}\\
&+\frac {34738436701499146652090368}{
46529006931033345}{a}^{\frac{101}{2}}+\frac {
174328760212840680730477047008}{15509668977011115}{a}^{\frac{61}{2}}+\frac {1688466955339563008}{15509668977011115}{a}^{\frac{119}{2}}+\frac {497158290836194024223147}{46529006931033345}{a}^{\frac{27}{
2}}\\
&+\frac {312764398606697372000140730368}{46529006931033345}{a}^{
\frac{77}{2}}+\frac {1209875107414016}{46529006931033345}{a}^{
\frac{125}{2}}+\frac {71283965626400484613576759456}{
15509668977011115}{a}^{\frac{55}{2}}+\frac {49099913781052240700765831168}{46529006931033345}{a}^{\frac{85}{2}}\\
&+\frac {38233385349739640837373952}{15509668977011115}{a}^{\frac{99}
{2}}+\frac {10035973055783223923\,}{1033977931800741}{a}^{\frac{19}{2}}+\frac {3197145137334020848222208}{15509668977011115}{a}^{\frac{103}{2
}}+\frac {3809952\,}{1418351072429}{a}^{\frac{3}{2}}+\frac {25969914892255232}{46529006931033345}{a}^{\frac{123}{2}}\\
&+\frac {419618082332069948495257899008}{46529006931033345}{a}^{\frac{75}{2}}+\frac {180243105937\,}{2836702144858}{a}^{\frac{9}{2}}+\frac {
526962212776409628068072324096}{46529006931033345}{a}^{\frac{73}{2}}
+\frac {36833639530496}{46529006931033345}{a}^{\frac{127}{2}}\\
&+\frac {94677809938011743\,}{76590957911166}{a}^{\frac{17}{2}}+\frac {5214542156552656\,}{38295478955583}{a}^{\frac{15}{2}}+\frac {
2463825296142995483854372864}{46529006931033345}{a}^{\frac{93}{2}}+
\frac {20997863111128412747625073}{31019337954022230}{a}^{\frac{33}{
2}}\\
&+\frac {7704249583566767268036608}{1033977931800741}{a}^{\frac{
97}{2}}+\frac {962732137091860887777574912}{46529006931033345}{a}^{
\frac{95}{2}}+\frac {24790840891247406265319616152}{
46529006931033345}{a}^{\frac{47}{2}}+\frac {8685759319100928242854972}{46529006931033345}{a}^{\frac{31}{2}}
\end{align*}
\begin{align*}
&+\frac {12089577733297066123987509806}{46529006931033345}{a}^{\frac{45}{2}}+\frac {81939042475704320}{9305801386206669}{a}^{\frac{121}{2
}}+\frac {1454598168894702388874699}{31019337954022230}{a}^{\frac{
29}{2}}+\frac {20885976543831443714128621}{9305801386206669}{a}^{
\frac{35}{2}}\\
&+\frac {635230038229728431368925273}{93058013862066690
}{a}^{\frac{37}{2}}+\frac {888683347357380981709347746}{46529006931033345}{a}^{\frac{39}{2}}+\frac {
2293364834214572783320577806}{46529006931033345}{a}^{\frac{41}{2}}+
\frac {5471008793804977527437466572}{46529006931033345}{a}^{\frac{43}{2}} \Big) \\
&+ \Big( \frac {36028797018963968\, ( 1/2+a
 ) ^{5} ( a+3/2 ) ^{9} ( {a}^{2}+2\,a+1/2
 ) ^{5} ( {a}^{2}+5/2\,a+5/4 ) ^{7} ( a+1
 ) ^{9} ( {a}^{2}+3/2\,a+1/4 ) ^{5}\sqrt {2\,{a}^{2}+
4\,a+1}}{6214043258290038234375} \\
&\Big( {a}^{8}+8\,{a}^{7}+28\,{a}^{6}
+56\,{a}^{5}+76\,{a}^{4}+80\,{a}^{3}
+63\,{a}^{2}+30\,a+\frac{105}{16}
 \Big) +\sqrt {2}\sqrt {{a}^{2}+3\,a+2}\sqrt {a+1} \Big( \frac {1288286823500979073414332416}{2071347752763346078125}{a}^{\frac{111}{
2}}\\
&+\frac {8935032163150251256650696230961152}{
2071347752763346078125}{a}^{\frac{81}{2}}+\frac {
839551898119154107328624\,}{1841198002456307625}{a}^{\frac{23}{2}} +\frac {371229768288716262827380170555392}{2071347752763346078125}{a}^{\frac{91}{2}}\\
&+\frac {913849879634598027664546373083136}{
690449250921115359375}{a}^{\frac{49}{2}}+\frac {
8164067143934275953913884357165056}{2071347752763346078125}{a}^{\frac
{53}{2}}+{a}^{\frac{11}{2}}
+\frac {1633851323087451063964056414027776}{
690449250921115359375}{a}^{\frac{51}{2}}\\
&+\frac {
43189862757173612265706971713241088}{2071347752763346078125}{a}^{
\frac{67}{2}}+\frac {60292644585\,}{20205190699109}{a}^{\frac{7}{2}}+
\frac {3338995139626653383131136}{2071347752763346078125}{a}^{{\frac{
117}{2}}}+\frac {2117413575\,}{20205190699109}{a}^{\frac{5}{2}}
\\
&+\frac {911256225356037430357120383975424}{138089850184223071875}{a}^{\frac{
79}{2}}+\frac {6131265682608049365461479693877248}{
690449250921115359375}{a}^{\frac{57}{2}}+\frac {
11734861729999873853970448384}{690449250921115359375}{a}^{\frac{107}{
2}}\\
&+\frac {51826136668193132779021533184}{690449250921115359375}{a}
^{\frac{105}{2}}+\frac {476192629116286924104925184}{
138089850184223071875}{a}^{\frac{109}{2}}+\frac {
5453341064484344228930998088237056}{2071347752763346078125}{a}^{\frac
{83}{2}}\\
&+\frac {13857882938167540266083194\,}{
5523594007368922875}{a}^{\frac{25}{2}}+\frac {93524065028843537170432}{
947118314020734375}{a}^{\frac{113}{2}}+\frac {533250\,}{
20205190699109}\sqrt {a}+\frac {38311298643647789519148389660360704}{
2071347752763346078125}{a}^{\frac{71}{2}}\\
&+\frac {
270421671676414119465208552357888}{690449250921115359375}{a}^{{\frac{
89}{2}}}+\frac {163057876019485508405021\,}{
2209437602947569150}{a}^{\frac{21}{2}}+\frac {8406507958530434181259944060256256}{
414269550552669215625}{a}^{\frac{69}{2}}\\
&+\frac {
923215837049551012695512078876672}{76716583435679484375}{a}^{\frac{59
}{2}}+\frac {548541331404397319980001194409984}{
690449250921115359375}{a}^{\frac{87}{2}}+\frac {9510361143154141\,
}{727386865167924}{a}^{\frac{13}{2}}+\frac {36028797018963968}{
2071347752763346078125}{a}^{\frac{129}{2}}\\
&+\frac {
28240943827444880352739328}{2071347752763346078125}{a}^{\frac{115}{2}
}+\frac {2880777548978561308424901499027456}{159334442520257390625}
{a}^{\frac{63}{2}}+\frac {4618011507571792344706555767160832}{
230149750307038453125}{a}^{\frac{65}{2}}\\
&+\frac {
2238763503046729775813824086016}{2071347752763346078125}{a}^{\frac{
101}{2}}+\frac {10529013623794787999533082298810368}{
690449250921115359375}{a}^{\frac{61}{2}}+\frac {
331244494883507601932288}{2071347752763346078125}{a}^{\frac{119}{2}}\\
&+\frac {1025507789564528843568643688}{82853910110533843125}{a}^{
\frac{27}{2}}+\frac {501089627561891102133336861048832}{53111480840085796875}{a}^{\frac{77}{2}}+\frac {6097873895459651584
}{159334442520257390625}{a}^{\frac{125}{2}}\\
&+\frac {
12682606898373656937196513212366848}{2071347752763346078125}{a}^{
\frac{55}{2}}+\frac {1034700961871398795883422607736832}{
690449250921115359375}{a}^{\frac{85}{2}}+\frac {
273184535623232651198367531008}{76716583435679484375}{a}^{\frac{99}{2
}}\\
&+\frac {1553964642453058067093\,}{147295840196504610}{a}^{\frac{19}{2}}+
\frac {619425670719478584198427574272}{2071347752763346078125}{a}^{
\frac{103}{2}}+\frac {48194100\,}{20205190699109}{a}^{\frac{3}{2}}+\frac
{188980047563720753152}{230149750307038453125}{a}^{\frac{123}{2}}
\end{align*}
\begin{align*}
&+\frac {26131912218100545548290921534062592}{2071347752763346078125}{a}
^{\frac{75}{2}}+\frac {10013252856065\,}{161641525592872}{a}^{\frac{9}{2}}
+\frac {10900296044337598744224540972285952}{690449250921115359375}{
a}^{\frac{73}{2}}+\frac {2413929400270585856}{
2071347752763346078125}{a}^{\frac{127}{2}}\\
&+\frac {
28779142733809718879\,}{21821605955037720}{a}^{\frac{17}{2}}+\frac {
155373628675526551\,}{1091080297751886}{a}^{\frac{15}{2}}+\frac {
157348903016613260931919368945664}{2071347752763346078125}{a}^{\frac{
93}{2}}+\frac {337787687394450811953879672544}{
414269550552669215625}{a}^{\frac{33}{2}}\\
&+\frac {
2471628909529830533671507984384}{230149750307038453125}{a}^{\frac{97}
{2}}+\frac {6847699507875135654425754861568}{230149750307038453125}
{a}^{\frac{95}{2}}+\frac {475567513562250952490172385869824}{
690449250921115359375}{a}^{\frac{47}{2}}\\
&+\frac {92009034576508033399748284256}{414269550552669215625}{a}^{\frac{31}{2
}}+\frac {689786473890065267897435807332352}{2071347752763346078125
}{a}^{\frac{45}{2}}+\frac {8938073504060482256896}{
690449250921115359375}{a}^{\frac{121}{2}}+\frac {
4562896572828267755381039264}{82853910110533843125}{a}^{\frac{29}{2}}\\
&+\frac {1133146212588316349305061267968}{414269550552669215625}{a}^
{\frac{35}{2}}+\frac {5808191821319646682524666280832}{
690449250921115359375}{a}^{\frac{37}{2}}+\frac {
16423961522469777653737235797504}{690449250921115359375}{a}^{\frac{39
}{2}}\\
&+\frac {8562452985086400375772161645568}{138089850184223071875
}{a}^{\frac{41}{2}}+\frac {309336635862909739523082892343296}{
2071347752763346078125}{a}^{\frac{43}{2}} \Big)  \Big) \sqrt {4
\,{a}^{2}+8\,a+3}\sqrt {2\,{a}^{3}+6\,{a}^{2}+5\,a+1}\\
&+\frac {
57218259482979\,\sqrt {{a}^{2}+3\,a+2}}{161641525592872}\Big( -
\frac {3959745050691481422249909551104}{1955255210769923015625}{a}^{
\frac{111}{2}}-\frac {15760906094621939392064468588102156288}{
17597296896929307140625}{a}^{\frac{81}{2}}\\
&-\frac {
14996976104471375780359168\,}{2234577383737054875}{a}^{\frac{23}{2}}
-{\frac {25608953002753860993187259625766912}{345045037194692296875}{a}^{\frac{91}{2}}}-\frac {942335379800076897701755446547841024}{
17597296896929307140625}{a}^{\frac{49}{2}}\\
&-\frac {3345860872631534458193488512956760064}{17597296896929307140625}{a}^{
\frac{53}{2}}-\frac {6150912360323068928\,}{
625681667446375365}{a}^{\frac{11}{2}}-\frac {204250842765454668501374933727838208}{
1955255210769923015625}{a}^{\frac{51}{2}}\\
&-\frac {
36597517711653299028760089579225088}{18621478197808790625}{a}^{\frac{
67}{2}}-\frac {51263227101184\,}{1986291007766271}{a}^{\frac{7}{2}}-\frac {10079904936945787481542885376}{703891875877172285625}{a}^{
\frac{117}{2}}-\frac {3943755456512\,}{4634679018121299}{a}^{\frac{5}{2}}\\
&-\frac {21317937989095636110332271274514972672}{
17597296896929307140625}{a}^{\frac{79}{2}}-\frac {
1809061718800521688524524975766372352}{3519459379385861428125}{a}^{
\frac{57}{2}}-\frac {191183682348483682950317651001344}{
5865765632309769046875}{a}^{\frac{107}{2}}\\
&-\frac {
9223372036854775808}{502779911340837346875}{a}^{\frac{133}{2}}-\frac {11084187345290226827264}{17597296896929307140625}{a}^{\frac{
131}{2}}-\frac {2001322618369935219468435126747136}{
17597296896929307140625}{a}^{\frac{105}{2}}-\frac {
4611686018427387904}{17597296896929307140625}{a}^{\frac{135}{2}}\\
&-\frac {150016474705073566888193951793152}{17597296896929307140625}{a}^
{\frac{109}{2}}-\frac {10927278465039712554730383069333684224}{
17597296896929307140625}{a}^{\frac{83}{2}}-\frac {
139469598724598835406482294784\,}{
3519459379385861428125}{a}^{\frac{25}{2}}\\
&-\frac {7645708832142297106974351818752}{
17597296896929307140625}{a}^{\frac{113}{2}}-\frac {97402880\,
}{514964335346811}\sqrt {a}-\frac {4207018162250559902567855949277560832}{1955255210769923015625}{a}^{
\frac{71}{2}}\\
&-\frac {2457243859062027311769312072262746112}{
17597296896929307140625}{a}^{\frac{89}{2}}-\frac {
79243632128095464775171072\,}{78210208430796920625}{a}^{\frac{21}{2}}
-\frac
{7468967431112083360540549713313660928}{3519459379385861428125}{a}^{
\frac{69}{2}}\\
&-\frac {4487233785670888126777149448788115456}{
5865765632309769046875}{a}^{{\frac{59}{2}}}-\frac {
4318773285079773133420803847207518208}{17597296896929307140625}{a}^{
\frac{87}{2}}-\frac {85831304802449306624\,}{
625681667446375365}{a}^{\frac{13}{2}}\\
&-\frac {248845424632837244452864}{
17597296896929307140625}{a}^{\frac{129}{2}}-\frac {
490483285454832443764416446464}{5865765632309769046875}{a}^{\frac{115
}{2}}-\frac {2725252497481863569462408894444732416}{
1955255210769923015625}{a}^{\frac{63}{2}}
\end{align*}
\begin{align*}
&-\frac {
30065865844097610871713917733775081472}{17597296896929307140625}{a}^{
\frac{65}{2}}-\frac {18826633114158762188746763779702784}{
17597296896929307140625}{a}^{\frac{101}{2}}-\frac {
750703318640357538512313396024573952}{703891875877172285625}{a}^{
\frac{61}{2}}\ \ \ \ \ \ \ \ \ \ \ \ \ \\
&-\frac {9205718471251038375510016}{4259815274008546875
}{a}^{\frac{119}{2}}-\frac {739198567122030975083750008832}{
3519459379385861428125}{a}^{\frac{27}{2}}-\frac {
3864969250656507885563951104063963136}{2513899556704186734375}{a}^{
\frac{77}{2}}\\
&-\frac {17731966463807316241154048}{
5865765632309769046875}{a}^{\frac{125}{2}}-\frac {
379259119929318406253345920076742656}{1173153126461953809375}{a}^{
\frac{55}{2}}-\frac {262946330404653410499516943216672768}{
651751736923307671875}{a}^{\frac{85}{2}}\\
&-\frac {
5678000856291327152773220572069888}{1955255210769923015625}{a}^{\frac
{99}{2}}-\frac {181320623106404849865728\,}{
1340746430242232925}{a}^{\frac{19}{2}}-\frac {426677876684484831021495759142912}{
1173153126461953809375}{a}^{\frac{103}{2}}-\frac {28104704000\,}{1544893006040433}{a}
^{\frac{3}{2}}\\
&-\frac {562528119726989076368195584}{
17597296896929307140625}{a}^{\frac{123}{2}}-\frac {
32229567063231094902620002797217644544}{17597296896929307140625}{a}^{
\frac{75}{2}}-\frac {119092461775287296\,}{
208560555815458455}{a}^{\frac{9}{2}}\\
&-\frac {288391572542969111223779762471174144}{
140778375175434457125}{a}^{\frac{73}{2}}-\frac {
4108474116306527376637952}{17597296896929307140625}{a}^{\frac{127}{2}}-\frac {148287862423081099752448\,}{9385225011695630475}{a}^{\frac{17}{2}}-\frac {998515100517383793664\,}{625681667446375365}{a}^{\frac{15}{2}}\\
&-\frac {647564356745536265492616008972632064}{17597296896929307140625}{
a}^{\frac{93}{2}}-\frac {20263121451488386626313810960384}{
1173153126461953809375}{a}^{\frac{33}{2}}-\frac {
18335619345818081625335306829955072}{2513899556704186734375}{a}^{
\frac{97}{2}}\\
&-\frac {33224321155528120873947010546270208}{
1955255210769923015625}{a}^{\frac{95}{2}}-\frac {
50012953532118990589388345061343232}{1955255210769923015625}{a}^{
\frac{47}{2}}-\frac {76801286494265762815349836894208}{
17597296896929307140625}{a}^{\frac{31}{2}}\\
&-\frac {
28577084740550244804070259304103936}{2513899556704186734375}{a}^{
\frac{45}{2}}-\frac {1665271960166937001405186048}{
5865765632309769046875}{a}^{\frac{121}{2}}-\frac {
17686550485028960662660396085248}{17597296896929307140625}{a}^{\frac{
29}{2}}\\
&-\frac {157210824206060190420577527455744}{
2513899556704186734375}{a}^{\frac{35}{2}}-\frac {
3657352634887339818533993208233984}{17597296896929307140625}{a}^{
\frac{37}{2}}-\frac {11190188928765772871374521626509312}{
17597296896929307140625}{a}^{\frac{39}{2}}\\
&-\frac {
31603041308845672256216177548427264}{17597296896929307140625}{a}^{
\frac{41}{2}}-\frac {27525320467362351253352630498000896}{
5865765632309769046875}{a}^{\frac{43}{2}}+ \Big( \frac {3856853439572223595676434432}{28155675035086891425}{a}^{\frac{111}{2}
}\\
&+\frac {505644835314940565206859873710505984}{
3519459379385861428125}{a}^{\frac{81}{2}}
+\frac {254236170116778534972492896\,}{100555982268167469375}{a}^{\frac{23}{2}}
+\frac {168057667698027950846882856801665024}{17597296896929307140625}{
a}^{\frac{91}{2}}\\
&+\frac {256656324561179185492498176722763776}{
17597296896929307140625}{a}^{\frac{49}{2}}
+\frac {860709635664547580550395494473662464}{17597296896929307140625}{a}^{\frac{53}{2}}+\frac {524749656146004022\,}{
125136333489275073}{a}^{\frac{11}{2}}\\
&+\frac {486704443185031035629619092605468672}{
17597296896929307140625}{a}^{\frac{51}{2}}+\frac {
7157867597894452412314131513745604608}{17597296896929307140625}{a}^{
\frac{67}{2}}+\frac {1430979518210668\,}{
125136333489275073}{a}^{\frac{7}{2}}\\
&+\frac {12616167506084273359918465024}{
17597296896929307140625}{a}^{\frac{117}{2}}+\frac {48054333036244
\,}{125136333489275073}{a}^{\frac{5}{2}}+\frac {
3557590633420414924309834165808791552}{17597296896929307140625}{a}^{
\frac{79}{2}}\\
&+\frac {2193113496346095341596709288106917888}{
17597296896929307140625}{a}^{\frac{57}{2}}+\frac {
45638077284505899429815902011392}{17597296896929307140625}{a}^{\frac{
107}{2}}+\frac {576460752303423488}{5865765632309769046875}{a}^{
\frac{133}{2}}\\
&+\frac {13258597302978740224}{1955255210769923015625}
{a}^{\frac{131}{2}}+\frac {171302523770744917562173941612544}{
17597296896929307140625}{a}^{\frac{105}{2}}
+\frac {
63080186320646332979383107584}{100555982268167469375}{a}^{\frac{109}{
2}}
\end{align*}
\begin{align*}
&+\frac {1682653084971115530357463749585010688}{
17597296896929307140625}{a}^{\frac{83}{2}}+\frac {
5709643926655248322305436456\,}{
391051042153984603125}{a}^{\frac{25}{2}}+\frac {94404166521309341766958186496}{
3519459379385861428125}{a}^{\frac{113}{2}}+\frac {11066496440\,
}{125136333489275073}\sqrt {a} \\
&+\frac {
1464047181433576403152239244957712384}{3519459379385861428125}{a}^{
\frac{71}{2}}+\frac {66321846595640169372493613537165312}{
3519459379385861428125}{a}^{\frac{89}{2}}
+\frac {
274591916710465739291097368\,}{703891875877172285625}{a}^{\frac{21}{2}}\\
&+\frac {498211804878544472902156152834359296}{1173153126461953809375}{a
}^{\frac{69}{2}}+\frac {452297307122284285137869158637633536}{
2513899556704186734375}{a}^{\frac{59}{2}}+\frac {
610071508491700942212685987632381952}{17597296896929307140625}{a}^{
\frac{87}{2}}\\
&+\frac {107700940495905797464\,}{
1877045002339126095}{a}^{\frac{13}{2}}+\frac {448630580480139329536}{
1955255210769923015625}{a}^{\frac{129}{2}}+\frac {
82221153144808505025518108672}{17597296896929307140625}{a}^{\frac{115
}{2}}+\frac {216735530181883821776942416806281216}{
703891875877172285625}{a}^{\frac{63}{2}}\\
&+\frac {
6431058950839603704065542039279763456}{17597296896929307140625}{a}^{
\frac{65}{2}}+\frac {14722030901750445370191634235392}{
140778375175434457125}{a}^{\frac{101}{2}}+\frac {
4279195750838027557257622351081111552}{17597296896929307140625}{a}^{
\frac{61}{2}}\\
&+\frac {561549554802732701247864832}{
5865765632309769046875}{a}^{\frac{119}{2}}+\frac {
266252012951161517031047663648}{3519459379385861428125}{a}^{\frac{27}
{2}}+\frac {670206844382548563281807000006033408}{
2513899556704186734375}{a}^{\frac{77}{2}}\\
&+\frac {
96821375179363497017344}{1173153126461953809375}{a}^{\frac{125}{2}}+\frac {1421143316347261237918846297961529344}{
17597296896929307140625}{a}^{\frac{55}{2}}+\frac {
1047923558623414464855177689273729024}{17597296896929307140625}{a}^{
\frac{85}{2}}\\
&+\frac {5308107985776005348348813474004992}{
17597296896929307140625}{a}^{\frac{99}{2}}+\frac {499226885818141578034408\,}{9385225011695630475}{a}^{\frac{19}{2}}+\frac {
117320304669188654505847949361152}{3519459379385861428125}{a}^{\frac{
103}{2}}+\frac {149419528880\,}{17876619069896439}{a}^{\frac{3}{2}}\\
&+\frac {6172181559929007857729536}{5865765632309769046875}{a}^{\frac{
123}{2}}+\frac {5800840586405778754014625454366916608}{
17597296896929307140625}{a}^{\frac{75}{2}}+\frac {310715729707828531\,}{1251363334892750730}{a}^{\frac{9}{2}}\\
&+\frac {
6727680393648428724029824632432361472}{17597296896929307140625}{a}^{
\frac{73}{2}}+\frac {1751143650309724700672}{345045037194692296875}
{a}^{\frac{127}{2}}+\frac {59558480142793463051198\,}{
9385225011695630475}{a}^{\frac{17}{2}}+\frac {1228043105038368077848\,}{
1877045002339126095}{a}^{\frac{15}{2}}\\
&+\frac {4663621192669418636340046543716352}{
1035135111584076890625}{a}^{\frac{93}{2}}+\frac {102010467004177131825393305577824}{17597296896929307140625}{a}^{\frac{33}{2}}+\frac {14122142508448885835560281098420224}{
17597296896929307140625}{a}^{\frac{97}{2}}\\
&+\frac {
34750460967315879383955249125392384}{17597296896929307140625}{a}^{
\frac{95}{2}}+\frac {126049417298227322675541899018190848}{
17597296896929307140625}{a}^{\frac{47}{2}}+\frac {
754347784262005906254560107648}{502779911340837346875}{a}^{\frac{31}{
2}}\\
&+\frac {57568734157466225535551729214889984}{
17597296896929307140625}{a}^{\frac{45}{2}}+\frac {
192753061829403353024561152}{17597296896929307140625}{a}^{\frac{121}{
2}}+\frac {889290956943758488991043563584}{2513899556704186734375}{
a}^{\frac{29}{2}}\\
&+\frac {360417026084504930046151467040384}{
17597296896929307140625}{a}^{\frac{35}{2}}
+\frac {1168339835852444918312117399606656}{17597296896929307140625}{a}^{
\frac{37}{2}}+\frac {697019967233935438463721016881664}
{3519459379385861428125}{a}^{\frac{39}{2}}\\
&+\frac {9591312545086985893399230866776064}{17597296896929307140625}{a}^{
\frac{41}{2}}+\frac {3487143841381149356981436617279488}{
2513899556704186734375}{a}^{\frac{43}{2}}
+\sqrt {2} \Big( \frac {
67449410081962263479905681408}{17597296896929307140625}{a}^{\frac{111
}{2}}\\
&+\frac {180831004347926501696060015924740096}{
17597296896929307140625}{a}^{\frac{81}{2}}+\frac {
23410403075249408743987396\,}{46926125058478152375}{a}^{\frac{23}{2}}+\frac
{9459958266543668969784046621032448}{17597296896929307140625}{a}^{
\frac{91}{2}}\\
&+\frac {317581118444994052606792859189248}{167593303780279115625}{a}^{\frac{49}{2}}+\frac {
104524968905573033166448428793987072}{17597296896929307140625}{a}^{
\frac{53}{2}}+{a}^{\frac{11}{2}}+\frac {20381609621776392733028822296887296
}{5865765632309769046875}{a}^{\frac{51}{2}}
\end{align*}
\begin{align*}
&+\frac {678015730602557946944527788941508608}{17597296896929307140625}{a}^{
\frac{67}{2}}
+\frac {166426821760\,}{57218259482979}{a}^{\frac{7}{2}}+\frac {252639511736572269886963712}{17597296896929307140625}{a}^{
\frac{117}{2}}+\frac {274997800\,}{2724679022999}{a}^{\frac{5}{2}}\\
&+\frac
{29521190093559390612050616302698496}{1955255210769923015625}{a}^{
\frac{79}{2}}+\frac {248611231569582440268647618351988736}{
17597296896929307140625}{a}^{\frac{57}{2}}+\frac {
508487987185544053063935852544}{5865765632309769046875}{a}^{\frac{107
}{2}}\\
&+\frac {288230376151711744}{17597296896929307140625}{a}^{\frac{131}{2}}+\frac {6199232645792463860283937390592}{
17597296896929307140625}{a}^{\frac{105}{2}}+\frac {
112926599004650467795960594432}{5865765632309769046875}{a}^{\frac{109
}{2}}\\
&+\frac {12789665091230751542739283839287296}{
1955255210769923015625}{a}^{\frac{83}{2}}+\frac {
131012309060200020704552528\,}{
46926125058478152375}{a}^{\frac{25}{2}}+\frac {477703665190099164545941504}{
703891875877172285625}{a}^{\frac{113}{2}}+\frac {474000\,
}{19072753160993}\sqrt {a}\\
&+\frac {642750707490399683403584879693135872}{
17597296896929307140625}{a}^{\frac{71}{2}}+\frac {
2183901341548830505187225991380992}{1955255210769923015625}{a}^{\frac
{89}{2}}+\frac {106495626089303645377952\,}{
1340746430242232925}{a}^{\frac{21}{2}}\\
&+\frac {97397031485515180768007648016596992}{
2513899556704186734375}{a}^{\frac{69}{2}}+\frac {
115521665728432067837768707230138368}{5865765632309769046875}{a}^{
\frac{59}{2}}+\frac {1407326482911818769656376713019392}{
651751736923307671875}{a}^{\frac{87}{2}}\\
&+\frac {2925070859520590\,
}{220699000862919}{a}^{\frac{13}{2}}+\frac {1152921504606846976}{
1035135111584076890625}{a}^{\frac{129}{2}}+\frac {
1862224592364205569155792896}{17597296896929307140625}{a}^{\frac{115}
{2}}+\frac {552297192064845288064983731070304256}{
17597296896929307140625}{a}^{\frac{63}{2}}\\
&+\frac {632097693638939424895061771134435328}{17597296896929307140625}{a}^{
\frac{65}{2}}+\frac {10988281102713155837908282769408}{
2513899556704186734375}{a}^{\frac{101}{2}}+\frac {
150703649257925976910341528069603328}{5865765632309769046875}{a}^{
\frac{61}{2}}\\
&+\frac {29361917076081287880507392}{
17597296896929307140625}{a}^{\frac{119}{2}}+\frac {9867008269096335580479132784}{703891875877172285625}{a}^{\frac{27}{2}
}+\frac {10439864615546125378626527797706752}{502779911340837346875
}{a}^{\frac{77}{2}}+\frac {2848148461943139598336}{
3519459379385861428125}{a}^{\frac{125}{2}}\\
&+\frac {166773392338463463128883180556255232}{17597296896929307140625}{a}^{
\frac{55}{2}}+\frac {68459551600788324932650127291580416}{
17597296896929307140625}{a}^{\frac{85}{2}}+\frac {
1548791772423307523419581448192}{115015012398230765625}{a}^{\frac{99}
{2}}\\
&+\frac {1398579084725018935621\,}{125136333489275073}{a}^{\frac{19}{2}}+\frac {22865513390129666880098013282304}{17597296896929307140625}{a
}^{\frac{103}{2}}+\frac {6187600\,}{2724679022999}{a}^{\frac{3}{2}}+
\frac {76040109174013156130816}{5865765632309769046875}{a}^{\frac{123
}{2}}\\
&+\frac {470662402808906734247716355607363584}{
17597296896929307140625}{a}^{\frac{75}{2}}+\frac {10495594012745\,
}{171654778448937}{a}^{\frac{9}{2}}+\frac {
113619498842656937202915220123549696}{3519459379385861428125}{a}^{
\frac{73}{2}}+\frac {217830106776656150528}{5865765632309769046875}
{a}^{\frac{127}{2}}\\
&+\frac {31886615307320249899\,}{
23173395090606495}{a}^{\frac{17}{2}}+\frac {678556681561031050\,}{
4634679018121299}{a}^{\frac{15}{2}}+\frac {1409543123480878730024798772002816}{
5865765632309769046875}{a}^{\frac{93}{2}}+\frac {
45844983676902276837720315392}{46926125058478152375}{a}^{\frac{33}{2}
}\\
&+\frac {74554627339239729504778102767616}{1955255210769923015625}{
a}^{\frac{97}{2}}+\frac {1751825588699915854574009301598208}{
17597296896929307140625}{a}^{\frac{95}{2}}+\frac {
16931912116614545002943623719534592}{17597296896929307140625}{a}^{
\frac{47}{2}}\\
&+\frac {306196046508398325804409281536}{
1173153126461953809375}{a}^{\frac{31}{2}}+\frac {
2664328292674600019788049865801728}{5865765632309769046875}{a}^{\frac
{45}{2}}+\frac {409209960452216055660544}{2513899556704186734375}{a
}^{\frac{121}{2}}+\frac {2629837347008374870093968448}{
41405404463363075625}{a}^{\frac{29}{2}}\\
&+\frac {
11767471199862137290071527524096}{3519459379385861428125}{a}^{\frac{
35}{2}}+\frac {26388921745029167764684920208384}{
2513899556704186734375}{a}^{\frac{37}{2}}+\frac {
59285742250105131563365072208896}{1955255210769923015625}{a}^{\frac{
39}{2}}\\
&+\frac {22566181503962582679084969152512}{
279322172967131859375}{a}^{\frac{41}{2}}+\frac {
3502187445113645961277322536214528}{17597296896929307140625}{a}^{
\frac{43}{2}} \Big) \sqrt {2\,{a}^{2}+4\,a+1} \Big) \sqrt {4\,{a
}^{2}+8\,a+3} \Big)  \Big)
\end{align*}
\end{landscape}
\end{appendix}

\bibliographystyle{abbrv}
\bibliography{Lit_CFMS}

\end{document}